\documentclass{article}
\usepackage{amsmath}
\usepackage{amssymb}
\usepackage{amsthm}
\usepackage{cite}
\usepackage[arrow,matrix,curve]{xy}
\begin{document}
\def\e#1\e{\begin{equation}#1\end{equation}}
\def\ea#1\ea{\begin{align}#1\end{align}}
\def\eq#1{{\rm(\ref{#1})}}
\theoremstyle{plain}
\newtheorem{thm}{Theorem}[section]
\newtheorem{lem}[thm]{Lemma}
\newtheorem{prop}[thm]{Proposition}
\newtheorem{cor}[thm]{Corollary}
\theoremstyle{definition}
\newtheorem{dfn}[thm]{Definition}
\newtheorem{ex}[thm]{Example}
\newtheorem{ass}[thm]{Assumption}
\def\Re{\mathop{\rm Re}}
\def\Im{\mathop{\rm Im}}
\def\Ker{\mathop{\rm Ker}}
\def\Hom{\mathop{\rm Hom}}
\def\Obj{\mathop{\rm Obj\kern .1em}\nolimits}
\def\End{\mathop{\rm End}}
\def\Mor{\mathop{\rm Mor}\nolimits}
\def\id{\mathop{\rm id}\nolimits}
\def\Quot{\mathop{\rm Quot}\nolimits}
\def\GL{\mathop{\rm GL}\nolimits}
\def\Ext{\mathop{\rm Ext}\nolimits}
\def\Aut{\mathop{\rm Aut}}
\def\gro{({\rm groupoids})}
\def\excat{({\rm exactcat})}
\def\Iso{\mathop{\rm Iso}\nolimits}
\def\Spec{\mathop{\rm Spec}\nolimits}
\def\Sch{\mathop{\rm Sch}\nolimits}
\def\coh{\mathop{\rm coh}}
\def\qcoh{\mathop{\rm qcoh}}
\def\fObj{\mathop{\mathfrak{Obj}\kern .05em}\nolimits}
\def\fSurj{\mathop{\mathfrak{Surj}\kern .05em}\nolimits}
\def\fExact{\mathop{\mathfrak{Exact}\kern .05em}\nolimits}
\def\fb{{\mathfrak b}}
\def\fe{{\mathfrak e}}
\def\fm{{\mathfrak m}}
\def\ge{\geqslant}
\def\le{\leqslant}
\def\pr{{\mathop{\preceq}\nolimits}}
\def\npr{{\mathop{\npreceq}\nolimits}}
\def\tl{\trianglelefteq\nobreak}
\def\ntl{\ntrianglelefteq}
\def\ps{\precsim\nobreak}
\def\nps{\not\precsim}
\def\ls{\mathop{\lesssim\kern .05em}\nolimits}
\def\nls{\mathop{\not\lesssim\kern .05em}}
\def\bu{\bullet}
\def\bdim{{\mathbin{\bf dim}}}
\def\modA{\text{\rm mod-$A$}}
\def\modKQ{\text{\rm mod-$\K Q$}}
\def\modKQI{\text{\rm mod-$\K Q/I$}}
\def\nilKQ{\text{\rm nil-$\K Q$}}
\def\nilKQI{\text{\rm nil-$\K Q/I$}}
\def\K{{\mathbin{\mathbb K}}}
\def\KP{{\mathbin{\mathbb{KP}}}}
\def\N{{\mathbin{\mathbb N}}}
\def\Q{{\mathbin{\mathbb Q}}}
\def\Z{{\mathbin{\mathbb Z}}}
\def\A{{\mathbin{\mathcal A}}}
\def\B{{\mathbin{\mathcal B}}}
\def\E{{\mathbin{\mathcal E}}}
\def\F{{\mathbin{\mathcal F}}}
\def\G{{\mathbin{\mathcal G}}}
\def\H{{\mathbin{\mathcal H}}}
\def\M{{\mathcal M}}
\def\O{{\mathbin{\mathcal O}}}
\def\fE{{\mathbin{\mathfrak E}}}
\def\fF{{\mathbin{\mathfrak F}}}
\def\fG{{\mathbin{\mathfrak G}}}
\def\fH{{\mathbin{\mathfrak H}}}
\def\fM{{\mathfrak M}}
\def\sIp{{\smash{\sst(I,\pr)}}}
\def\sJp{{\smash{\sst(J,\pr)}}}
\def\sKp{{\smash{\sst(K,\pr)}}}
\def\sIt{{\smash{\sst(I,\tl)}}}
\def\sJt{{\smash{\sst(J,\tl)}}}
\def\sKt{{\smash{\sst(K,\tl)}}}
\def\sLt{{\smash{\sst(L,\tl)}}}
\def\sIl{{\smash{\sst(I,\ls)}}}
\def\sJl{{\smash{\sst(J,\ls)}}}
\def\sKl{{\smash{\sst(K,\ls)}}}
\def\sIb{{\smash{\sst(I,\bu)}}}
\def\al{\alpha}
\def\be{\beta}
\def\ga{\gamma}
\def\de{\delta}
\def\ep{\epsilon}
\def\io{\iota}
\def\ka{\kappa}
\def\la{\lambda}
\def\ze{\zeta}
\def\th{\theta}
\def\si{\sigma}
\def\up{\upsilon}
\def\De{\Delta}
\def\Th{\Theta}
\def\ts{\textstyle}
\def\sst{\scriptscriptstyle}
\def\sm{\setminus}
\def\op{\oplus}
\def\ot{\otimes}
\def\bigop{\bigoplus}
\def\iy{\infty}
\def\ra{\rightarrow}
\def\longra{\longrightarrow}
\def\dashra{\dashrightarrow}
\def\t{\times}
\def\bs{\boldsymbol}
\def\ci{\circ}
\def\ti{\tilde}
\def\ov{\overline}
\def\md#1{\vert #1 \vert}
\def\bmd#1{\big\vert #1 \big\vert}
\def\ab{\allowbreak}
\title{Configurations in abelian categories. I. \\
Basic properties and moduli stacks}
\author{Dominic Joyce}
\date{}
\maketitle

\begin{abstract}
This is the first in a series of papers on {\it configurations\/}
in an abelian category $\A$. Given a finite partially ordered
set $(I,\pr)$, an $(I,\pr)$-{\it configuration\/} $(\si,\io,\pi)$
is a finite collection of objects $\si(J)$ and morphisms $\io(J,K)$
or $\pi(J,K):\si(J)\ra\si(K)$ in $\A$ satisfying some axioms, where
$J,K$ are subsets of $I$. Configurations describe how an object
$X$ in $\A$ decomposes into subobjects, and are useful for
studying {\it stability conditions\/} on~$\A$.

We define and motivate the idea of configurations, and explain
some natural operations upon them --- subconfigurations, quotient
configurations, substitution, refinements and improvements. Then
we study {\it moduli spaces} of $(I,\pr)$-configurations in $\A$,
and natural morphisms between them, using the theory of
{\it Artin stacks}. We prove well-behaved moduli stacks exist
when $\A$ is the abelian category of coherent sheaves on a
projective scheme $P$, or of representations of a quiver~$Q$.

In the sequels, given a stability condition $(\tau,T,\le)$
on $\A$, we will show the moduli spaces of $\tau$-(semi)stable
objects or configurations are constructible subsets in the
moduli stacks of all objects or configurations. We associate
infinite-dimensional algebras of constructible functions to
a quiver $Q$ using the method of Ringel--Hall algebras, and
define systems of invariants of $P$ that `count'
$\tau$-(semi)stable coherent sheaves on $P$ and satisfy
interesting identities.
\end{abstract}

\section{Introduction}
\label{aa1}

This is the first of a series of papers \cite{Joyc1,Joyc2,Joyc3}
developing the concept of {\it configuration} in an abelian category.
Given an abelian category $\A$ and a finite partially ordered
set (poset) $(I,\pr)$, we define an $(I,\pr)$-{\it configuration}
$(\si,\io,\pi)$ in $\A$ to be a collection of objects $\si(J)$ and
morphisms $\io(J,K)$ or $\pi(J,K):\si(J)\ra\si(K)$ in $\A$ satisfying
certain axioms, where $J,K$ are subsets of $I$. Configurations are a
tool for describing {\it how an object\/ $X$ in $\A$ decomposes into
subobjects}. They are especially useful for studying {\it stability
conditions} on~$\A$.

This paper introduces configurations, studies their basic properties,
and develops the theory of {\it moduli stacks of configurations}. We
begin in \S\ref{aa2} with background material on abelian categories
and Artin stacks. Section \ref{aa3} refines the Jordan--H\"older
Theorem for abelian categories in the case when the simple factors
$S_1,\ldots,S_n$ of $X\in\A$ are nonisomorphic. We find that the set
of all {\it subobjects} of $X$ may be classified using a {\it
partial order} $\pr$ on $I=\{1,\ldots,n\}$, the indexing
set for the simple factors of $X$. We also classify quotient
objects and composition series for $X$ using~$(I,\pr)$.

Motivated by this, \S\ref{aa4} defines the notion of
$(I,\pr)$-{\it configuration} $(\si,\io,\pi)$ in $\A$,
and proves that it captures the properties of the set
of all subobjects of $X\in\A$ when $X$ has nonisomorphic
simple factors $\{S^i:i\in I\}$. Section \ref{aa5}
considers some elementary operations on configurations.
Given an $(I,\pr)$-configuration we can make {\it sub-}
and {\it quotient\/ $(K,\tl)$-configurations}, where
$(K,\tl)$ comes from $(I,\pr)$ with $K\subseteq I$ or
using a surjective $\phi:I\ra K$. We also construct new
configurations by {\it substituting} one configuration
into another.

Let $\tl,\pr$ be partial orders on $I$, with $i\pr j$ implies
$i\tl j$. Then each $(I,\pr)$-configuration $(\si,\io,\pi)$
has a quotient $(I,\tl)$-configuration $(\ti\si,\ti\io,\ti\pi)$.
Call $(\si,\io,\pi)$ an $(I,\pr)$-{\it improvement\/} of
$(\ti\si,\ti\io,\ti\pi)$. Call $(\ti\si,\ti\io,\ti\pi)$
{\it best\/} if it has no strict improvements. Section \ref{aa6}
shows that improvements can be divided into a sequence of
{\it steps}, classifies {\it one step improvements}, and gives
a criterion for best configurations in terms of whether short
exact sequences split.

Fix an algebraically closed field $\K$, and a $\K$-{\it linear}
abelian category $\A$. To form moduli spaces of configurations
in $\A$ we need some extra data, on algebraic families of
objects and morphisms in $\A$ parametrized by a base $\K$-scheme
$U$. We encode this in a {\it stack in exact categories} $\fF_\A:
\Sch_\K\ra\excat$, which must satisfy conditions given in
Assumptions \ref{aa7ass} and \ref{aa8ass} below.

Section \ref{aa7} defines {\it moduli stacks of objects}
$\fObj_\A$ and $(I,\pr)$-{\it configurations} $\fM(I,\pr)_\A$
in $\A$, and {\it substacks} $\fObj_\A^\al,\fM(I,\pr,\ka)_\A$
of objects and configurations with prescribed classes in $K(\A)$.
There are many natural 1-{\it morphisms} between these stacks.
Section \ref{aa8} shows these are {\it algebraic} ({\it Artin})
$\K$-stacks, {\it locally of finite type}, and some of the
1-morphisms are {\it representable}, or of {\it finite type}.

We finish with some examples. Section \ref{aa9} takes $\A$ to
be the abelian category $\coh(P)$ of {\it coherent sheaves} on
a projective $\K$-scheme $P$, and \S\ref{aa10} considers the
abelian category $\modKQ$ of {\it representations of a quiver}
$Q$ and some variants $\nilKQ,\modKQI,\nilKQI$, and the abelian
category $\modA$ of representations of a {\it finite-dimensional\/
$\K$-algebra} $A$. We define the data $\A,K(\A),\fF_\A$ and prove
it satisfies Assumptions \ref{aa7ass} and \ref{aa8ass} in each example.

The second paper \cite{Joyc1} defines and studies
infinite-dimensional algebras of constructible functions on
$\M(I,\pr)_\A$, motivated by the idea of {\it Ringel--Hall
algebras}. The sequels \cite{Joyc2,Joyc3} concern {\it stability
conditions} $(\tau,T,\le)$ on $\A$, such as Gieseker stability on
$\coh(P)$, or slope stability on~$\modKQ$.

We shall regard the set $\Obj_{\rm ss}^\al(\tau)$ of
$\tau$-{\it semistable} objects in $\A$ with class $\al$ in
$K(\A)$ not as a moduli scheme under S-equivalence, but
as a {\it constructible subset\/} in the stack $\fObj_\A^\al$.
One of our goals is to understand the relationship between
$\Obj_{\rm ss}^\al(\tau)$ and $\Obj_{\rm ss}^\al(\ti\tau)$
for {\it two different\/} stability conditions $(\tau,T,\le),
(\ti\tau,\ti T,\le)$. Our key idea is that this is best done
using the moduli stacks~$\fM(I,\pr,\ka)_\A$.

Write $\M_{\rm ss}(I,\pr,\ka,\tau)_\A$ for the subset
of points $[(\si,\io,\pi)]$ in $\M(I,\pr,\ka)_\A$ with
$\si(\{i\})$ $\tau$-semistable for all $i\in I$.
We shall express $\Obj_{\rm ss}^\al(\ti\tau)$ and
$\M_{\rm ss}(I,\pr,\ka,\ti\tau)_\A$ in terms of
projections of $\M_{\rm ss}(K,\tl,\mu,\tau)_\A$ for other
finite posets $(K,\tl)$. We will then define systems of
invariants of $\A,(\tau,T,\le)$ by taking weighted Euler
characteristics of $\M_{\rm ss}(I,\pr,\ka,\tau)_\A$, and
determine identities the invariants satisfy, and their
transformation laws as $(\tau,T,\le)$ changes.
\medskip

\noindent{\it Acknowledgements.} I would like to thank
Tom Bridgeland for many inspiring conversations and for
being interested, Frances Kirwan and Burt Totaro for
help with moduli spaces and stacks, and Bernd Siebert
for explaining Quot-schemes over a base. I also want
to thank Ian Grojnowski, Alastair King, Andrew Kresch,
Paul Seidel, and Richard Thomas for useful conversations.
I was supported by an EPSRC Advanced Research Fellowship
whilst writing this paper.

\section{Background material}
\label{aa2}

We review {\it abelian categories} and {\it Artin stacks}.
Some useful references for \S\ref{aa21} and \S\ref{aa22} are
Popescu \cite{Pope} and Gelfand and Manin \cite[\S II.5--\S
II.6]{GeMa}, and for \S\ref{aa23} are G\'omez \cite{Gome},
Behrend et al.\ \cite{BEFF}, and Laumon and Moret-Bailly~\cite{LaMo}.

\subsection{Abelian and exact categories}
\label{aa21}

Here is the definition of abelian category, taken
from~\cite[\S II.5]{GeMa}.

\begin{dfn} A category $\A$ is called {\it abelian} if
\begin{itemize}
\setlength{\itemsep}{0pt}
\setlength{\parsep}{0pt}
\item[(i)] $\Hom(X,Y)$ is an abelian group for all
$X,Y\in\A$, and composition of morphisms is biadditive.
\item[(ii)] There exists a {\it zero object\/} $0\in\A$ such
that~$\Hom(0,0)=0$.
\item[(iii)] For any $X,Y\in\A$ there exists $Z\in\A$ and
morphisms $\io_X:X\ra Z$, $\io_Y:Y\ra Z$, $\pi_X:Z\ra X$,
$\pi_Y:Z\ra Y$ with $\pi_X\ci\io_X=\id_X$, $\pi_Y\ci\io_Y=
\id_Y$, $\io_X\ci\pi_X+\io_Y\ci\pi_Y=\id_Z$ and
$\pi_X\ci\io_Y=\pi_Y\ci\io_X=0$. We write $Z=X\op Y$,
the {\it direct sum} of $X$ and~$Y$.
\item[(iv)] For any morphism $f:X\ra Y$ there is a sequence
$\smash{K{\buildrel k\over\ra}X{\buildrel i\over\ra}I{\buildrel
j\over\ra}Y{\buildrel c\over\ra}C}$ in $\A$ such that $j\ci i=f$,
and $K$ is the kernel of $f$, and $C$ the cokernel of $f$, and
$I$ is both the cokernel of $k$ and the kernel of~$c$.
\end{itemize}
An abelian category $\A$ is called $\K$-{\it linear} over a field
$\K$ if $\Hom(X,Y)$ is a $\K$-vector space for all $X,Y\in\A$, and
composition maps are bilinear.
\label{aa2def1}
\end{dfn}

We will often use the following properties of abelian categories:
\begin{itemize}
\setlength{\itemsep}{0pt}
\setlength{\parsep}{0pt}
\item If $i\ci f=i\ci g$ and $i$ is injective, then $f=g$
($i$ is {\it left cancellable}).
\item If $f\ci\pi=g\ci\pi$ and $\pi$ is surjective, then $f=g$
($\pi$ is {\it right cancellable}).
\item If $f:X\ra Y$ is injective and surjective, then it is
an isomorphism.
\end{itemize}
In an abelian category $\A$ we can define {\it exact
sequences}~\cite[\S II.6]{GeMa}.

\begin{dfn} Let $\smash{X{\buildrel f\over\ra}Y{\buildrel g\over\ra}Z}$
be a sequence in $\A$ with $g\ci f=0$. Let $k:K\ra Y$ be the kernel of
$g$ and $c:Y\ra C$ the cokernel of $f$. Then there exist unique morphisms
$a:X\ra K$ and $b:C\ra Z$ with$ f=k\ci a$ and $g=b\ci c$. We say
$\smash{X{\buildrel f\over\ra}Y{\buildrel g\over\ra}Z}$ is {\it
exact at\/} $Y$ if $a$ is surjective, or equivalently if $b$ is
injective.

A short exact sequence $0\ra X\ra Y\ra Z\ra 0$ in $\A$ is called
{\it split\/} if there exists a compatible isomorphism $X\op Z\ra Y$.
The {\it Grothendieck group} $K_0(\A)$ of $\A$ is the abelian group
generated by $\Obj(\A)$, with a relation $[Y]=[X]+[Z]$ for
each short exact sequence $0\!\ra\!X\!\ra\!Y\!\ra\!Z\!\ra\!0$
in $\A$. Throughout the paper $K(\A)$ will mean {\it the quotient
of\/ $K_0(\A)$ by some fixed subgroup}.
\label{aa2def2}
\end{dfn}

{\it Exact categories} were introduced by Quillen \cite[\S 2]{Quil},
and are discussed in Gelfand and Manin~\cite[Ex.~IV.3.3, p.~275]{GeMa}.

\begin{dfn} Let $\skew8\hat\A$ be an abelian category, and $\A$ be
a full additive subcategory of $\skew8\hat\A$, which is closed under
extensions. Let $\E$ be the class of exact sequences $0\ra X\ra Y\ra
Z\ra 0$ in $\skew8\hat\A$ with $X,Y,Z\in\A$. Then the pair $(\A,\E)$
is called an {\it exact category}. Usually we refer to $\A$ as the
exact category, taking $\E$ to be implicitly given. Quillen
\cite[\S 2]{Quil} gives necessary and sufficient conditions
on $\A,\E$ for $\A$ to be embedded in an abelian category
$\skew8\hat\A$ in this way, and we take this to be the
{\it definition} of an exact category. An {\it exact functor}
$F:(\A,\E)\ra(\A',\E')$ of exact categories is a functor
$F:\A\ra\A'$ taking exact sequences $\E$ in $\A$ to exact
sequences $\E'$ in~$\A'$.
\label{aa2def3}
\end{dfn}

\subsection{Subobjects and the Jordan--H\"older Theorem}
\label{aa22}

{\it Subobjects} of objects in $\A$ are analogous to
subgroups of an abelian group.

\begin{dfn} Let $\A$ be an abelian category. Two injective
morphisms $i:S\ra X$, $i':S'\ra X$ in $\A$ are {\it equivalent\/}
if there exists an isomorphism $h:S\ra S'$ with $i=i'\ci h$. Then
$h$ is unique. A {\it subobject\/} of $X\in\A$ is an equivalence
class of injective morphisms $i:S\ra X$. Usually we refer to $S$
as the subobject, and write $S\subset X$ to mean $S$ is a subobject
of $X$. We write $0,X$ for the subobjects of $X$ which are
equivalence classes of $0\ra X$ and~$\id_X:X\ra X$.

Similarly, surjective morphisms $\pi:X\ra Q$, $\pi':X\ra Q'$ in
$\A$ are {\it equivalent\/} if there is an isomorphism $h:Q\ra Q'$
with $\pi'=h\ci\pi$. A {\it quotient object\/} of $X\in\A$ is an
equivalence class of surjective $\pi:X\ra Q$. If $S,T\subset X$
are represented by $i:S\ra X$ and $j:T\ra X$, we write $S\subset
T\subset X$ if there exists $a:S\ra T$ with $i=j\ci a$. Then $a$
fits into an exact sequence $0\ra S\smash{{\buildrel a\over\ra}
T{\buildrel b\over\ra}}F\ra 0$. We write $F=T/S$, and call $F$
a {\it factor} of~$X\in\A$.
\label{aa2def4}
\end{dfn}

We define operations $\cap,+$ on subobjects, following Popescu
\cite[\S 2.6]{Pope}. The notation comes from the intersection
and sum of {\it subgroups of abelian groups}.

\begin{dfn} Let $\A$ be an abelian category, let $X\in\A$, and
suppose injective maps $i:S\ra X$, $j:T\ra X$ define subobjects
$S,T$ of $X$. Apply Definition \ref{aa2def1}(iv) to
$f=i\ci\pi_S\!+\!j\ci\pi_T:S\op T\ra X$. This yields $U,V\in\A$
and morphisms $k:U\ra S\op T$, $l:S\op T\ra V$ and $e:V\ra X$
such that $i\ci\pi_S\!+\!j\ci\pi_T=e\ci l$, and $k$ is the kernel
of $i\ci\pi_S\!+\!j\ci\pi_T$, and $l$ is the cokernel of $k$, and
$e$ is the {\it image} (the kernel of the cokernel)
of~$i\ci\pi_S\!+\!j\ci\pi_T$.

Define $a:U\ra S$ by $a=k\ci\pi_S$, and $b:U\ra T$ by $b=-k\ci\pi_T$
and $c:S\ra V$ by $c=f\ci\io_S$, and $d:T\ra V$ by $d=f\ci\io_T$. Then
$k=\io_S\ci a-\io_T\ci b$, $l=c\ci\pi_S+d\ci\pi_T$, $i=e\ci c$ and
$j=e\ci d$. Now $0\ra U{\buildrel k\over\longra}S\op T{\buildrel
l\over\longra}V\ra 0$ is exact. So $i\ci a=j\ci b$, and
\e
\xymatrix{
0 \ar[r] & U \ar[rr]^{\io_S\ci a-\io_T\ci b\,\,\,\,}
&& S\op T \ar[rr]^{\,\,\,\,c\ci\pi_S+d\ci\pi_T}
&& V \ar[r] & 0
}
\quad\text{is exact.}
\label{aa2eq1}
\e

As $i,a$ are injective $i\ci a=j\ci b:U\ra X$ is too, and
defines a subobject $S\cap T$ of $X$. Also $e:V\ra X$ is
injective, and defines a subobject $S+T$ of $X$. Then
$S\cap T,S+T$ depend only on $S,T\subset X$, with inclusions
$S\!\cap\!T\!\subset\!S,T\!\subset\!S+T\!\subset\!X$. Popescu
\cite[Prop.~2.6.4, p.~39]{Pope} gives canonical isomorphisms
\e
S/(S\cap T)\cong(S+T)/T \quad\text{and}\quad
T/(S\cap T)\cong(S+T)/S.
\label{aa2eq2}
\e
These operations $\cap,+$ are {\it commutative} and {\it associative},
so we can form multiple sums and intersections. We write $\sum_{j\in J}
T_j$ for the multiple sum $+$ of a finite set of subobjects $T_j\subset
X$, in the obvious way.
\label{aa2def5}
\end{dfn}

\begin{dfn} We call an abelian category $\A$ {\it artinian} if all
descending chains of subobjects $\cdots\!\subset\!A_2\!\subset
\!A_1\!\subset\!X$ stabilize, that is, $A_{n+1}=A_n$ for $n\gg 0$.
We call $\A$ {\it noetherian} if all ascending chains of subobjects
$A_1\!\subset\!A_2\!\subset\!\cdots\!\subset\!X$ stabilize. We
call $\A$ of {\it finite length\/} if it is artinian and noetherian.

A nonzero object $X$ in $\A$ is called {\it simple} if
it has no nontrivial proper subobjects. Let $X\in\A$
and consider {\it filtrations} of subobjects
\e
0=A_0\subset A_1\subset\cdots\subset A_n=X.
\label{aa2eq3}
\e
We call \eq{aa2eq3} a {\it composition series} if the
factors $S_k\!=\!A_k/A_{k-1}$ are all {\it simple}.
\label{aa2def6}
\end{dfn}

Here is the {\it Jordan--H\"older Theorem} in an abelian
category,~\cite[Th.~2.1]{Sesh}.

\begin{thm} Let\/ $\A$ be an abelian category of finite length.
Then every filtration $0=A_0\subset A_1\subset\cdots\subset A_n=X$
without repetitions can be refined to a composition series for $X$.
Suppose $0=A_0\subset A_1\subset\cdots\subset A_m=X$ and\/ $0=B_0
\subset B_1\subset\cdots\subset B_n=X$ are two composition series
for $X\in\A$, with simple factors $S_k=A_k/A_{k-1}$ and\/
$T_k=B_k/B_{k-1}$. Then $m=n$, and for some permutation $\si$ of\/
$1,\ldots,n$ we have $S_k\cong T_{\si(k)}$ for~$k=1,\ldots,n$.
\label{aa2thm1}
\end{thm}

\subsection{Introduction to algebraic $\K$-stacks}
\label{aa23}

Fix an {\it algebraically closed field\/} $\K$ throughout.
There are four main classes of `spaces' over $\K$ used in
algebraic geometry, in increasing order of generality:
\begin{equation*}
\text{$\K$-varieties}\subset
\text{$\K$-schemes}\subset
\text{algebraic $\K$-spaces}\subset
\text{algebraic $\K$-stacks}.
\end{equation*}
{\it Algebraic stacks} (also known as Artin stacks) were
introduced by Artin, generalizing {\it Deligne--Mumford stacks}.
For an introduction see G\'omez \cite{Gome}, and for a thorough
treatment see Behrend et al.\ \cite{BEFF} or Laumon and
Moret-Bailly~\cite{LaMo}.

We write our definitions in the language of 2-{\it categories}
\cite[App.~B]{Gome}. A 2-category has {\it objects} $X,Y$,
1-{\it morphisms} $f,g:X\ra Y$ between objects, and
2-{\it morphisms} $\al:f\ra g$ between 1-morphisms. An example
to keep in mind is a 2-{\it category of categories}, where
{\it objects} are categories, 1-{\it morphisms} are functors,
and 2-{\it morphisms} are isomorphisms (natural transformations)
of functors.

As in G\'omez \cite[\S 2.1--\S 2.2]{Gome} there are two
different but equivalent ways of defining $\K$-stacks.
We shall work with the first \cite[Def.~2.10]{Gome}, even
though the second is more widely used, as it is more
convenient for our applications.

\begin{dfn} A {\it groupoid\/} is a category with all
morphisms isomorphisms. Let $\gro$ be the 2-category
whose {\it objects} are groupoids, 1-{\it morphisms}
functors of groupoids, and 2-{\it morphisms} natural
transformations of functors.

Let $\K$ be an algebraically closed field, and $\Sch_\K$
the {\it category of\/ $\K$-schemes}. We make $\Sch_\K$
into a 2-{\it category} by taking 1-morphisms to be
morphisms, and the only 2-morphisms to be identities
$\id_f$ for each 1-morphism $f$. To define $\K$-stacks
we need to choose a {\it Grothendieck topology} on
$\Sch_\K$, as in \cite[App.~A]{Gome}, and we choose
the {\it \'etale topology}.

A {\it prestack in groupoids\/} on $\Sch_\K$ is a {\it
contravariant\/ $2$-functor} $\fF:\Sch_\K\ra\gro$. As
in \cite[App.~B]{Gome}, this comprises the following
data, satisfying conditions we shall not give:
\begin{itemize}
\setlength{\itemsep}{0pt}
\setlength{\parsep}{0pt}
\item For each object $U$ in $\Sch_\K$, an object (groupoid)
$\fF(U)$ in~$\gro$.
\item For each 1-morphism $f:U\ra V$ in $\Sch_\K$, a 1-morphism
(functor) $\fF(f):\fF(V)\ra\fF(U)$ in~$\gro$.
\item For each 2-morphism $\al:f\ra f'$ in $\Sch_\K$, a
2-morphism $\fF(\al):\fF(f')\ra\fF(f)$. As the only 2-morphisms
in $\Sch_\K$ are $\id_f$, this data is {\it trivial}.
\item If $f:U\ra V$ and $g:V\ra W$ are 1-morphisms in $\Sch_\K$,
a 2-{\it isomorphism} $\ep_{g,f}:\fF(f)\ci\fF(g)\ra\fF(g\ci f)$,
that is, an {\it isomorphism of functors}. Thus, $\fF$ only
respects composition of 1-morphisms up to 2-isomorphism.
The 2-isomorphisms $\ep_{g,f}$ are often omitted in proofs.
\end{itemize}

A $\K$-{\it stack\/} is a {\it stack in groupoids\/} on
$\Sch_\K$. That is, it is a prestack $\fF$ satisfying the
following axioms. Let $\{f_i:U_i\ra V\}_{i\in I}$ be an
open cover of $V$ in the site $\Sch_\K$. Write $U_{ij}=
U_i\t_{f_i,V,f_j}U_j$ for the fibre product scheme and
$f_{ij}:U_{ij}\ra V$, $f_{ij,i}:U_{ij}\ra U_i$,
$f_{ij,j}:U_{ij}\ra U_j$ for the projections, and
similarly for `triple intersections' $U_{ijk}$. Then
\begin{itemize}
\setlength{\itemsep}{0pt}
\setlength{\parsep}{0pt}
\item[(i)] (Glueing of morphisms). If $X,Y\in\Obj(\fF(V))$ and
$\phi_i:\fF(f_i)X\ra\fF(f_i)Y$ are morphisms for $i\in I$ such
that
\end{itemize}
\e
\begin{split}
&\ep_{f_i,f_{ij,i}}(Y)\ci
\bigl(\fF(f_{ij,i})\phi_i\bigr)
\ci\ep_{f_i,f_{ij,i}}(X)^{-1}=\\
&\ep_{f_j,f_{ij,j}}(Y)\ci
\bigl(\fF(f_{ij,j})\phi_j\bigr)
\ci\ep_{f_j,f_{ij,j}}(X)^{-1}
\end{split}
\label{aa2eq4}
\e
\begin{itemize}
\setlength{\itemsep}{0pt}
\setlength{\parsep}{0pt}
\item[]in $\Hom\bigl(\fF(f_{ij})X,\fF(f_{ij})Y\bigr)$ for all
$i,j$, then there exists a morphism $\eta:X\ra Y$ in
$\Mor(\fF(V))$ with~$\fF(f_i)\eta=\phi_i$.
\item[(ii)] (Monopresheaf). If $X,Y$ lie in $\Obj(\fF(V))$ and
$\phi,\psi:X\!\ra\!Y$ in $\Mor(\fF(V))$ with $\fF(f_i)\phi=\fF(f_i)\psi$
for all $i\in I$, then~$\phi=\psi$.
\item[(iii)] (Glueing of objects). If $X_i\!\in\!\Obj(\fF(U_i))$
and $\phi_{ij}:\fF(f_{ij,j})X_j\!\ra\!\fF(f_{ij,i})X_i$ are
morphisms for all $i,j$ satisfying the {\it cocycle condition}
\e
\begin{split}
&\bigl[\ep_{f_{ij,i},f_{ijk,ij}}(X_i)\ci
\bigl(\fF(f_{ijk,ij})\phi_{ij}\bigr)\ci
\ep_{f_{ij,j},f_{ijk,ij}}(X_j)^{-1}\bigr]\ci\\
&\bigl[\ep_{f_{jk,j},f_{ijk,jk}}(X_j)\ci
\bigl(\fF(f_{ijk,jk})\phi_{jk}\bigr)
\ep_{f_{jk,k},f_{ijk,jk}}(X_k)^{-1}\bigr]=\\
&\bigl[\ep_{f_{ik,i},f_{ijk,ik}}(X_i)\ci
\bigl(\fF(f_{ijk,ik})\phi_{ik}\bigr)
\ep_{f_{ik,k},f_{ijk,ik}}(X_k)^{-1}\bigr]
\end{split}
\label{aa2eq5}
\e
in $\Hom\bigl(\fF(f_{ijk,k})X_k,\fF(f_{ijk,k})X_i\bigr)$
for all $i,j,k$, then there exists $X\in\Obj(\fF(V))$ and
isomorphisms $\phi_i:\fF(f_i)X\ra X_i$ in $\Mor(\fF(U_i))$
such that
\e
\phi_{ji}\ci\bigl(\fF(f_{ij,i})\phi_i\bigr)\ci\ep_{f_i,f_{ij,i}}(X)^{-1}
=\bigl(\fF(f_{ij,j})\phi_j\bigr)\ci\ep_{f_j,f_{ij,j}}(X)^{-1}
\label{aa2eq6}
\e
in $\Hom\bigl(\fF(f_{ij})X,\fF(f_{ij,j})X_j\bigr)$ for all~$i,j$.
\end{itemize}
Using notation from \cite[\S 3]{LaMo}, an {\it algebraic
$\K$-stack\/} is a $\K$-stack $\fF$ such that
\begin{itemize}
\setlength{\itemsep}{0pt}
\setlength{\parsep}{0pt}
\item[(i)] The {\it diagonal\/} $\De_\fF$ is
representable, quasicompact and separated.
\item[(ii)] There exists a scheme $U$ and a smooth
surjective 1-morphism $u:U\ra\fF$. We call $U,u$ an
{\it atlas} for~$\fF$.
\end{itemize}
\label{aa2def7}
\end{dfn}

As in \cite[\S 2.2]{Gome}, \cite[\S 3]{LaMo}, $\K$-stacks form
a 2-{\it category}. A 1-{\it morphism} $\phi:\fF\ra\fG$ is a
natural transformation between the 2-functors $\fF,\fG$. A
2-{\it morphism} $\al:\phi\ra\psi$ of 1-morphisms $\phi,\psi:
\fF\ra\fG$ is an isomorphism of the natural transformations.
{\it Representable} and {\it finite type} 1-morphisms are
defined in \cite[Def.s~3.9, 4.16]{LaMo}, and will be
important in~\S\ref{aa82}.

We define the set of $\K$-{\it points} of a stack.

\begin{dfn} Let $\fF$ be a $\K$-stack. Regarding
$\Spec\K$ as a $\K$-stack, we can form $\Hom(\Spec\K,\fF)$.
This is a {\it groupoid}, with {\it objects} 1-morphisms
$\Spec\K\ra\fF$, and {\it morphisms} 2-morphisms. Define
$\fF(\K)$ to be the set of {\it isomorphism classes of
objects} in $\Hom(\Spec\K,\fF)$. Elements of $\fF(\K)$ are
called $\K$-{\it points}, or {\it geometric points}, of
$\fF$. If $\phi:\fF\ra\fG$ is a 1-morphism of $\K$-stacks
then composition with $\phi$ yields a morphism of
groupoids $\Hom(\Spec\K,\fF)\ra\Hom(\Spec\K,\fG)$, and
therefore induces a map of sets~$\phi_*:\fF(\K)\ra\fG(\K)$.
\label{aa2def8}
\end{dfn}

One important difference in working with 2-categories rather
than ordinary categories is that in diagram-chasing one only
requires 1-morphisms to be 2-{\it isomorphic} rather than
{\it equal}. The simplest kind of {\it commutative diagram} is:
\begin{equation*}
\xymatrix@R=10pt{
& \fG \ar@{=>}[d]^F \ar[dr]^\psi \\
\fF \ar[ur]^\phi \ar[rr]_\chi && \fH,
}
\end{equation*}
by which we mean that $\fF,\fG,\fH$ are $\K$-{\it stacks},
$\phi,\psi,\chi$ are 1-{\it morphisms}, and $F:\psi\ci\phi\ra\chi$
is a 2-{\it isomorphism}. Sometimes we omit $F$, and mean
that~$\psi\ci\phi\cong\chi$.

\begin{dfn} Let $\phi:\fF\ra\fH$, $\psi:\fG\ra\fH$ be 1-morphisms
of $\K$-stacks. Define the {\it fibre product\/} $\fF\t_{\phi,\fH,
\psi}\fG$, or $\fF\t_\fH\fG$ for short, such that $(\fF\t_\fH\fG)(U)$
for $U\in\Sch_\K$ is the groupoid with {\it objects} triples
$(A,B,\al)$ for $A\in\Obj(\fF(U))$, $B\in\Obj(\fG(U))$ and
$\al:\phi(U)A\ra\psi(U)B$ in $\Mor(\fH(U))$, and
{\it morphisms} $(\be,\ga):(A,B,\al)\ra(A',B',\al')$
for $\be:A\ra A'$ in $\Mor(\fF(U))$ and $\ga:B\ra B'$ in
$\Mor(\fG(U))$ with $\al'\ci(\phi(U)\be)=(\psi(U)\ga)\ci\al$
in~$\Mor(\fH(U))$.

The rest of the data to make a {\it contravariant\/ $2$-functor}
$\fF\t_\fH\fG:\Sch_\K\ra\gro$ is defined in the obvious way. Then
$\fF\t_\fH\fG$ is a {\it $\K$-stack}, with projection 1-morphisms
$\pi_\fF:\fF\t_\fH\fG\ra\fF$, $\pi_\fG:\fF\t_\fH\fG\ra\fG$ fitting
into a {\it commutative diagram}:
\e
\begin{gathered}
\xymatrix@R=2pt{
& \fF \ar[dr]^\phi \ar@{=>}[dd]^A \\
\fF\t_\fH\fG
\ar[dr]_{\pi_\fG} \ar[ur]^{\pi_\fF} && \fH,\\
& \fG \ar[ur]_\psi \\
}
\end{gathered}
\label{aa2eq7}
\e
where $A:\phi\ci\pi_\fF\ra\psi\ci\pi_\fG$ is a 2-{\it isomorphism}.
Any commutative diagram
\e
\begin{gathered}
\xymatrix@R=2pt{
& \fF \ar[dr]^\phi \ar@{=>}[dd]^B \\
\fE
\ar[dr]_\eta \ar[ur]^\th && \fH\\
& \fG \ar[ur]_\psi \\
}
\end{gathered}
\label{aa2eq8}
\e
extends to a commutative diagram with $\rho$ unique up to 2-isomorphism
\e
\begin{gathered}
\xymatrix@R=5pt@C=15pt{
&& \fF \ar[dr]^\phi \ar@{=>}[dd]^A \\
\fE
\ar@(dr,l)[drr]_(0.2)\eta \ar@(ur,l)[urr]^(0.2)\th
\ar[r]^\rho
&
\fF\t_\fH\fG \ar@{=>}[u]^C \ar@{=>}[d]_D
\ar[dr]^(0,6){\pi_\fG} \ar[ur]_(0.6){\pi_\fF} && \fH,\\
&& \fG \ar[ur]_\psi \\
}
\end{gathered}
\;\>
\begin{gathered}
\text{with}\\
\text{commuting}\\
\text{2-morphisms}
\end{gathered}
\;\>
\begin{gathered}
\xymatrix@R=17pt@C=18pt{
\phi\ci\pi_\fF\ci\rho \ar[r]_{\phi\ci C} \ar[d]_{A\ci\rho} & \phi\ci\th
\ar[d]^B \\
\psi\ci\pi_\fG\ci\rho \ar[r]^{\psi\ci D} & \psi\ci\eta.
}
\end{gathered}
\label{aa2eq9}
\e

We call \eq{aa2eq8} a {\it Cartesian square} if $\rho$ in
\eq{aa2eq9} is a 1-{\it isomorphism}, so that $\fE$ is
1-isomorphic to $\fF\t_\fH\fG$. Cartesian squares may
also be characterized by a {\it universal property}, as
in \cite[\S 3]{BEFF}. By the {\it product\/} $\fF\t\fG$
of $\K$-stacks $\fF,\fG$ we mean their fibre product
$\fF\t_{\Spec\K}\fG$ over $\Spec\K$, using the natural
projections~$\fF,\fG\ra\Spec\K$.
\label{aa2def9}
\end{dfn}

Here are some properties of these that we will need later.
\begin{itemize}
\setlength{\itemsep}{0pt}
\setlength{\parsep}{0pt}
\item If $\fF,\fG,\fH$ are {\it algebraic} $\K$-stacks, then
their fibre product $\fF\t_\fH\fG$ is also algebraic, by
\cite[Prop.~4.5(i)]{LaMo}. Hence, in any Cartesian square
\eq{aa2eq8}, if $\fF,\fG,\fH$ are algebraic then $\fE$ is
also algebraic, as it is 1-isomorphic to $\fF\t_\fH\fG$. If
also $\fF,\fG,\fH$ are {\it locally of finite type}, then
so are $\fF\t_\fH\fG$ and~$\fE$.
\item Let $\fF,\fG,\fH$ in \eq{aa2eq7} be {\it algebraic
$\K$-stacks}. We may think of $\fG$ as a stack over a
{\it base} $\fH$, and then $\pi_\fF:\fF\t_\fH\fG\ra\fF$
is obtained from $\psi:\fG\ra\fH$ by {\it base extension}.
Therefore, for any property P of morphisms of algebraic
$\K$-stacks that is {\it stable under base extension},
if $\psi$ has P then $\pi_\fF$ has~P.

But in any Cartesian square \eq{aa2eq8}, $\fE$ is
1-isomorphic to $\fF\t_\fH\fG$. Thus if $\fE,\fF,\fG,\fH$
are algebraic and $\psi$ has P in \eq{aa2eq8}, then $\th$ has
P, as in \cite[Rem.~4.14.1]{LaMo}. This holds when $P$ is
{\it of finite type} \cite[Rem.~4.17(2)]{LaMo}, or {\it
locally of finite type} \cite[p.~33-4]{LaMo}, or {\it
representable}~\cite[Lem.~3.11]{LaMo}.
\end{itemize}

\section{Refining the Jordan--H\"older Theorem}
\label{aa3}

We shall study the following situation.

\begin{dfn} Let $\A$ be an abelian category of finite length,
and $X\in\A$. Then $X$ admits a composition series
$0=A_0\subset A_1\subset\cdots\subset A_n=X$ by Theorem
\ref{aa2thm1}, and the simple factors $S_k=A_k/A_{k-1}$
for $k=1,\ldots,n$ of $X$ are independent of choices, up to
isomorphism and permutation of $1,\ldots,n$. Suppose
$S_k\not\cong S_l$ for $1\le k<l\le n$. Then we say that $X$
has {\it nonisomorphic simple factors}.

Let $X$ have nonisomorphic simple factors, and let $I$ be an
{\it indexing set\/} for $\{S_1,\ldots,S_n\}$, so that $\md{I}=n$,
and write $\{S_1,\ldots,S_n\}=\{S^i:i\in I\}$. Then Theorem
\ref{aa2thm1} implies that for every composition series
$0=B_0\subset B_1\subset\cdots\subset B_n=X$ for $X$ with
simple factors $T_k=B_k/B_{k-1}$, there exists a unique
bijection $\phi:I\ra\{1,\ldots,n\}$ such that $S^i\cong
T_{\phi(i)}$ for all~$i\in I$.

Define a {\it partial order} $\pr$ on $I$ by $i\pr j$
for $i,j\in I$ if and only if $\phi(i)\le\phi(j)$ for all
bijections $\phi:I\ra\{1,\ldots,n\}$ constructed from a
composition series $0=B_0\subset B_1\subset\cdots\subset B_n=X$
for $X$ as above. Then $(I,\pr)$ is a {\it partially ordered
set}, or {\it poset} for short.
\label{aa3def1}
\end{dfn}

The point of this definition is to treat all the Jordan--H\"older
composition series $0=B_0\subset B_1\subset\cdots\subset B_n=X$
for $X$ on an equal footing. Writing the simple factors of
$X$ as $S_1,\ldots,S_n$ gives them a preferred order, and
favours one composition series over the rest. So instead we
write the simple factors as $S^i$ for $i\in I$, some arbitrary
indexing set. Here is some more notation.

\begin{dfn} Let $(I,\pr)$ be a finite poset. Define
$J\subseteq I$ to be
\begin{itemize}
\setlength{\itemsep}{0pt}
\setlength{\parsep}{0pt}
\item[(i)] an {\it s-set\/} if $i\in I$, $j\in J$ and
$i\pr j$ implies~$i\in J$,
\item[(ii)] a {\it q-set\/} if $i\in I$, $j\in J$ and
$j\pr i$ implies $i\in J$, and
\item[(iii)] an {\it f-set\/} if $i\in I$ and $h,j\in J$
and $h\pr i\pr j$ implies~$i\in J$.
\end{itemize}
The motivation for this is that below s-sets will correspond to
{\it subobjects} $S\subset X$, q-sets to {\it quotient objects}
$X/S$, and f-sets to {\it factors} $T/S$ for~$S\subset T\subset X$.
\label{aa3def2}
\end{dfn}

Here are some properties of s-sets, q-sets and f-sets.

\begin{prop} Let\/ $(I,\pr)$ be a finite poset. Then
\begin{itemize}
\setlength{\itemsep}{0pt}
\setlength{\parsep}{0pt}
\item[{\rm(a)}] $I$ and\/ $\emptyset$ are s-sets. If\/ $J,K$ are
s-sets then $J\cap K$ and $J\cup K$ are s-sets.
\item[{\rm(b)}] $J$ is an s-set if and only if\/ $I\sm J$
is a q-set.
\item[{\rm(c)}] If\/ $J\subset K$ are s-sets then $K\sm J$
is an f-set. Every f-set is of this form.
\item[{\rm(d)}] If\/ $J,K$ are f-sets then $J\cap K$ is an f-set,
but\/ $J\cup K$ may not be an f-set.
\end{itemize}
\label{aa3prop1}
\end{prop}

The proof is elementary, and left as an exercise. For the
last part of (c), if $F\subseteq I$ is an f-set, define
$K=\{i\in I:i\pr j$ for some $j\in F\}$ and $J=K\sm F$.
Then $J\subset K$ are s-sets with $F=K\sm J$. Note that
(a) and (b) imply the collections of s-sets and q-sets are
both {\it topologies} on $I$, but (d) shows the f-sets may
not be. Also $\pr$ can be reconstructed from the set of
s-sets on $I$, as $i\pr j$ if and only if $i\in J$ for
every s-set $J\subset I$ with~$j\in J$.

\begin{lem} In the situation of Definition \ref{aa3def1}, suppose
$S\subset X$ is a subobject. Then there exists a unique s-set
$J\subseteq I$ such that the simple factors in any composition
series for $S$ are isomorphic to $S^i$ for~$i\in J$.
\label{aa3lem1}
\end{lem}

\begin{proof} Let $0=B_0\subset\cdots\subset B_l=S$ be a
composition series for $S$, with simple factors $T_k=B_k/B_{k-1}$
for $k=1,\ldots,l$. Then $0=B_0\subset\cdots\subset B_l\subset X$ is
a filtration of $X$ without repetitions, and can be refined to a
composition series by Theorem \ref{aa2thm1}. As $T_k$ is simple,
no extra terms are inserted between $B_{k-1}$ and $B_k$. Thus
$X$ has a composition series $0=B_0\subset\cdots\subset B_l
\subset\cdots\subset B_n=X$, with simple factors $T_k=B_k/B_{k-1}$
for~$k=1,\ldots,n$.

By Definition \ref{aa3def1} there is a unique bijection
$\phi:I\ra\{1,\ldots,n\}$ such that $S^i\cong T_{\phi(i)}$
for all $i\in I$. Define $J=\phi^{-1}(\{1,\ldots,l\})$.
Then $J\subseteq I$, and the simple factors $T_k$ of the
composition series $0=B_0\subset\cdots\subset B_l=S$
are isomorphic to $S^i$ for $i\in J$. Theorem \ref{aa2thm1} then
implies that the simple factors in {\it any} composition
series for $S$ are isomorphic to $S^i$ for~$i\in J$.

Uniqueness of $J$ is now clear, as a different $J$ would
give different simple factors for $S$. Suppose $j\in J$ and
$i\in I\sm J$. Then $1\le\phi(j)\le l$ and $l+1\le\phi(i)\le n$,
so $\phi(j)<\phi(i)$, which implies that $i\npr j$ by
Definition \ref{aa3def1}. Hence if $j\in J$ and $i\in I$ with
$i\pr j$ then $i\in J$, and $J$ is an s-set.
\end{proof}

\begin{lem} Suppose $S,T\subset X$ correspond to s-sets
$J,K\subseteq I$, as in Lemma \ref{aa3lem1}. Then $S\cap T$
corresponds to $J\cap K$, and\/ $S+T$ corresponds to~$J\cup K$.
\label{aa3lem2}
\end{lem}

\begin{proof} Let $S\cap T$ correspond to the s-set
$L\subseteq I$, and $S+T$ to the s-set $M$. We must
show $L=J\cap K$ and $M=J\cup K$. By Theorem \ref{aa2thm1}
we may refine the filtration $0\subset S\cap T\subset S$ to
a composition series for $S$ containing one for $S\cap T$.
Thus the simple factors of $S$ contain those of $S\cap T$,
and $L\subseteq J$. Similarly $L\subseteq K$, so
$L\subseteq J\cap K$, and $J\cup K\subseteq M$
as~$S,T\subseteq S\cap T$.

Now the simple factors of $S/(S\cap T)$ are $S^i$ for $i\in
J\sm L$, and the simple factors of $(S+T)/T$ are $S^i$ for
$i\in M\sm K$. As $S/(S\cap T)\cong(S+T)/T$ by \eq{aa2eq2} we
see that $J\sm L=M\sm K$. Together with $L\subseteq J\cap K$ and
$J\cup K\subseteq M$ this implies that $L=J\cap K$ and~$M=J\cup K$.
\end{proof}

\begin{lem} Suppose $S,T\subset X$ correspond to s-sets
$J,K\subseteq I$. Then $J\subseteq K$ if and only if\/
$S\subset T\subset X$, and\/ $J=K$ if and only if\/~$S=T$.
\label{aa3lem3}
\end{lem}

\begin{proof} If $S\subset T\subset X$ we can refine $0\subset
S\subset T\subset X$ to a composition series $0=B_0\subset
B_1\subset\cdots\subset B_n=X$ with $S=B_k$ and $T=B_l$ for
$0\le k\le l\le n$. Let $T_m=B_m/B_{m-1}$. Then the simple factors
of $S$ are $T_1,\ldots,T_k$ and of $T$ are $T_1,\ldots,T_l$.
Hence $J\subseteq K$, as $k\le l$. This proves the first~`if'.

Now suppose $S,T\subset X$ and $J\subseteq K$. Then $J\cap K=J$,
so $S\cap T$ corresponds to the s-set $J$. But $S\cap T
\subset S$, so $S/(S\cap T)$ has no simple factors, and
$S=S\cap T$. Thus $S\subset T$, proving the first `only if'.
The second part is immediate.
\end{proof}

\begin{lem} Let\/ $j\in I$ and define $J^j=\{i\in I:i\pr j\}$.
Then $J^j$ is an s-set, and there exists a subobject\/ $D^j\subset X$
corresponding to~$J^j$.
\label{aa3lem4}
\end{lem}

\begin{proof} Clearly $J^j$ is an s-set. By Definition \ref{aa3def1}
each composition series $0=B_0\subset B_1\subset\cdots\subset B_n=X$
for $X$ gives a bijection $\phi:I\ra\{1,\ldots,n\}$. Let $\phi_1,
\ldots,\phi_r$ be the distinct bijections $\phi:I\ra\{1,\ldots,n\}$
realized by composition series for $X$. For each $k=1,\ldots,r$ choose
a composition series $0=B_0\subset B_1\subset\cdots\subset B_n=X$
with bijection $\phi_k$, and define $C_k$ to be the subobject
$B_{\phi_k(j)}\subset X$.

This defines subobjects $C_1,\ldots,C_r\subset X$, where $C_k$
corresponds to the s-set $\phi_k^{-1}(\{1,\ldots,\phi_k(j)\})
\subseteq I$. Define $D^j=C_1\cap C_2\cap\cdots\cap C_r$.
Then $S\subset X$, and Lemma \ref{aa3lem2} shows that $D^j$ corresponds
to the s-set
\begin{gather*}
\bigcap_{k=1}^r\phi_k^{-1}\bigl(\{1,\ldots,\phi_k(j)\}\bigr)
=\bigcap_{k=1}^r\bigl\{i\in I:\phi_k(i)\le\phi_k(j)\bigr\}=\\
\bigl\{i\in I:\text{$\phi_k(i)\le\phi_k(j)$ for all $k=1,\ldots,r$}\bigr\}
=\{i\in I:i\pr j\}=J^j,
\end{gather*}
by definition of $\pr$.
\end{proof}

We can now {\it classify subobjects} of $X$ in terms of s-sets.

\begin{prop} In the situation of Definitions \ref{aa3def1} and
\ref{aa3def2}, for each s-set\/ $J\subseteq I$ there exists a
unique subobject\/ $S\subset X$ such that the simple factors
in any composition series for $S$ are isomorphic to $S^i$
for $i\in J$. This defines a $1$-$1$ correspondence between
subobjects $S\subset X$ and s-sets~$J\subseteq I$.
\label{aa3prop2}
\end{prop}

\begin{proof} For each $j\in J$ define $J^j$ and $D^j$ as in Lemma
\ref{aa3lem4}. Then $j\in J^j\subseteq J$, so $J=\bigcup_{j\in J}J^j$.
Set $S=\sum_{j\in J}D^j$. Then $S\subset X$ corresponds to the s-set
$\bigcup_{j\in J}J^j=J$ by Lemma \ref{aa3lem2}. Uniqueness follows
from Lemma~\ref{aa3lem3}.
\end{proof}

The dual proof classifies {\it quotient objects} of $X$ in terms
of q-sets.

\begin{prop} In the situation of Definitions \ref{aa3def1} and
\ref{aa3def2}, for each q-set\/ $K\subseteq I$ there exists a
unique quotient object\/ $Q=X/S$ of $X$ such that the simple
factors in any composition series for $Q$ are isomorphic to
$S^i$ for $i\in K$. This defines a $1$-$1$ correspondence
between quotient objects and q-sets.
\label{aa3prop3}
\end{prop}

We can also classify {\it composition series} for~$X$.

\begin{prop} In the situation of Definition \ref{aa3def1}, for
each bijection $\phi:I\ra\{1,\ldots,n\}$ there exists a unique
composition series $0=B_0\subset\cdots\subset B_n=X$ with\/
$S^i\cong B_{\phi(i)}/B_{\phi(i)-1}$ for all\/ $i\in I$ if and
only if\/ $i\pr j$ implies~$\phi(i)\le\phi(j)$.
\label{aa3prop4}
\end{prop}

\begin{proof} The `only if' part follows from Definition
\ref{aa3def1}. For the `if' part, let $\phi:I\ra\{1,\ldots,n\}$
be a bijection for which $i\pr j$ implies that $\phi(i)\le
\phi(j)$. Then $\phi^{-1}(\{1,\ldots,k\})$ is an s-set
for each $k=0,1,\ldots,n$. Let $B_k\subset X$ be the unique
subobject corresponding to $\phi^{-1}(\{1,\ldots,k\})$,
which exists by Proposition \ref{aa3prop2}. It easily follows
that $0=B_0\subset B_1\subset\cdots\subset B_n=X$ is the unique
composition series with $B_k/B_{k-1}\cong S^{\phi^{-1}(k)}$
for $k=1,\ldots,n$, and the result follows.
\end{proof}

This implies that composition series for $X$ up to isomorphism
are in 1-1 correspondence with {\it total orders} on $I$ compatible
with the partial order $\pr$. In Definition \ref{aa3def1} we
defined the partial order $\pr$ on $I$ to be the intersection
of all the total orders on $I$ coming from composition series
for $X$. We now see that {\it every} total order on $I$
compatible with $\pr$ comes from a composition series.

\section{Posets $(I,\pr)$ and $(I,\pr)$-configurations in $\A$}
\label{aa4}

Although a {\it subobject\/} of $X$ is an equivalence class
of injective $i:S\ra X$, in \S\ref{aa3} we for simplicity
suppressed the morphisms $i$, and just wrote $S\subset X$.
We shall now change our point of view, and investigate the
natural morphisms between the factors $T/S$ of $X$. Therefore
we adopt some new notation, which stresses morphisms between
objects. The following definition encodes the properties we
expect of the factors of $X$, and their natural morphisms.

\begin{dfn} Let $(I,\pr)$ be a finite poset, and use the
notation of Definition \ref{aa3def2}. Define $\F_\sIp$ to be the
set of f-sets of $I$. Define $\G_\sIp$ to be the subset of
$(J,K)\in\F_\sIp\t\F_\sIp$ such that $J\subseteq K$, and if
$j\in J$ and $k\in K$ with $k\pr j$, then $k\in J$.
Define $\H_\sIp$ to be the subset of $(J,K)\in\F_\sIp\t\F_\sIp$ such
that $K\subseteq J$, and if $j\in J$ and $k\in K$ with
$k\pr j$, then $j\in K$. It is easy to show that $\G_\sIp$
and $\H_\sIp$ have the following properties:
\begin{itemize}
\setlength{\itemsep}{0pt}
\setlength{\parsep}{0pt}
\item[(a)] $(J,K)$ lies in $\G_\sIp$ if and only if $(K,K\sm J)$ lies
in~$\H_\sIp$.
\item[(b)] If $(J,K)\in\G_\sIp$ and $(K,L)\in\G_\sIp$ then $(J,L)\in\G_\sIp$.
\item[(c)] If $(J,K)\in\H_\sIp$ and $(K,L)\in\H_\sIp$ then $(J,L)\in\H_\sIp$.
\item[(d)] If $(J,K)\!\in\!\G_\sIp$, $(K,L)\!\in\!\H_\sIp$ then
$(J,J\!\cap\!L)\!\in\!\H_\sIp$, $(J\!\cap\!L,L)\!\in\!\G_\sIp$.
\end{itemize}

Let $\A$ be an {\it abelian category}, or more generally an {\it exact
category}, as in \S\ref{aa21}. Define an $(I,\pr)$-{\it configuration}
$(\si,\io,\pi)$ in $\A$ to be maps $\si:\F_\sIp\ra\Obj(\A)$,
$\io:\G_\sIp\ra\Mor(\A)$, and $\pi:\H_\sIp\ra\Mor(\A)$, where
\begin{itemize}
\setlength{\itemsep}{0pt}
\setlength{\parsep}{0pt}
\item[(i)] $\si(J)$ is an object in $\A$ for $J\in\F_\sIp$,
with~$\si(\emptyset)=0$.
\item[(ii)] $\io(J,K):\si(J)\!\ra\!\si(K)$ is {\it injective}
for $(J,K)\!\in\!\G_\sIp$, and~$\io(J,J)\!=\!\id_{\si(J)}$.
\item[(iii)] $\pi(J,K)\!:\!\si(J)\!\ra\!\si(K)$ is {\it surjective}
for $(J,K)\!\in\!\H_\sIp$, and~$\pi(J,J)\!=\!\id_{\si(J)}$.
\end{itemize}
These should satisfy the conditions:
\begin{itemize}
\setlength{\itemsep}{0pt}
\setlength{\parsep}{0pt}
\item[(A)] Let $(J,K)\in\G_\sIp$ and set $L=K\sm J$. Then the
following is exact in~$\A$:
\e
\xymatrix@C=40pt{ 0 \ar[r] &\si(J) \ar[r]^{\io(J,K)} &\si(K)
\ar[r]^{\pi(K,L)} &\si(L) \ar[r] & 0. }
\label{aa4eq1}
\e
\item[(B)] If $(J,K)\in\G_\sIp$ and $(K,L)\in\G_\sIp$
then~$\io(J,L)=\io(K,L)\ci\io(J,K)$.
\item[(C)] If $(J,K)\in\H_\sIp$ and $(K,L)\in\H_\sIp$
then~$\pi(J,L)=\pi(K,L)\ci\pi(J,K)$.
\item[(D)] If $(J,K)\in\G_\sIp$ and $(K,L)\in\H_\sIp$ then
\e
\pi(K,L)\ci\io(J,K)=\io(J\cap L,L)\ci\pi(J,J\cap L).
\label{aa4eq2}
\e
\end{itemize}
Note that (A)--(D) make sense because of properties (a)--(d),
respectively.

A {\it morphism} $\al:(\si,\io,\pi)\ra(\si',\io',\pi')$ of
$(I,\pr)$-configurations in $\A$ is a collection of morphisms
$\al(J):\si(J)\ra\si'(J)$ for each $J\in\F_\sIp$ satisfying
\e
\begin{aligned}
\al(K)\ci\io(J,K)&=\io'(J,K)\ci\al(J)&&
\text{for all $(J,K)\in\G_\sIp$, and}\\
\al(K)\ci\pi(J,K)&=\pi'(J,K)\ci\al(J)&&
\text{for all $(J,K)\in\H_\sIp$.}
\end{aligned}
\label{aa4eq3}
\e
It is an {\it isomorphism} if $\al(J)$ is an isomorphism for
all $J\in\F_\sIp$. Morphisms compose in the obvious way.
\label{aa4def}
\end{dfn}

We now show that Definition \ref{aa4def} captures the properties
of the families of subobjects $S^J\subset X$ considered in~\S\ref{aa3}.

\begin{thm} Let\/ $(I,\pr)$ be a finite poset, $\A$ an abelian
category, and\/ $X\in\A$. Suppose that for each s-set\/ $J\subseteq I$
we are given a subobject\/ $S^J\subset X$, such that
\e
S^\emptyset=0,\quad S^I=X,\quad S^A\cap S^B=S^{A\cap B}
\quad\text{and\/}\quad S^A+S^B=S^{A\cup B}
\label{aa4eq4}
\e
for all s-sets $A,B\subseteq I$. Then there exists an
$(I,\pr)$-configuration $(\si,\io,\pi)$ in $\A$ with\/
$\si(I)=X$ such that\/ $\io(J,I):\si(J)\ra X$ represents
$S^J\subset X$ for each s-set\/ $J\subseteq I$. This
$(\si,\io,\pi)$ is unique up to canonical isomorphism
in~$\A$.
\label{aa4thm1}
\end{thm}

\begin{proof} Throughout (i)--(iii) and (A)--(D) will refer to
Definition \ref{aa4def}. We divide the proof into the following
seven steps:
\begin{list}{}{
\setlength{\itemsep}{1pt}
\setlength{\parsep}{1pt}
\setlength{\labelwidth}{50pt}
\setlength{\leftmargin}{50pt}
}
\item[{\bf Step 1.}] Define $\si$ and $\io$ on s-sets,
and prove (B) for s-sets.
\item[{\bf Step 2.}] For $J,K$ s-sets with $J\cap K=\emptyset$,
show~$\si(J\cup K)\cong\si(J)\op\si(K)$.
\item[{\bf Step 3.}] Define $\si$ on f-sets and $\pi(J,L)$ for
s-sets~$J$.
\item[{\bf Step 4.}] Complete the definitions of $\io,\pi$, and
prove~(A).
\item[{\bf Step 5.}] Prove partial versions of (C), (D), mixing
s-sets and f-sets.
\item[{\bf Step 6.}] Prove (B), (C), and $\io(J,J)=\pi(J,J)=
\id_{\si(J)}$ in (ii) and~(iii).
\item[{\bf Step 7.}] Prove (D).
\end{list}

\noindent{\bf Step 1.} For each s-set $J\subseteq I$, choose
$\si(J)\in\A$ and an injective morphism $\io(J,I):\si(J)\ra X$
representing $S^J\subset X$. Then $\si(J)$ and $\io(J,I)$ are
unique up to canonical isomorphism. In particular, choose
$\si(\emptyset)=0$ as in (i), $\si(I)=X$, and $\io(I,I)=\id_X$.
Suppose $J\subseteq K$ are s-sets. Then \eq{aa4eq4} implies
that $S^J\subset S^K\subset X$. Hence there exists a unique,
injective $\io(J,K):\si(J)\ra\si(K)$ such that
\e
\io(J,I)=\io(K,I)\ci\io(J,K)
\quad\text{for $J\subseteq K$ s-sets, as in (B).}
\label{aa4eq5}
\e
By uniqueness the two definitions of $\io(J,I)$ coincide,
and~$\io(J,J)=\id_{\si(J)}$.

Suppose $J\subseteq K\subseteq L$ are s-sets. Applying
\eq{aa4eq5} to $(K,L),(J,K),(J,L)$ gives
\begin{equation*}
\io(L,I)\ci\io(K,L)\ci\io(J,K)
=\io(K,I)\ci\io(J,K)=\io(J,I)
=\io(L,I)\ci\io(J,L).
\end{equation*}
Since $\io(L,I)$ is injective we can cancel it from both
sides, so that
\e
\io(J,L)=\io(K,L)\ci\io(J,K),
\quad\text{for $J\subseteq K\subseteq L$ s-sets, as in (B).}
\label{aa4eq6}
\e

\noindent{\bf Step 2.} Let $J,K$ be s-sets with $J\cap K=\emptyset$.
We shall show that
\e
\io(J,J\cup K)\ci\pi_{\si(J)}+\io(K,J\cup K)\ci\pi_{\si(K)}:
\si(J)\op\si(K)\ra\si(J\cup K)
\label{aa4eq7}
\e
is an {\it isomorphism}. Apply Definition \ref{aa2def5} with
$\io(J,I):\si(J)\ra X$ in place of $i:S\ra X$, and $\io(K,I):
\si(K)\ra X$ in place of $j:T\ra X$. By \eq{aa4eq4} we may
take $U=\si(I\cap J)=\si(\emptyset)=0$, $V=\si(J\cup K)$ and
$e=\io(J\cup K,I)$. The definition gives $c:\si(J)\ra\si(J\cup K)$
with $\io(J,I)=\io(J\cup K,I)\ci c$, so $c=\io(J,J\cup K)$ by
\eq{aa4eq5} and injectivity of $\io(J\cup K,I)$. Similarly
$d=\io(K,J\cup K)$. Thus \eq{aa4eq7} is the second map in
\eq{aa2eq1}. As $U=0$, exactness implies \eq{aa4eq7} is an
isomorphism.
\medskip

\noindent{\bf Step 3.} Let $L\subseteq I$ be an f-set which is
not an s-set or a q-set, and define $J'=\{i\in I\sm L:l\npr i$
for all $l\in L\}$ and $K'=J'\cup L$. Then $J'\subset K'$ are s-sets
with $L=K'\sm J'$. Choose $\si(L)\in\A$ and a surjective $\pi(K',L):
\si(K')\ra\si(L)$ which is a {\it cokernel\/} for $\io(J',K'):\si(J')\ra
\si(K')$. Then $\si(L)$, $\pi(K',L)$ are unique up to canonical isomorphism.

If $L$ is an s-set then $J',L$ are s-sets with $J'\cap L=\emptyset$, and
Step 2 shows that $\si(K')\cong\si(J')\op\si(L)$, and we take $\pi(K',L)$
to be the natural projection with~$\pi(K',L)\ci\io(L,K')=\id_{\si(L)}$.

Now let $J\subset K$ be s-sets in $I$ with $K\sm J=L$. Then
$J\subseteq J'$ and $K\subseteq K'$, as $J',K'$ are defined to
be as large as possible, and $J'\cap K=J$, $J'\cup K'=K'$.
Let $c:\si(K)\ra C$ be a cokernel for $\io(J,K)$, and consider
the commutative diagram with rows short exact sequences
\e
\begin{gathered}
\xymatrix@R=13pt{
0 \ar[r] & \si(J) \ar[r]_{\io(J,K)} \ar[d]_{\io(J,J')} &
\si(K) \ar[r]_c \ar[d]^(0.45){\io(K,K')}
& C \ar[r] \ar@{-->}[d]^h & 0 \\
0 \ar[r] & \si(J') \ar[r]^{\io(J',K')} & \si(K')
\ar[r]^{\pi(K',L)} & \si(L) \ar[r] &0,
}
\end{gathered}
\label{aa4eq8}
\e
where $h$ is not yet constructed. The first square of \eq{aa4eq8}
commutes by \eq{aa4eq6}, so
\begin{equation*}
0=\pi(K',L)\ci\io(J',K')\ci\io(J,J')=\pi(K',L)\ci\io(K,K')\ci\io(J,K).
\end{equation*}
But $c$ is the cokernel of $\io(J,K)$, so there exists a unique
$h:C\ra\si(L)$ such that $h\ci c=\pi(K',L)\ci\io(K,K')$, that is,
the second square in \eq{aa4eq8} commutes.

As $S^J=S^{J'}\cap S^K$ and $S^{K'}=S^{J'}+S^K$ by \eq{aa4eq4},
equation \eq{aa2eq1} implies that
\begin{equation*}
\xymatrix@C=10pt{
0 \ar[r] & \si(J) \ar[rrrrrr]
\ar@{}@<.5ex>[rrrrrr]^(0.45){\io_{\si(J')}\ci\io(J,J')-
\io_{\si(K)}\ci\io(J,K)}
&&&&&& \si(J')\!\op\!\si(K) \ar[rrrrrr]
\ar@{}@<.5ex>[rrrrrr]^(0.58){\io(J',K')\ci\pi_{\si(J')}+
\io(K,K')\ci\pi_{\si(K)}} &&&&&& \si(K') \ar[r] & 0
}
\end{equation*}
is exact. As the composition of the first map with $c\ci\pi_{\si(K)}$
is zero we see that
\begin{equation*}
c\ci\pi_{\si(K)}=l\ci\bigl(\io(J',K')\ci\pi_{\si(J')}+
\io(K,K')\ci\pi_{\si(K)}\bigr)
\end{equation*}
for some unique $l:\si(K')\ra C$, by definition of cokernel.
Composing with $\io_{\si(J')}$ gives $l\ci\io(J',K')=0$, so
$l=m\ci\pi(K',L)$ for some unique $m:\si(L)\ra C$, by
exactness of the bottom line of \eq{aa4eq8}.

Then $m\ci h\ci c=m\ci\pi(K',L)\ci\io(K,K')=l\ci\io(K,K')=c=
\id_C\ci c$, so as $c$ is surjective we have $m\ci h=\id_C$.
Also $\pi(K',L)=h\ci l$, since
\begin{align*}
&h\ci l\ci\bigl(\io(J',K')\ci\pi_{\si(J')}+
\io(K,K')\ci\pi_{\si(K)}\bigr)=h\ci c\ci\pi_{\si(K)}=\\
&\pi(K',L)\!\ci\!\io(K,K')\!\ci\!\pi_{\si(K)}\!=\!
\pi(K',L)\!\ci\!\bigl(\io(J',K')\!\ci\!\pi_{\si(J')}\!+\!
\io(K,K')\!\ci\!\pi_{\si(K)}\bigr),
\end{align*}
and $\io(J',K')\ci\pi_{\si(J')}+\io(K,K')\ci\pi_{\si(K)}$
is surjective. Hence $h\ci m\ci\pi(K',L)=h\ci l=\pi(K',L)=
\id_{\si(L)}\ci\pi(K',L)$, and $h\ci m=\id_{\si(L)}$ as
$\pi(K',L)$ is surjective. Thus $m=h^{-1}$, and $h$ is
an {\it isomorphism}.

Define $\pi(K,L)=\pi(K',L)\ci\io(K,K')$, so that
$\pi(K,L)=h\ci c$ as \eq{aa4eq8} is commutative.
As $h$ is an isomorphism and $c$ is a cokernel for
$\io(J,K)$, we see that $\pi(K,L):\si(K)\ra\si(L)$ is
a {\it cokernel\/} for $\io(J,K):\si(J)\ra\si(K)$. Hence
$\pi(K,L)$ is {\it surjective}, and \eq{aa4eq1} is exact
when $J\subseteq K$ are s-sets.

Suppose now that $J,K$ are s-sets and $L$ is an f-set with
$(J,K)\in\G_\sIp$ and $(J,L),(K,L)\in\H_\sIp$. Define $K'$ using $L$
as above. Then $J\subseteq K\subseteq K'$ and
\begin{equation*}
\pi(J,L)=\pi(K',L)\ci\io(J,K')=\pi(K',L)\ci\io(K,K')
\ci\io(J,K)=\pi(K,L)\ci\io(J,K),
\end{equation*}
by \eq{aa4eq6} and the definitions of $\pi(J,L),\pi(K,L)$. Hence
\e
\pi(J,L)=\pi(K,L)\ci\io(K,J)
\quad\text{when $J,K$ are s-sets.}
\label{aa4eq9}
\e

\noindent{\bf Step 4.} Let $(J,K)\in\G_\sIp$, and define $A'=\{i\in I
\sm K:k\npr i$ for all $k\in K\}$, $B'=A'\cup J$ and $C'=A'\cup K$.
Then $A'\subseteq B'\subseteq C'$ are s-sets with $J=B'\sm A'$, $K=C'\sm A'$,
and $K\sm J=C'\sm B'$, and they are the {\it largest\/} s-sets with this
property. Consider the commutative diagram with exact rows and columns:
\e
\begin{gathered}
\xymatrix@R=1pt{
& & 0\ar[dd] & 0 \ar[dd] \\ \\
0 \ar[r] &\si(A') \ar[r]^{\io(A',B')} \ar@{=}[ddd]_{\id_{\si(A')}}
&\si(B') \ar[r]^{\pi(B',J)} \ar[ddd]_(0.4){\io(B',C')}
&\si(J) \ar[r] \ar@{-->}[ddd]^{\io(J,K)} & 0 \\ \\ \\
0 \ar[r] &\si(A') \ar[r]^{\io(A',C')}
&\si(C') \ar[r]^{\pi(C',K)} \ar[ddd]_{\pi(C',K\sm J)}
&\si(K) \ar[r] \ar@{-->}[ddd]^{\pi(K,K\sm J)} & 0 \\ \\ \\
&&\si(K\sm J) \ar@{=}[r]^{\id_{\si(K\!\sm\!J)}} \ar[dd]
&\si(K\sm J) \ar[dd] \\ \\
&& 0 & {\phantom{.}}0.
}
\end{gathered}
\label{aa4eq10}
\e
Here solid arrows `$\ra$' have already been defined, and
dashed arrows `$\dashra$' remain to be constructed. The
left hand square commutes by~\eq{aa4eq6}.

Now $\pi(C',K\sm J)\ci\io(A',C')=\pi(C',K\sm J)\ci\io(B',C')
\ci\io(A',B')=0$ as the middle column is exact. Since $\pi(C',K)$
is the cokernel of $\io(A',C')$, there exists a unique $\pi(K,
K\sm J):\si(K)\ra\si(K\sm J)$ with $\pi(K,K\sm J)\ci\pi(C',K)=
\pi(C',K\sm J)$. As $\pi(C',K\sm J)$ is surjective, $\pi(K,K\sm J)$
is surjective as in (iii). Thus in \eq{aa4eq10} the lower dashed
arrow exists, and the lower square commutes.

Suppose $f:D\ra\si(B')$ with $\pi(C',K)\ci\io(B',C')\ci f=0$. Then
there is a unique $h:D\ra\si(A')$ with $\io(A',C')\ci h=\io(B',C')
\ci f$, as $\io(A',C')$ is the kernel of $\pi(C',K)$. But then
$\io(B',C')\ci\io(A',B')\ci h=\io(B',C')\ci f$, so $\io(A',B')\ci h=f$
as $\io(B',C')$ is injective. Thus $\io(A',B')$ is the kernel
of $\pi(C',K)\ci\io(B',C')$. Similarly, $\pi(K,K\sm J)$ is the
cokernel of~$\pi(C',K)\ci\io(B',C')$.

Now apply Definition \ref{aa2def1}(iv) to $\pi(C',K)\ci\io(B',C')$.
As it has kernel $\io(A',B')$ and cokernel $\pi(K,K\sm J)$, and
$\pi(B',J)$ is the cokernel of $\io(A',B')$, this gives a unique
$\io(J,K):\si(J)\ra\si(K)$ with $\io(J,K)\ci\pi(B',J)=\pi(C',K)\ci
\io(B',C')$, such that $\io(J,K)$ is the kernel of $\pi(K,K\sm J)$.
Thus $\io(J,K)$ is injective, as in (ii), and in \eq{aa4eq10} the
upper dashed arrow exists, the upper right square commutes, and
the right hand column is exact, proving~(A).

We should also check that if $J,K$ are s-sets, the definition
above gives the same answer for $\io(J,K)$ as Step 1, and for
$\pi(K,K\sm J)$ as Step 3. If $J,K$ are s-sets then $A',J,K$ are
s-sets with $A'\cap J=A'\cap K=\emptyset$, so Step 2 gives $\si(B')
\cong\si(A')\op\si(J)$ and $\si(C')\cong\si(A')\op\si(K)$. Substituting
these into \eq{aa4eq10}, we find the definitions are consistent.
\medskip

\noindent{\bf Step 5.} Let $C$ be an s-set and $D,E$ f-sets with
$(C,D),(C,E),(D,E)\in\H_\sIp$, so that $C\supseteq D\supseteq E$. Apply
Step 4 with $J=D\sm E$ and $K=D$. This gives $C'$ which is the
largest s-set with $(C',D)\in\H_\sIp$, so $C\subseteq C'$. Therefore
\begin{equation*}
\pi(C,E)\!=\!\pi(C',E)\!\ci\!\io(C,C')\!=\!\pi(D,E)\!\ci\!
\pi(C',D)\!\ci\!\io(C',C)\!=\!\pi(D,E)\!\ci\!\pi(C,D),
\end{equation*}
using \eq{aa4eq9} for the first and third steps, and commutativity
of the bottom square in \eq{aa4eq10} for the second. Hence
\e
\pi(C,E)=\pi(D,E)\ci\pi(C,D)
\quad\text{for $C$ an s-set, as in (C).}
\label{aa4eq11}
\e

Suppose $J,K$ are s-sets and $L$ an f-set with $(J,K)\in\G_\sIp$ and
$(K,L)\in\H_\sIp$. Then $J\cap L$ is an s-set with $(J,J\cap L)\in\H_\sIp$
and $(J\cap L,L)\in\G_\sIp$. As in Step 4 with $J,K$ replaced by $J\cap
L,L$, define $A'=\{i\in I\sm L:l\npr i$ for all $l\in L\}$,
$B'=A'\cup(J\cap L)$ and $C'=A'\cup L$. Then $J\subseteq B'$ and
$K\subseteq C'$. Therefore
\begin{align*}
&\io(J\cap L,L)\ci\pi(J,J\cap L)=
\io(J\cap L,L)\ci\pi(B',J\cap L)\ci\io(J,B')=\\
&\pi(C',L)\ci\io(B',C')\ci\io(J,B')=
\pi(C',L)\ci\io(B,C')=\\
&\pi(C',L)\ci\io(K,C')\ci\io(J,K)=
\pi(K,L)\ci\io(J,K),
\end{align*}
using \eq{aa4eq9} at the first and fifth steps, commutativity
of the upper right square in \eq{aa4eq10} at the second, and
\eq{aa4eq6} at the third and fourth. This proves
\e
\pi(K,L)\ci\io(J,K)=\io(J\cap L,L)\ci\pi(J,J\cap L)
\quad\text{for $J,K$ s-sets, as in (D).}
\label{aa4eq12}
\e

\noindent{\bf Step 6.} Suppose $(J,K),(K,L)\in\G_\sIp$, and define
$D=\{i\in I:i\pr l$ for some $l\in L\}$, $A=D\sm L$, $B=A\cup J$,
and $C=A\cup K$. Then $A\subseteq B\subseteq C\subseteq D$ are s-sets,
with $J=B\sm A$, $K=C\sm A$, and $L=D\sm A$. Therefore
\begin{align*}
\io(K,L)\ci\io(J,K)\ci\pi(B,J)&=
\io(K,L)\ci\pi(C,K)\ci\io(B,C)=\\
\pi(D,L)\ci\io(C,D)\ci\io(B,C)&=
\pi(D,L)\ci\io(B,D)=\io(J,L)\ci\pi(B,J),
\end{align*}
using \eq{aa4eq12} at the first, second and fourth steps with
$J=B\cap K$, $K=C\cap L$ and $J=B\cap L$ respectively, and
\eq{aa4eq6} at the third. As $\pi(B,J)$ is surjective this
implies that $\io(K,L)\ci\io(J,K)=\io(J,L)$, proving~(B).

Similarly, suppose $(J,K),(K,L)\in\H_\sIp$, and define $D=\{i\in I:
i\pr j$ for some $j\in J\}$, $A=D\sm J$, $B=D\sm K$, and
$C=D\sm L$. Then $A\subseteq B\subseteq C\subseteq D$ are s-sets,
with $J=D\sm A$, $K=D\sm B$, and $L=D\sm C$. Therefore
\begin{equation*}
\pi(J,L)\!\ci\!\pi(D,J)\!=\!\pi(D,L)\!=\!\pi(K,L)\!\ci\!
\pi(D,K)\!=\!\pi(K,L)\!\ci\!\pi(J,K)\!\ci\!\pi(D,J),
\end{equation*}
using \eq{aa4eq11} three times. As $\pi(D,J)$ is surjective,
this proves (C). Applying (B), (C) with $J=K=L$ gives
$\io(J,J)=\pi(J,J)=\id_{\si(J)}$, as in (ii) and~(iii).
\medskip

\noindent{\bf Step 7.} Suppose $(J,K)\in\G_\sIp$ and $(K,L)\in\H_\sIp$,
and define $C=\{i\in I:i\pr k$ for some $k\in K\}$, $A=C\sm K$ and
$B=A\cup J$. Then $A\subseteq B\subseteq C$ are s-sets, with $J=B\sm A$
and $K=C\sm A$. Therefore
\begin{align*}
&\pi(K,L)\!\ci\!\io(J,K)\!\ci\!\pi(B,J)\!=\!
\pi(K,L)\!\ci\!\pi(C,K)\!\ci\!\io(B,C)\!=\!
\pi(C,L)\!\ci\!\io(B,C)\!=\\
&\io(J\cap L,L)\ci\pi(B,J\cap L)
=\io(J\cap L,L)\ci\pi(J,J\cap L)\ci\pi(B,J),
\end{align*}
using \eq{aa4eq12} at the first and third steps with $J=B\cap K$,
$B\cap L=J\cap L$ respectively, and (C) at the second and fourth.
As $\pi(B,J)$ is surjective this proves~(D).

Hence $(\si,\io,\pi)$ is an $(I,\pr)$-{\it configuration}, in
the sense of Definition \ref{aa4def}. It remains only to show that
$(\si,\io,\pi)$ is {\it unique up to canonical isomorphism} in $\A$.
At each stage in the construction the objects and morphisms were
determined either uniquely up to canonical isomorphism, or uniquely.
Thus, if $(\si,\io,\pi)$, $(\si',\io',\pi')$ both satisfy the
conditions of the theorem, one can go through the steps above and
construct a canonical isomorphism between them.
\end{proof}

Applying the theorem to the situation of \S\ref{aa3} gives:

\begin{cor} Let\/ $\A,X,I$ and\/ $\pr$ be as in Definition
\ref{aa3def1}. Then there exists an $(I,\pr)$-configuration
$(\si,\io,\pi)$ in $\A$ with\/ $\si(I)=X$ and\/ $\si(\{i\})
\cong S^i$ for $i\in I$, such that\/ $\io(J,I):\si(J)\ra X$
represents the subobject of\/ $X$ corresponding to $J$ under
the $1$-$1$ correspondence of Proposition \ref{aa3prop2} for
each s-set\/ $J\subseteq I$. This $(\si,\io,\pi)$ is unique
up to canonical isomorphism in~$\A$.
\label{aa4cor1}
\end{cor}

As the s-sets of $(\{1,\ldots,n\},\le)$ are $\{1,\ldots,j\}$
for $0\le j\le n$, the corresponding subobjects $S^J$ form a
{\it filtration} $0\!=\!A_0\!\subset\!\cdots\!\subset\!A_n\!=\!X$,
so Theorem \ref{aa4thm1} gives:

\begin{cor} Let\/ $0=A_0\subset A_1\subset\cdots\subset A_n=X$
be a filtration in an abelian category $\A$. Then there is a
$(\{1,\ldots,n\},\le)$-configuration $(\si,\io,\pi)$ in $\A$,
unique up to canonical isomorphism, such that\/ $\io(\{1,\ldots,j\},
\{1,\ldots,n\}):\si(\{1,\ldots,j\})\ra X$ represents $A_j\subset X$
for~$j=0,\ldots,n$.
\label{aa4cor2}
\end{cor}

This shows we can regard $(I,\pr)$-configurations as
{\it generalized filtrations}. Here is the converse to
Theorem~\ref{aa4thm1}.

\begin{thm} Let\/ $(\si,\io,\pi)$ be an $(I,\pr)$-configuration
in an abelian category $\A$. Define $X=\si(I)$, and let\/
$S^J\subset X$ be represented by $\io(J,I):\si(J)\ra X$ for each
s-set\/ $J\subseteq I$. Then the $S^J$ satisfy~\eq{aa4eq4}.
\label{aa4thm2}
\end{thm}

\begin{proof} The first two equations of \eq{aa4eq4} are obvious.
So suppose $A,B\subseteq I$ are s-sets. Definition \ref{aa2def5}
with $S=\si(A)$, $T=\si(B)$, $i=\io(A,I)$ and $j=\io(B,I)$ gives
$U\in\A$ and $a:U\ra\si(A)$, $b:U\ra\si(B)$ with $i\ci a=j\ci b$,
such that $i\ci a:U\ra X$ represents $S^A\cap S^B$. As
$i\ci\io(A\cap B,A)=\io(A\cap B,I)=j\ci(A\cap B,B)$, by exactness
in \eq{aa2eq1} there is a unique $h:\si(A\cap B)\ra U$ with
\begin{equation*}
\bigl(\io_{\si(A)}\ci a-\io_\si(B)\ci b\bigr)\ci h=
\io_{\si(A)}\ci\io(A\cap B,A)-\io_\si(B)\ci\io(A\cap B,B).
\end{equation*}
Composing $\pi_{\si(A)},\pi_{\si(B)}$ gives $\io(A\cap B,A)=
a\ci h$ and $\io(A\cap B,B)=b\ci h$. Now
\begin{equation*}
\io(A\cup B,I)\ci\io(A,A\cup B)\ci a=i\ci a=j\ci b=
\io(A\cup B,I)\ci\io(B,A\cup B)\ci b,
\end{equation*}
by Definition \ref{aa4def}(B). Thus $\io(A,A\cup B)\ci a=
\io(B,A\cup B)\ci b$, as $\io(A\cup B,I)$ is injective. Hence
\begin{align*}
\pi(A,A\sm B)\ci a=\io(A\sm B,A\sm B)\ci\pi(A,A\sm B)\ci a&=\\
\pi(A\cup B,A\sm B)\ci\io(A,A\cup B)\ci a=
\pi(A\cup B,A\sm B)\ci\io(B,A\cup B)\ci b&=0,
\end{align*}
using Definition \ref{aa4def}(D) at the second step
and exactness in (A) at the fourth.

But $\io(A\cap B,A)$ is the kernel of $\pi(A,A\sm B)$, so there
is a unique $h':U\ra\si(A\cap B)$ with $a=\io(A\cap B,A)\ci h'$.
As $\io(A\cap B,A)=a\ci h$ and $a,\io(A\cap B,A)$ are injective we
see that $h,h'$ are {\it inverse}, so $h$ is {\it invertible}. This
implies that
\begin{equation*}
\io(A\cap B,I)=\io(A,I)\ci\io(A\cap B,A)=\io(A,I)\ci a\ci h:
\si(A\cap B)\ra X
\end{equation*}
represents $S^A\cap S^B$, so that $S^{A\cap B}=S^A\cap S^B$.
We prove $S^{A\cup B}=S^A+S^B$ in a similar way.
\end{proof}

Combining Theorems \ref{aa4thm1} and \ref{aa4thm2} we deduce:

\begin{cor} For $(I,\pr)$ a finite poset and\/ $\A$ an abelian
category, there is an equivalence of categories between the groupoid
of\/ $(I,\pr)$-configurations in $\A$, and the groupoid of collections
$(X\in\A,$ subobjects $S^J\subset X$ for s-sets $J\subset I)$
satisfying \eq{aa4eq4}, with the obvious notion of isomorphism.
\label{aa4cor3}
\end{cor}

Finally, for an $(I,\pr)$-configuration $(\si,\io,\pi)$ we show
how the classes $[\si(J)]$ in the {\it Grothendieck group}
$K_0(\A)$ are related.

\begin{prop} Let\/ $(\si,\io,\pi)$ be an $(I,\pr)$-configuration
in an abelian category $\A$. Then there exists a unique map
$\ka:I\ra K_0(\A)$ such that\/ $[\si(J)]=\sum_{j\in J}\ka(j)$
in $K_0(\A)$ for all f-sets~$J\subseteq I$.
\label{aa4prop}
\end{prop}

\begin{proof} Combining Definitions \ref{aa2def2} and
\ref{aa3def1}(A) shows that
\e
\bigl[\si(K)\bigr]=\bigl[\si(J)\bigr]+\bigl[\si(K\sm J)\bigr]
\quad\text{for all $(J,K)\in\G_\sIp$.}
\label{aa4eq13}
\e
Define $\ka:I\ra K_0(\A)$ by $\ka(i)=[\si(\{i\})]$. As
$\{i\}\in\F_\sIp$ for all $i\in I$ this is unique and
well-defined. Suppose $K\in\F_\sIp$ with $\md{K}\ge 1$. Let
$j\in K$ be $\pr$-{\it minimal}. Then $(\{j\},K)\in\G_\sIp$,
so \eq{aa4eq13} gives $[\si(K)]=\ka(j)+[\si(K\sm\{j\})]$.
Thus $[\si(J)]=\sum_{j\in J}\ka(j)$ for all $J\in\F_\sIp$
by induction on~$\md{J}$.
\end{proof}

\section{New $(I,\pr)$-configurations from old}
\label{aa5}

Let $(\si,\io,\pi)$ be an $(I,\pr)$-configuration in an abelian
or exact category $\A$. Then we can derive $(K,\tl)$-configurations
$(\ti\si,\ti\io,\ti\pi)$ in $\A$ from $(\si,\io,\pi)$ for other,
simpler finite posets $(K,\tl)$, by forgetting some of the
information in $(\si,\io,\pi)$. The next two definitions give two
ways to do this. We use the notation of \S\ref{aa3} and~\S\ref{aa4}.

\begin{dfn} Let $(I,\pr)$ be a finite poset and $J\in\F_\sIp$. Then
$(J,\pr)$ is also a finite poset, and $K\subseteq J$ is an f-set in
$(J,\pr)$ if and only if it is an f-set in $(I,\pr)$. Hence
$\F_\sJp\!\subseteq\!\F_\sIp$. We also have $\G_\sJp\!=\!\G_\sIp\!
\cap\!(\F_\sJp\!\t\!\F_\sJp)$ and $\H_\sJp\!=\!\H_\sIp\!\cap\!
(\F_\sJp\!\t\!\F_\sJp)$, so that $\G_\sJp\!\subseteq\!\G_\sIp$
and~$\H_\sJp\!\subseteq\!\H_\sIp$.

Let $(\si,\io,\pi)$ be an $(I,\pr)$-configuration in an abelian
or exact category $\A$, and define $\si':\F_\sJp\ra\Obj(\A)$,
$\io':\G_\sJp\ra\Mor(\A)$ and $\pi':\H_\sJp\ra\Mor(\A)$ by
$\si'=\si\vert_{\F_\sJp}$, $\io'=\io\vert_{\G_\sJp}$ and
$\pi'=\pi\vert_{\H_\sJp}$. Then (A)--(D) of Definition
\ref{aa4def} for $(\si,\io,\pi)$ imply (A)--(D) for
$(\si',\io',\pi')$, so $(\si',\io',\pi')$ is a
$(J,\pr)$-{\it configuration} in~$\A$. We call
$(\si',\io',\pi')$ a {\it subconfiguration}
of~$(\si,\io,\pi)$.
\label{aa5def1}
\end{dfn}

\begin{dfn} Let $(I,\pr),(K,\tl)$ be finite posets, and
$\phi:I\ra K$ be surjective with $\phi(i)\tl\phi(j)$ when
$i,j\in I$ with $i\pr j$. If $J\subseteq K$ is an f-set
in $K$ then $\phi^{-1}(J)\subseteq I$ is an f-set in $I$. Hence
$\phi^*(\F_\sKt)\!\subseteq\!\F_\sIp$, where $\phi^*$ pulls back
subsets of $K$ to subsets of $I$. Similarly, $\phi^*(\G_\sKt)\!
\subseteq\!\G_\sIp$ and~$\phi^*(\H_\sKt)\!\subseteq\!\H_\sIp$.

Let $(\si,\io,\pi)$ be an $(I,\pr)$-configuration in an abelian
or exact category $\A$, and define $\ti\si:\F_\sKt\ra\Obj(\A)$,
$\ti\io:\G_\sKt\ra\Mor(\A)$ and $\ti\pi:\H_\sKt\ra\Mor(\A)$ by
$\ti\si(A)=\si\bigl(\phi^{-1}(A)\bigr)$, $\ti\io(A,B)=
\io\bigl(\phi^{-1}(A),\phi^{-1}(B)\bigr)$, and
$\ti\pi(A,B)=\pi\bigl(\phi^{-1}(A),\phi^{-1}(B)\bigr)$. Then
$(\ti\si,\ti\io,\ti\pi)$ is a $(K,\tl)$-{\it configuration}
in $\A$, the {\it quotient configuration} of $(\si,\io,\pi)$.
We also call $(\si,\io,\pi)$ a {\it refinement\/}
of~$(\ti\si,\ti\io,\ti\pi)$.
\label{aa5def2}
\end{dfn}

Compositions of these constructions all behave in the
obvious ways. Next we explain a method to {\it glue two
configurations $(\si',\io',\pi')$, $(\ti\si,\ti\io,\ti\pi)$
together}, to get $(\si,\io,\pi)$ containing $(\si',\io',\pi')$
as a {\it subconfiguration}, and $(\ti\si,\ti\io,\ti\pi)$ as a
{\it quotient configuration}. Consider the following situation.

\begin{dfn} Let $(J,\ls),(K,\tl)$ be finite posets and $L\subset K$
an f-set, with $J\cap(K\sm L)=\emptyset$. Suppose $\psi:J\ra L$
is a surjective map with $\psi(i)\tl\psi(j)$ when $i,j\in J$ with
$i\ls j$. Set $I=J\cup(K\sm L)$, and define a binary relation
$\pr$ on $I$ by
\begin{equation*}
i\pr j\quad\text{for $i,j\in I$ if}\quad
\begin{cases}
i\ls j, & i,j\in J, \\
i\tl j, & i,j\in K\sm L, \\
\psi(i)\tl j, & i\in J,\quad j\in K\sm L, \\
i\tl\psi(j), & i\in K\sm L,\quad j\in J.
\end{cases}
\end{equation*}
One can show $\pr$ is a partial order on $I$, and $J\subseteq I$
an f-set in $(I,\pr)$. The restriction of $\pr$ to $J$ is $\ls$.
Define $\phi:I\ra K$ by $\phi(i)=\psi(i)$ if $i\in J$ and $\phi(i)=i$
if $i\in K\sm L$. Then $\phi$ is surjective, with $\phi(i)\tl\phi(j)$
when $i,j\in I$ with~$i\pr j$.
\label{aa5def3}
\end{dfn}

An $(I,\pr)$-configuration gives the same
$(L,\tl)$-configuration in two ways.

\begin{lem} In the situation of Definition \ref{aa5def3}, suppose
$(\si,\io,\pi)$ is an $(I,\pr)$-configuration in an abelian or exact
category. Let\/ $(\si',\io',\pi')$ be its $(J,\ls)$-subconfiguration,
and\/ $(\ti\si,\ti\io,\ti\pi)$ its quotient\/ $(K,\tl)$-configuration
from $\phi$. Let\/ $(\hat\si,\hat\io,\hat\pi)$ be the quotient\/
$(L,\tl)$-configuration from $(\si',\io',\pi')$ and\/ $\psi$,
and\/ $(\check\si,\check\io,\check\pi)$ the
$(L,\tl)$-subconfiguration from $(\ti\si,\ti\io,\ti\pi)$.
Then~$(\hat\si,\hat\io,\hat\pi)=(\check\si,\check\io,\check\pi)$.
\label{aa5lem1}
\end{lem}

Our third construction is a kind of converse to Lemma \ref{aa5lem1}.
In categorical notation, the last part says there is an {\it equivalence
of categories} between the category of $(I,\pr)$-configurations, and the
{\it fibred product\/} of the categories of $(J,\ls)$-configurations and
$(K,\tl)$-configurations over $(L,\tl)$-configurations.

\begin{thm} In the situation of Definition \ref{aa5def3}, let\/ $\A$ be
an abelian or exact category, $(\si',\io',\pi')$ a $(J,\ls)$-configuration
in $\A$, and\/ $(\ti\si,\ti\io,\ti\pi)$ a $(K,\tl)$-configuration in
$\A$. Define $(\hat\si,\hat\io,\hat\pi)$ to be the quotient\/
$(L,\tl)$-configuration from $(\si',\io',\pi')$ and\/ $\psi$, and\/
$(\check\si,\check\io,\check\pi)$ to be the $(L,\tl)$-subconfiguration
from~$(\ti\si,\ti\io,\ti\pi)$.

Suppose $(\hat\si,\hat\io,\hat\pi)\!=\!(\check\si,\check\io,\check\pi)$.
Then there exists an $(I,\pr)$-configuration $(\si,\io,\pi)$ in $\A$,
unique up to canonical isomorphism, such that\/ $(\si',\io',\pi')$
is its $(J,\ls)$-subconfiguration, and\/ $(\ti\si,\ti\io,\ti\pi)$
its quotient\/ $(K,\tl)$-configuration from~$\phi$.

More generally, given $\al:(\check\si,\check\io,\check\pi)
\,\smash{{\buildrel\cong\over\longra}}\,(\hat\si,\hat\io,
\hat\pi)$ there is an $(I,\pr)$-configuration $(\si,\io,\pi)$
in $\A$, unique up to canonical isomorphism, with\/
$(J,\ls)$-subconfiguration isomorphic to $(\si',\io',\pi')$,
and quotient\/ $(K,\tl)$-configuration from $\phi$ isomorphic
to $(\ti\si,\ti\io,\ti\pi)$, such that the equality in Lemma
\ref{aa5lem1} corresponds to~$\al$.
\label{aa5thm1}
\end{thm}

\begin{proof} Let $\A$ be an {\it abelian category}. We prove
the first part in five steps:
\begin{list}{}{
\setlength{\itemsep}{1pt}
\setlength{\parsep}{1pt}
\setlength{\labelwidth}{50pt}
\setlength{\leftmargin}{50pt}}
\item[{\bf Step 1.}] Characterize $(I,\pr)$ s-sets.
\item[{\bf Step 2.}] Define $\si(B)$ for all $(I,\pr)$ s-sets,
and some morphisms~$\io(B,C)$.
\item[{\bf Step 3.}] Define $\io(B,B')$ for all $(I,\pr)$
s-sets $B\subseteq B'\subseteq I$, and prove~$\io=\io\ci\io$.
\item[{\bf Step 4.}] Let $S^B$ be the subobject represented
by $\io(B,I):\si(B)\ra\si(I)=X$ for all $(I,\pr)$ s-sets $B$.
Show that the $S^B$ satisfy~\eq{aa4eq4}.
\item[{\bf Step 5.}] Apply Theorem \ref{aa4thm1} to construct
$(\si,\io,\pi)$, and complete the proof.
\end{list}

\noindent{\bf Step 1.} The proof of the next lemma is elementary,
and left as an exercise.

\begin{lem} In the situation above, let\/ $B\subseteq I$ be an
s-set in $(I,\pr)$. Define $P=\bigl\{k\in K:$ if\/ $i\in I$ and\/
$\phi(i)\tl k$, then $i\in B\bigr\},$ and\/ $R=\bigl\{k\in K:k\tl
\phi(i)$ for some $i\in B\bigr\}$. Define $A=\phi^{-1}(P)$ and\/
$C=\phi^{-1}(R)$. Then $P\subseteq R$ are $(K,\tl)$ s-sets, and\/
$A\subseteq B\subseteq C$ are $(I,\pr)$ s-sets, with\/ $P\sm L=
R\sm L=B\sm J$. Define $D=A\cap J$, $E=B\cap J$ and\/ $F=C\cap J$.
Then $D\subseteq E\subseteq F$ are $(J,\ls)$ s-sets, with
$(P,R)\in\G_\sKt$ and\/ $(D,E),(D,F),(E,F)\in\G_\sJl$. Define
$U=P\cap L$ and\/ $W=R\cap L$. Then $U\subseteq W$ are $(L,\tl)$
s-sets with\/ $\phi^{-1}(U)=\psi^{-1}(U)=D$ and\/
$\phi^{-1}(W)=\psi^{-1}(W)=F$. Hence
\e
\begin{gathered}
\si'(D)\!=\!\hat\si(U)\!=\!\check\si(U)\!=\!\ti\si(U),
\;\>\text{and similarly}\;\>
\si'(F\sm D)=\ti\si(W\sm U),\\
\si'(F)=\ti\si(W),\;\>
\io'(D,F)=\ti\io(U,W),\;\>
\pi'(F,F\sm D)=\ti\pi(W,W\sm U).
\end{gathered}
\label{aa5eq1}
\e
\label{aa5lem2}
\end{lem}

\noindent Here $P,R$ are the largest, smallest $(K,\tl)$ s-sets
with~$\phi^{-1}(P)\!\subseteq\!B\!\subseteq\!\phi^{-1}(R)$.
\medskip

\noindent{\bf Step 2.} Let $B$ be an $(I,\pr)$ s-set, and use
the notation of Lemma \ref{aa5lem2}. As $R\sm P\!=\!W\sm U$
we have $\ti\si(R\sm P)\!=\!\ti\si(W\sm U)\!=\!\si'(F\sm D)$
by \eq{aa5eq1}. Consider
\e
\pi'(F\sm D,F\sm E)\ci\ti\pi(R,R\sm P):\ti\si(R)\ra\si'(F\sm E).
\label{aa5eq2}
\e
Choose $\si(B)\in\A$ and $\io(B,C):\si(B)\ra\si(C)=\ti\si(R)$
to be a kernel for \eq{aa5eq2}. If $B=\phi^{-1}(Q)$ for some
$(K,\tl)$ s-set $Q$ then $B=C$ and \eq{aa5eq2} is zero, and
we choose $\si(B)=\ti\si(R)$ and $\io(B,C)=\id_{\si(B)}$.
Define~$\io(B,I)=\ti\io(R,K)\ci\io(B,C)$.
\medskip

\noindent{\bf Step 3.} Let $B\subseteq B'$ be $(I,\pr)$ s-sets.
Use the notation of Lemma \ref{aa5lem2} for $B$, and $P',\ldots,F'$
for $B'$. Then $P\subseteq P'$, $R\subseteq R'$, and so on. We have
\begin{align*}
\pi'(F'\sm D',F'\sm E')\ci\ti\pi(R',R'\sm P')\ci
\ti\io(R,R')\ci\io(B,C)&=\\
\pi'(F'\sm D',F'\sm E')\ci\ti\io(R\sm P',R'\sm P')\ci
\ti\pi(R,R\sm P')\ci\io(B,C)&=\\
\pi'(F'\!\sm\!D',F'\!\sm\!E')\!\ci\!\io'(F\!\sm\!D',F'\!\sm\!D')
\!\ci\!\ti\pi(R\!\sm\!P,R\!\sm\!P')\!\ci\!\ti\pi(R,R\!\sm\!P)
\!\ci\!\io(B,C)&=\\
\io'(F\!\sm\!E',F'\!\sm\!E')\ci\pi'(F\!\sm\!D',F\!\sm\!E')\!\ci\!
\pi'(F\!\sm\!D,F\!\sm\!D')\!\ci\!\ti\pi(R,R\!\sm\!P)\!\ci\!\io(B,C)&=\\
\io'(F\!\sm\!E',F'\!\sm\!E')\!\ci\!
\pi'(F\!\sm\!E,F\!\sm\!E')\!\ci\!
\pi'(F\!\sm\!D,F\!\sm\!E)\!\ci\!
\ti\pi(R,R\!\sm\!P)\ci\io(B,C)&=0,
\end{align*}
using Definition \ref{aa4def}(C), (D), and the definition
of~$\io(B,C)$.

Thus, as $\io(B',C')$ is the kernel of $\pi'(F'\sm D',F'\sm E')\ci
\ti\pi(R',R'\sm P')$, there exists a unique $\io(B,B')$ with
$\ti\io(R,R')\ci\io(B,C)=\io(B',C')\ci\io(B,B')$. Hence
\begin{align*}
\io(B',I)\ci\io(B,B')=\ti\io(R',K)\ci\io(B',C')\ci\io(B,B')&=\\
\ti\io(R',K)\ci\ti\io(R,R')\ci\io(B,C)=\ti\io(R,K)\ci\io(B,C)&=\io(B,I).
\end{align*}
The proof of \eq{aa4eq6} from \eq{aa4eq5} then gives
\e
\io(B,B'')=\io(B',B'')\ci\io(B,B')\quad
\text{when $B\subseteq B'\subseteq B''$ are $(I,\pr)$ s-sets.}
\label{aa5eq3}
\e

\noindent{\bf Step 4.} Set $X=\si(I)=\ti\si(K)$, and for each
$(I,\pr)$ s-set $B$ let $S^B\subset X$ be subobject represented
by $\io(B,I):\si(B)\ra\si(I)=X$. We must prove that these $S^B$
satisfy \eq{aa4eq4}. The first two equations of \eq{aa4eq4}
are immediate. Let $B',B''$ be $(I,\pr)$ s-sets, and
$B=B'\cap B''$. We shall show that~$S^B=S^{B'}\cap S^{B''}$.

Use the notation of Lemma \ref{aa5lem2} for $B$, and
$P',R',\ldots$ for $B'$ and $P'',R'',\ldots$ for $B''$ in the
obvious way. Apply Definition \ref{aa2def5} with $i=\io(B',I)$
and $j=\io(B'',I)$, giving $U,V\in\A$ and morphisms $a,b,c,d,e$
with $i\ci a=j\ci b$, $i=e\ci c$ and $j=e\ci d$, such that
$i\ci a:U\ra X$ represents~$S^{B'}\cap S^{B''}$.

Set $\hat C=C'\cap C''$, $\hat R=R'\cap R''$ and $\hat F=F'\cap F''$,
so that $C\subseteq\hat C$, $R\subseteq\hat R$, $F\subseteq\hat F$,
$\hat C=\phi^{-1}(\hat R)$, and $\hat F=\hat C\cap J$. Define
$\hat a:U\ra\si(C')$ and $\hat b:U\ra\si(C'')$ by $\hat a=\io(B',C')
\ci a$ and $\hat b=\io(B'',C'')\ci b$. Then
\begin{align*}
\ti\io(R',K)\ci\hat a&=\ti\io(R',K)\ci\io(B',C')\ci a=
\io(B',I)\ci a=i\ci a=\\
j\ci b&=\io(B'',I)\ci b=\ti\io(R'',K)\ci\io(B'',C'')\ci b=
\ti\io(R'',K)\ci\hat b.
\end{align*}
As $(\ti\si,\ti\io,\ti\pi)$ is a configuration we see that
$S^{\hat C}=S^{C'}\cap S^{C''}$ by Theorem \ref{aa4thm2}.
Thus there is a unique $\hat h:U\ra\ti\si(\hat R)$ with
$\hat a=\ti\io(\hat R,R')\ci\hat h$ and~$\hat b=
\ti\io(\hat R,R'')\ci\hat h$.

As $\io(B',C')$ is the kernel of $\pi'(F'\sm D',F'\sm E')\ci
\ti\pi(R',R'\sm P')$, we have
\begin{align*}
0=\pi'(F'\sm D',F'\sm E')\ci\ti\pi(R',R'\sm P')\ci\io(B',C')\ci a&=\\
\pi'(F'\sm D',F'\sm E')\ci\ti\pi(R',R'\sm P')\ci\hat a&=\\
\pi'(F'\sm D',F'\sm E')\ci\ti\pi(R',R'\sm P')\ci\ti\io(\hat R,R')
\ci\hat h=\cdots&=\\
\io'(\hat F\!\sm\!E',F'\!\sm\!E')\!\ci\!\pi'(\hat F\!\sm\!E,
\hat F\!\sm\!E')\!\ci\!\pi'(F\!\sm\!D,\hat F\!\sm\!E)
\!\ci\!\ti\pi(\hat R,\hat R\!\sm\!P)\!\ci\!\hat h&,
\end{align*}
by Definition \ref{aa4def}(C), (D). Since
$\io'(\hat F\sm E',F'\sm E')$ is injective this gives
\e
\begin{split}
\pi'(\hat F\sm E,\hat F\sm E')\ci\pi'(\hat F\sm D,\hat F\sm E)\ci
\ti\pi(\hat R,\hat R\sm P)\ci\hat h&=0,\\
\text{and}\quad
\pi'(\hat F\sm E,\hat F\sm E'')\ci\pi'(\hat F\sm D,\hat F\sm E)\ci
\ti\pi(\hat R,\hat R\sm P)\ci\hat h&=0,
\end{split}
\label{aa5eq4}
\e
proving the second equation in the same way using~$B'',C'',\ldots$.

As $(\si',\io',\pi')$ is a configuration and $E=E'\cap E''$, one
can show that
\begin{equation*}
\io_{\si'(\hat F\sm E')}\ci\pi'(\hat F\sm E,\hat F\sm E')+
\io_{\si'(\hat F\sm E'')}\ci\pi'(\hat F\sm E,\hat F\sm E'')
\end{equation*}
is an injective morphism $\si'(\hat F\sm E)\ra\si'(\hat F\sm E')
\op\si'(\hat F\sm E'')$. Therefore
\e
\pi'(\hat F\sm D,\hat F\sm E)\ci\ti\pi(\hat R,\hat R\sm P)\ci\hat h=0,
\label{aa5eq5}
\e
by \eq{aa5eq4}. Composing \eq{aa5eq5} with $\pi'(\hat F\sm E,
\hat F\sm F)$ and using Definition \ref{aa4def}(C) shows that
$\ti\pi(\hat R,\hat R\sm R)\ci\hat h=0$. But $\ti\io(R,\hat R)$
is the kernel of $\ti\pi(\hat R,\hat R\sm R)$, so
$\hat h=\ti\io(R,\hat R)\ci\ti h$ for some
unique~$\ti h:U\ra\ti\si(R)=\si(C)$.

Substituting $\hat h=\ti\io(R,\hat R)\ci\ti h$ into
\eq{aa5eq5} and using Definition \ref{aa4def}(D) gives
\begin{equation*}
\io'(F\sm E,\hat F\sm E)\ci\pi'(F\sm D,F\sm E)\ci
\ti\pi(R,R\sm P)\ci\ti h=0.
\end{equation*}
Hence $\pi'(F\sm D,F\sm E)\ci\ti\pi(R,R\sm P)\ci\ti h=0$, as
$\io'(F\sm E,\hat F\sm E)$ is injective. Thus, as $\io(B,C)$ is
the kernel of \eq{aa5eq2}, there is a unique $h:U\ra\si(B)$
with $\ti h=\io(B,C)\ci h$. Then
\begin{gather*}
\io(B',C')\ci a=\hat a=\ti\io(\hat R,R')\ci\hat h=
\ti\io(\hat R,R')\ci\ti\io(R,\hat R)\ci\ti h=\\
\ti\io(R,R')\ci\io(B,C)\ci h=
\io(B,C')\ci h=\io(B',C')\ci\io(B,B')\ci h
\end{gather*}
by \eq{aa5eq3}, so $a=\io(B,B')\ci h$ as $\io(B',C')$ is
injective, and similarly~$b=\io(B,B'')\ci h$.

Recall the definition of $i,j,U,V,a,\ldots,e$ above. By
\eq{aa5eq3} we have
\begin{align*}
e\ci c\ci\io(B,B')&=i\ci\io(B,B')=\io(B',I)\ci\io(B,B')=\io(B,I)=\\
\io(B'',I)\ci\io(B,B'')&=j\ci\io(B,B'')=e\ci d\ci\io(B,B'').
\end{align*}
Since $e$ is injective this gives $c\ci\io(B,B')=d\ci\io(B,B'')$,
and hence
\begin{equation*}
\bigl(c\ci\pi_{\si(B')}\op d\ci\pi_{\si(B'')}\bigr)
\ci\bigl(\io_{\si(B')}\ci\io(B,B')-
\io_{\si(B'')}\ci\io(B,B'')\bigr)=0,
\end{equation*}
factoring via $\si(B')\!\op\!\si(B'')$. So by \eq{aa2eq1}
there is a unique $m:\si(B)\!\ra\!U$ with
\begin{equation*}
\io_{\si(B')}\ci\io(B,B')-\io_{\si(B'')}\ci\io(B,B'')=
\bigl(\io_{\si(B')}\ci a-\io_{\si(B'')}\ci b\bigr)\ci m.
\end{equation*}
Composing with $\pi_{\si(B')}$ gives $\io(B,B')=a\ci m$.
As $a=\io(B,B')\ci h$ and $a,\io(B,B')$ are injective, we
see $m$ and $h$ are inverse, so $h$ is an isomorphism.

Since $S^{B'}\cap S^{B''}$ is represented by $\io(B',I)\ci a=
\io(B',I)\ci\io(B,B')\ci h=\io(B,I)\ci h$ and $S^B$ by $\io(B,I)$,
this proves that $S^B=S^{B'}\cap S^{B''}$ for all $(I,\pr)$ s-sets
$B',B''$ and $B=B'\cap B''$. A similar proof shows that
$S^B=S^{B'}+S^{B''}$ when $B=B'\cup B''$. Hence the $S^B$
satisfy~\eq{aa4eq4}.
\medskip

\noindent{\bf Step 5.} Theorem \ref{aa4thm1} now constructs an
$(I,\pr)$-configuration $(\si,\io,\pi)$, unique up to canonical
isomorphism, from the $S^B$. It follows from the construction
of the $S^B$ that the $(J,\ls)$-subconfiguration of $(\si,\io,\pi)$
is canonically isomorphic to $(\si',\io',\pi')$, and the quotient
$(K,\tl)$-configuration from $\phi$ is canonically isomorphic to
$(\ti\si,\ti\io,\ti\pi)$. It is not difficult to see that we can
choose $(\si,\io,\pi)$ so that these sub- and quotient configurations
are equal to $(\si',\io',\pi')$ and~$(\ti\si,\ti\io,\ti\pi)$.

For the last part, define $\be:\F_\sKt\!\ra\!\Mor(\A)$ by
$\be(A)\!=\!\al(A)$ if $A\!\in\!\F_\sLt$, and $\be(A)\!=\!\id_{\ti\si(A)}$
if $A\!\notin\!\F_\sLt$. Define a $(K,\tl)$-configuration
$(\dot\si,\dot\io,\dot\pi)$ by
\begin{equation*}
\dot\si(A)=\begin{cases} \hat\si(A) & A\in\F_\sLt, \\
\ti\si(A) & A\notin\F_\sLt, \end{cases} \;\>\text{and}\;\>
\begin{aligned}
\dot\io(A,B)&=\be(B)^{-1}\ci\ti\io(A,B)\ci\be(A),\\
\dot\pi(A,B)&=\be(B)^{-1}\ci\ti\pi(A,B)\ci\be(A).
\end{aligned}
\end{equation*}
Then $\be:(\ti\si,\ti\io,\ti\pi)\ra(\dot\si,\dot\io,\dot\pi)$
is an isomorphism. The $(L,\tl)$-subconfiguration of
$(\dot\si,\dot\io,\dot\pi)$ is $(\hat\si,\hat\io,\hat\pi)$,
so we may apply the first part with $(\ti\si,\ti\io,\ti\pi)$
replaced by $(\dot\si,\dot\io,\dot\pi)$. This proves Theorem
\ref{aa5thm1} when $\A$ is an {\it abelian category}.

If $\A$ is only an {\it exact category} we have more work to
do, as Steps 2,5 above involve choosing kernels and
cokernels, which may not exist in $\A$. So suppose $\A$
is an exact category, contained in an abelian category
$\skew8\hat\A$ as in \S\ref{aa21}. In the situation of the
first part of the theorem, the proof above yields an
$(I,\pr)$-configuration $(\si,\io,\pi)$ in $\skew8\hat\A$
with the properties we want. We must show $\si(J)\in\A$
for all $J\in\F_\sIp$, so that $(\si,\io,\pi)$ is a
configuration in~$\A$.

As $(\si',\io',\pi')$ is the $(J,\ls)$-subconfiguration of
$(\si,\io,\pi)$ we have $\si(\{i\})=\si'(\{i\})\in\A$ for
$i\in J$. And as $(\ti\si,\ti\io,\ti\pi)$ is the quotient
$(K,\tl)$-configuration of $(\si,\io,\pi)$ from $\phi$ we
have $\si(\{i\})=\ti\si(\{i\})\in\Obj(\A)$ for
$i\in I\sm J=K\sm L$. Hence $\si(\{i\})\in\A$ for all
$i\in I$. Also $\si(\emptyset)=0\in\A$.

Suppose by induction that $\si(A)\in\A$ for all
$A\in\F_\sIp$ with $\md{A}\le k$, for $1\le k<\md{I}$. Let
$B\in\F_\sIp$ with $\md{B}=k+1$, let $i$ be $\pr$-maximal
in $B$, and set $A=B\sm\{i\}$. Then $(A,B)\in\G_\sIp$, so
\eq{aa4eq1} gives a short exact sequence $0\ra\si(A)\ra
\si(B)\ra\si(\{i\})\ra 0$. Now $\si(A)\in\A$ by induction,
$\si(\{i\})\in\A$ from above, and $\A$ is closed under
extensions in $\skew8\hat\A$, so $\si(B)\in\A$. Thus by
induction $\si(A)\in\A$ for all $A\in\F_\sIp$, and
$(\si,\io,\pi)$ is an $(I,\pr)$-configuration in $\A$.
This proves the first part for $\A$ an exact category,
and the last part follows as above.
\end{proof}

The case when $L=\{l\}$ is one point will be particularly useful.

\begin{dfn} Let $(J,\ls)$ and $(K,\tl)$ be nonempty finite posets
with $J\cap K=\emptyset$, and $l\in K$. Set $I=J\cup(K\sm\{l\})$,
and define a partial order $\pr$ on $I$ by
\begin{equation*}
i\pr j\quad\text{for $i,j\in I$ if}\quad
\begin{cases}
i\ls j, & i,j\in J, \\
i\tl j, & i,j\in K\sm\{l\}, \\
l\tl j, & i\in J,\quad j\in K\sm\{l\}, \\
i\tl l, & i\in K\sm\{l\},\quad j\in J,
\end{cases}
\end{equation*}
and a surjective map $\phi:I\ra K$ by $\phi(i)=l$ if $i\in J$,
and $\phi(i)=i$ if~$i\!\in\!K\!\sm\!\{l\}$.

Let $\A$ be an abelian or exact category, $(\si',\io',\pi')$ a
$(J,\ls)$-configuration in $\A$, and $(\ti\si,\ti\io,\ti\pi)$
a $(K,\tl)$-configuration in $\A$ with $\si'(J)=\ti\si(\{l\})$.
Then by Theorem \ref{aa5thm1} there exists an
$(I,\pr)$-configuration $(\si,\io,\pi)$ in $\A$, unique up to
canonical isomorphism, such that $(\si',\io',\pi')$ is its
$(J,\ls)$-subconfiguration, and $(\ti\si,\ti\io,\ti\pi)$ its
quotient $(K,\tl)$-configuration from $\phi$. We call
$(\si,\io,\pi)$ the {\it substitution of\/ $(\si',\io',\pi')$
into}~$(\ti\si,\ti\io,\ti\pi)$.
\label{aa5def4}
\end{dfn}

\section{Improvements and best configurations}
\label{aa6}

We now study {\it quotient configurations} from $(I,\pr)$,
$(K,\tl)$ when $\phi:I\ra K$ is a {\it bijection}. So we
identify $I,K$ and regard $\pr,\tl$ as two partial orders on~$I$.

\begin{dfn} Let $I$ be a finite set and $\pr,\tl$ partial orders
on $I$ such that if $i\pr j$ then $i\tl j$ for $i,j\in I$. Then
we say that $\tl$ {\it dominates} $\pr$, and $\tl$ {\it strictly
dominates} $\pr$ if $\pr,\tl$ are distinct. Let $s$ be the number
of pairs $(i,j)\in I\t I$ with $i\tl j$ but $i\npr j$. Then we
say that $\tl$ {\it dominates $\pr$ by $s$ steps}. Also
\e
\F_\sIt\subseteq\F_\sIp,\quad
\G_\sIt\subseteq\G_\sIp\quad\text{and}\quad
\H_\sIt\subseteq\H_\sIp.
\label{aa6eq1}
\e

For each $(I,\pr)$-configuration $(\si,\io,\pi)$ in an abelian
or exact category we have a quotient $(I,\tl)$-configuration
$(\ti\si,\ti\io,\ti\pi)$, as in Definition \ref{aa5def2}
with $\phi=\id:I\ra I$. We call $(\si,\io,\pi)$ an
{\it improvement\/} or an $(I,\pr)$-{\it improvement} of
$(\ti\si,\ti\io,\ti\pi)$, and a {\it strict improvement\/} if
$\pr\ne\tl$. We call an $(I,\tl)$-configuration $(\ti\si,\ti\io,
\ti\pi)$ {\it best\/} if there exists no strict improvement
$(\si,\io,\pi)$ of $(\ti\si,\ti\io,\ti\pi)$. Improvements are
a special kind of {\it refinement}, in the sense of
Definition~\ref{aa5def2}.
\label{aa6def1}
\end{dfn}

Our first result is simple. An $(I,\tl)$-configuration $(\si,\io,\pi)$
cannot have an infinite sequence of strict improvements, as $I$ has
finitely many partial orders. So after finitely many improvements
we reach a {\it best\/} configuration, giving:

\begin{lem} Let\/ $(\si,\io,\pi)$ be an $(I,\tl)$-configuration in
an abelian or exact category. Then $(\si,\io,\pi)$ can be improved
to a best\/ $(I,\pr)$-configuration $(\si',\io',\pi')$, for some
partial order $\pr$ on $I$ dominated by~$\tl$.
\label{aa6lem1}
\end{lem}

After some preliminary results on partial orders in \S\ref{aa61},
section \ref{aa62} proves a criterion for {\it best configurations}
in terms of {\it split\/} short exact sequences.

\subsection{Partial orders $\tl,\pr$ where $\tl$ dominates $\pr$}
\label{aa61}

We study partial orders $\tl,\pr$ on $I$ where $\tl$ strictly
dominates~$\pr$.

\begin{lem} Let\/ $\tl,\pr$ be partial orders on a finite set\/
$I$, where $\tl$ strictly dominates $\pr$. Then there exist
$i,j\in I$ with\/ $i\tl j$ and\/ $i\npr j$, such that there
exists no $k\in I$ with\/ $i\ne k\ne j$ and\/ $i\tl k\tl j$.
Also $(\{j\},\{i,j\})\in\G_\sIp\sm\G_\sIt$ and\/~$(\{i,j\},
\{i\})\in\H_\sIp\sm\H_\sIt$.
\label{aa6lem2}
\end{lem}

\begin{proof} As $\tl$ strictly dominates $\pr$ there exist
$i,j\in I$ with $i\tl j$ and $i\npr j$. Suppose there exists
$k\in I$ with $i\ne k\ne j$ and $i\tl k\tl j$. Then as $i\npr j$
either (a) $i\npr k$, or (b) $k\npr j$. In case (a) we replace
$j$ by $k$, and in case (b) we replace $i$ by $k$. Then the new
$i,j$ satisfy the original conditions, but are `closer together'
than the old $i,j$. After finitely many steps we reach $i,j$
satisfying the lemma.
\end{proof}

This implies that if $\tl$ strictly dominates $\pr$ then
$\G_\sIt\subseteq\G_\sIp$ and $\H_\sIt\subseteq\H_\sIp$
in \eq{aa6eq1} are {\it strict\/} inclusions. But $\F_\sIt
\subseteq\F_\sIp$ need not be strict. For example, if
$I=\{1,2\}$ and $\tl=\le$ then $\F_\sIt=\F_\sIp$ is the set
of subsets of~$I$.

The following elementary lemma characterizes $\tl,\pr$ differing
by one step.

\begin{lem} Let\/ $(I,\tl)$ be a finite poset, and suppose
$i\ne j\in I$ with\/ $i\tl j$ but there exists no $k\in I$
with\/ $i\ne k\ne j$ and\/ $i\tl k\tl j$. Define $\pr$ on $I$
by $a\pr b$ if and only if\/ $a\tl b$ and\/ $a\ne i$, $b\ne j$.
Then $\pr$ is a partial order and\/ $\tl$ dominates $\pr$ by
one step. Conversely, if\/ $\pr$ is a partial order and $\tl$
dominates $\pr$ by one step then $\pr$ arises as above for
some unique~$i,j\in I$.
\label{aa6lem3}
\end{lem}

If $\tl$ dominates $\pr$, we can interpolate a chain of partial
orders differing by one step. The proof is elementary, using
Lemmas \ref{aa6lem2} and~\ref{aa6lem3}.

\begin{prop} Let\/ $I$ be a finite set and\/ $\pr,\tl$ partial
orders on $I$, where $\tl$ dominates $\pr$ by $s$ steps. Then
there exist partial orders $\tl=\ls_0,\ls_1,\ldots,\ls_s=\pr$
on $I$ such that\/ $\ls_{r-1}$ dominates $\ls_r$ by one step,
for~$r=1,\ldots,s$.
\label{aa6prop1}
\end{prop}

\subsection{Best $(I,\pr)$-configurations and split sequences}
\label{aa62}

We now prove a criterion for {\it best\/} $(I,\tl)$-configurations.
First we decompose certain objects $\si(J\cup K)$ as {\it direct
sums}~$\si(J)\op\si(K)$.

\begin{prop} Suppose $(\si,\io,\pi)$ is an $(I,\pr)$-configuration
in an abelian or exact category. Let\/ $J,K\in\F_\sIp$ with\/
$j\npr k$ and\/ $k\npr j$ for all\/ $j\in J$ and\/ $k\in K$. Then
$J\cup K\in\F_\sIp$ is an f-set and there is a canonical isomorphism
$\si(J)\op\si(K)\cong\si(J\cup K)$ identifying $\io_{\si(J)},
\io_{\si(K)},\pi_{\si(J)},\pi_{\si(K)}$ with\/ $\io(J,J\cup K),
\io(K,J\cup K),\pi(J\cup K,J),\pi(J\cup K,K)$ respectively. Hence
\begin{equation*}
\io(J,J\cup K)\ci\pi(J\cup K,J)+
\io(K,J\cup K)\ci\pi(J\cup K,K)=\id_{\si(J\cup K)}.
\end{equation*}
\label{aa6prop2}
\end{prop}

\begin{proof} The conditions on $J,K$ imply that
$J\cap K=\emptyset$ and $J\cup K\in\F_\sIp$ with $(J,J\cup K),
(K,J\cup K)\in\G_\sIp$ and $(J\cup K,J),(J\cup K,K)\in\H_\sIp$.
Definition \ref{aa4def}(A) applied to $(J,J\cup K)$
shows that $\pi(J\cup K,K)\ci\io(J,J\cup K)=0$, and similarly
$\pi(J\cup K,J)\ci\io(K,J\cup K)=0$. Parts (ii), (iii) and
(D) of Definition \ref{aa4def} with $J,J\cup K,J$ in place of
$J,K,L$ give $\pi(J\cup K,J)\ci\io(J,J\cup K)=\io(J,J)
\ci\pi(J,J)=\id_{\si(J)}$, and similarly $\pi(J\cup K,K)\ci
\io(K,J\cup K)=\id_{\si(K)}$. The proposition then follows
from Popescu~\cite[Cor.~2.7.4, p.~48]{Pope}.
\end{proof}

Recall that a short exact sequence $0\ra X\ra Y\ra Z\ra 0$
in $\A$ is called {\it split\/} if there is a compatible
isomorphism~$Y\cong X\op Z$.

\begin{prop} Suppose $(\si,\io,\pi)$ is an $(I,\tl)$-configuration
in an abelian or exact category, which is not best. Then there
exist\/ $i\ne j\in I$ with\/ $i\tl j$ but there exists no $k\in I$
with\/ $i\ne k\ne j$ and\/ $i\tl k\tl j$, such that the following
short exact sequence is split:
\e
\xymatrix@C=25pt{
0 \ar[r] & \si\bigl(\{i\}\bigr)
\ar[rr]^(0.45){\io(\{i\},\{i,j\})} &&
\si\bigl(\{i,j\}\bigr)
\ar[rr]^{\pi(\{i,j\},\{j\})} &&
\si\bigl(\{j\}\bigr) \ar[r] & 0.
}
\label{aa6eq2}
\e
\label{aa6prop3}
\end{prop}

\begin{proof} As $(\si,\io,\pi)$ is not best it has a strict
$(I,\pr)$-improvement $(\si',\io',\pi')$, for some $\pr$
dominated by $\tl$. Let $i,j$ be as in Lemma \ref{aa6lem2}.
Then $i\ne j$ as $i\npr j$, and there exists no $k\in I$ with
$i\ne k\ne j$ and $i\tl k\tl j$. As $i\npr j$, $j\npr i$
Proposition \ref{aa6prop2} shows that $\si'(\{i,j\})\cong
\si'(\{i\})\op\si'(\{j\})$. But $\si'(\{i\})=\si(\{i\})$,
$\si'(\{i,j\})=\si(\{i,j\})$, $\si'(\{j\})=\si(\{j\})$,
so~$\si(\{i,j\})\cong\si(\{i\})\op\si(\{j\})$.

Proposition \ref{aa6prop2} and equalities between $\io,\io'$
and $\pi,\pi'$ show that the diagram
\begin{equation*}
\xymatrix@C=20pt@R=11pt{
0 \ar[r] & \si\bigl(\{i\}\bigr)
\ar[rr]_(0.4){\io_{\si(\{i\})}}
\ar[d]_{\id_{\si(\{i\})}} &&
\si\bigl(\{i\}\bigl)\op\si\bigl(\{j\}\bigl)
\ar[rr]_(0.6){\pi_{\si(\{j\})}} \ar[d]_h
&& \si\bigl(\{j\}\bigr)\ar[r]
\ar[d]^{\id_{\si(\{j\})}} & 0 \\
0 \ar[r] & \si\bigl(\{i\}\bigr) \ar[rr]^{\io(\{i\},\{i,j\})}
&& \si\bigl(\{i,j\}\bigr) \ar[rr]^{\pi(\{i,j\},\{j\})}
&& \si\bigl(\{j\}\bigr) \ar[r] & 0
}
\end{equation*}
commutes, where $h=\io(\{i\},\{i,j\})\ci\pi_{\si(\{i\})}+\io'(
\{j\},\{i,j\})\ci\pi_{\si(\{j\})}$ is an isomorphism. Therefore
the short exact sequence \eq{aa4eq1} is split.
\end{proof}

We {\it classify improvements} for a two point indexing set
$K=\{i,j\}$. The 1-1 correspondence below is {\it not canonical},
but depends on a choice of base point; canonically, the
$(K,\ls)$-improvements form a $\Hom\bigl(\si(\{j\}),
\si(\{i\})\bigr)$-{\it torsor}.

\begin{lem} Define partial orders $\tl,\ls$ on $K=\{i,j\}$ by
$i\tl i$, $i\tl j$, $j\tl j$, $i\ls i$ and\/ $j\ls j$. Let\/
$(\si,\io,\pi)$ be a $(K,\tl)$-configuration in an abelian or
exact category. Then there exists a $(K,\ls)$-improvement\/
$(\si',\io',\pi')$ of\/ $(\si,\io,\pi)$ if and only if the
short exact sequence \eq{aa6eq2} is split, and then such\/
$(K,\ls)$-improvements $(\si',\io',\pi')$ are in $1$-$1$
correspondence with\/~$\Hom\bigl(\si(\{j\}),\si(\{i\})\bigr)$.
\label{aa6lem4}
\end{lem}

\begin{proof} If there exists a $(K,\ls)$-improvement of
$(\si,\io,\pi)$ then \eq{aa6eq2} is split by Proposition
\ref{aa6prop3}, which proves the `only if' part. For the
`if' part, suppose \eq{aa6eq2} is split. Then we can
choose morphisms $\io'(\{j\},K):\si(\{j\})\ra\si(K)$ and
$\pi'(K,\{i\}):\si(K)\ra\si(\{i\})$ with
\e
\pi'(K,\{i\})\ci\io(\{i\},K)=\id_{\si(\{i\})}
\;\>\text{and}\;\>
\pi(K,\{j\})\ci\io'(\{j\},K)=\id_{\si(\{j\})}.
\label{aa6eq3}
\e
Defining $\si'=\si$, $\io'\vert_{\G_\sKt}=\io$, $\pi'
\vert_{\H_\sKt}=\pi$ then gives a $(K,\ls)$-improvement
$(\si',\io',\pi')$ of $(\si,\io,\pi)$, proving the `if' part.

Finally, fix $\io_0'(\{j\},K),\pi_0'(K,\{i\})$ satisfying
\eq{aa6eq3}. We can easily prove that every $(K,\ls)$-improvement
$(\si',\io',\pi')$ of $(\si,\io,\pi)$ is defined uniquely by
$\si'=\si$, $\io'\vert_{\G_\sKt}=\io$, $\pi'\vert_{\H_\sKt}=\pi$ and
\begin{equation*}
\io'(\{j\},K)\!=\!\io_0'(\{j\},K)\!+\!\io(\{i\},K)\ci f,
\;\>
\pi'(K,\{i\})\!=\!\pi_0'(K,\{i\})\!-\!f\ci\pi(K,\{j\})
\end{equation*}
for some unique $f\in\Hom(\si(\{j\}),\si(\{i\}))$, and
every $f\in\Hom(\si(\{j\}),\si(\{i\}))$ gives  a
$(K,\ls)$-improvement. This establishes a 1-1 correspondence
between $(K,\ls)$-improvements $(\si',\io',\pi')$
and~$f\in\Hom\bigl(\si(\{j\}),\si(\{i\})\bigr)$.
\end{proof}

Here is the converse to Proposition~\ref{aa6prop3}.

\begin{prop} Suppose $(\si,\io,\pi)$ is an $(I,\tl)$-configuration
in an abelian or exact category. Let\/ $i\ne j\in I$ with\/ $i\tl j$
but there exists no $k\in I$ with\/ $i\ne k\ne j$ and\/ $i\tl k\tl j$,
such that\/ \eq{aa6eq2} is split. Define $\pr$ on $I$ by $a\pr b$ if\/
$a\tl b$ and\/ $a\ne i$, $b\ne j$, so that\/ $\tl$ dominates $\pr$ by
one step. Then there exists an $(I,\pr)$-improvement\/ $(\ti\si,\ti\io,
\ti\pi)$ of\/ $(\si,\io,\pi)$. Such improvements up to canonical
isomorphism are in $1$-$1$ correspondence with\/~$\Hom\bigl(\si(\{j\}),
\si(\{i\})\bigr)$.
\label{aa6prop4}
\end{prop}

\begin{proof} Set $K=\{i,j\}$, and let $(\check\si,\check\io,
\check\pi)$ be the $(K,\tl)$-subconfiguration of $(\si,\io,\pi)$.
As \eq{aa6eq2} is split, Lemma \ref{aa6lem4} shows that there
exists a $(K,\ls)$-improvement $(\si',\io',\pi')$ of
$(\check\si,\check\io,\check\pi)$. Then $(\si,\io,\pi)$
and $(\si',\io',\pi')$ satisfy the conditions of Theorem
\ref{aa5thm1} with $\phi=\id$, $I$ in place of $K$, and $K$ in
place of both $J$ and $L$. Therefore Theorem \ref{aa5thm1} gives
the $(I,\pr)$-improvement $(\ti\si,\ti\io,\ti\pi)$ that we want.

For the last part, note that every $(I,\pr)$-improvement
$(\ti\si,\ti\io,\ti\pi)$ of $(\si,\io,\pi)$ may be
constructed this way, taking $(\si',\io',\pi')$ to be the
$(K,\pr)$-subconfiguration of $(\ti\si,\ti\io,\ti\pi)$. Thus,
uniqueness up to canonical isomorphism in Theorem \ref{aa5thm1}
shows that such improvements $(\ti\si,\ti\io,\ti\pi)$ up to
canonical isomorphism are in 1-1 correspondence with
$(K,\ls)$-improvements $(\si',\io',\pi')$ of $(\check\si,\check\io,
\check\pi)$. But Lemma \ref{aa6lem4} shows that these are in 1-1
correspondence with~$\Hom\bigl(\si(\{j\}),\si(\{i\})\bigr)$.
\end{proof}

Propositions \ref{aa6prop3} and \ref{aa6prop4} imply
a {\it criterion for best configurations:}

\begin{thm} An $(I,\pr)$-configuration $(\si,\io,\pi)$ in an
abelian or exact category is best if and only if for all\/
$i\ne j\in I$ with\/ $i\pr j$ but there exists no $k\in I$
with\/ $i\ne k\ne j$ and\/ $i\pr k\pr j$, the short exact
sequence \eq{aa6eq2} is not split.
\label{aa6thm1}
\end{thm}

If this criterion holds, it also holds for any
{\it subconfiguration} of $(\si,\io,\pi)$, giving:

\begin{cor} Suppose $(\si,\io,\pi)$ is a best\/
$(I,\pr)$-configuration in an abelian or exact category.
Then all subconfigurations of\/ $(\si,\io,\pi)$ are also best.
\label{aa6cor}
\end{cor}

\section{Moduli stacks of configurations}
\label{aa7}

Let $\A$ be an abelian category. We wish to study
{\it moduli stacks of configurations} $\fM(I,\pr)_\A$,
$\fM(I,\pr,\ka)_\A$ in $\A$. To do this we shall need some
{\it extra structure} on $\A$, which is described in
Assumption \ref{aa7ass} below, and encodes information
about {\it families of objects and morphisms} in $\A$
over a {\it base scheme}~$U$.

This section will construct $\fM(I,\pr)_\A$, $\fM(I,\pr,\ka)_\A$ just
as $\K$-{\it stacks}, and some 1-{\it morphisms of\/ $\K$-stacks}
between them. But this is not enough to do algebraic geometry with.
So under some additional conditions Assumption \ref{aa8ass}, section
\ref{aa8} will prove $\fM(I,\pr)_\A$, $\fM(I,\pr,\ka)_\A$ are
{\it algebraic} $\K$-stacks, and various morphisms between them
are {\it representable} or {\it of finite type}.

\subsection{Stacks in exact categories and stacks of configurations}
\label{aa71}

Here is our first assumption, which uses ideas from
Definition~\ref{aa2def7}.

\begin{ass} Fix an algebraically closed field $\K$, and let
$\A$ be an abelian category with $\Hom(X,Y)=\Ext^0(X,Y)$
and $\Ext^1(X,Y)$ finite-dimensional $\K$-vector spaces for
all $X,Y\in\A$, and all composition maps $\Ext^i(Y,Z)\t
\Ext^j(X,Y)\ra\Ext^{i+j}(X,Z)$ bilinear for $i,j,i+j=0$ or 1.
Let $K(\A)$ be the quotient of the Grothendieck group $K_0(\A)$
by some fixed subgroup. Suppose that if $X\in\A$ with $[X]=0$
in $K(\A)$ then~$X\cong 0$.

Let $\excat$ be the 2-category whose objects are
exact categories, as in Definition \ref{aa2def3},
1-morphisms are exact functors between exact
categories, and 2-morphisms are natural transformations
between these functors. Regard $\Sch_\K$ as a 2-category
as in Definition \ref{aa2def7}, and also as a {\it site}
with the {\it \'etale topology}.

Suppose $\fF_\A:\Sch_\K\ra\excat$ is a {\it contravariant\/
$2$-functor} which is a {\it stack in exact categories} on
$\Sch_\K$ with its Grothendieck topology, that is, Definition
\ref{aa2def7}(i)--(iii) hold for $\fF_\A$, satisfying the
following conditions:
\begin{itemize}
\setlength{\itemsep}{0pt}
\setlength{\parsep}{0pt}
\item[(i)] $\fF_\A(\Spec\K)=\A$.
\item[(ii)] Let $\{f_i:U_i\ra V\}_{i\in I}$ be an open cover
of $V$ in the site $\Sch_\K$. Then a sequence $0\ra X\,\smash{
{\buildrel\phi\over\longra}\,Y\,{\buildrel\psi\over\longra}}
\,Z\ra 0$ is exact in $\fF_\A(V)$ if its images under
$\fF_\A(f_i)$ in $\fF_\A(U_i)$ are exact for all~$i\in I$.
\item[(iii)] For all $U\in\Sch_\K$ and $X\in\Obj(\fF_\A(U))$,
the map $\Hom(\Spec\K,U)\ra K(\A)$ given by $u\mapsto[\fF_\A(u)X]$
is {\it locally constant\/} in the {\it Zariski topology} on
$\Hom(\Spec\K,U)$. Here $\fF_\A(u)X\in\Obj(\fF_\A(\Spec\K))=
\Obj(\A)$ by (i), so $[\fF_\A(u)X]$ is well-defined in~$K(\A)$.
\item[(iv)] Let $X,Y\in\A$, and regard $\Hom(X,Y)$ as an
{\it affine $\K$-scheme}, with projection morphism
$\pi:\Hom(X,Y)\ra\Spec\K$. Then there should exist a
{\it tautological morphism\/} $\th_{X,Y}:\fF_\A(\pi)X\ra
\fF_\A(\pi)Y$ in $\fF_\A(\Hom(X,Y))$ such that if
$f\in\Hom(X,Y)$ and $\io_f:\Spec\K\ra\Hom(X,Y)$ is the
corresponding morphism then the following commutes in~$\A$:
\e
\begin{gathered}
\xymatrix@C=80pt@R=15pt{
\fF_\A(\io_f)\ci\fF_\A(\pi)X \ar[d]^{\fF_\A(\io_f)\th_{X,Y}}
\ar[r]_{\qquad\ep_{\pi,\io_f}(X)} & X \ar[d]_f \\
\fF_\A(\io_f)\ci\fF_\A(\pi)Y \ar[r]^{\qquad\ep_{\pi,\io_f}(Y)} & Y.
}
\end{gathered}
\label{aa7eq1}
\e
\end{itemize}
\label{aa7ass}
\end{ass}

Here is some explanation of all this.
\begin{itemize}
\setlength{\itemsep}{0pt}
\setlength{\parsep}{0pt}
\item The 2-functor $\fF_\A$ contains information
about families of objects and morphisms in $\A$. For
$U\in\Sch_\K$, objects $X$ in $\fF_\A(U)$ should
be interpreted as {\it families of objects} $X_u$ in
$\A$ parametrized by $u\in U$, which are {\it flat over}
$U$. Morphisms $\phi:X\ra Y$ in $\fF_\A(U)$ should be
interpreted as {\it families of morphisms} $\phi_u:
X_u\ra Y_u$ in $\A$ parametrized by~$u\in U$.
\item For nontrivial $U$, the condition that objects of
$\fF_\A(U)$ be {\it flat over} $U$ means that $\fF_\A(U)$ is
not an abelian category, but only an {\it exact category}, as
(co)kernels of morphisms between flat families may not be flat.
\item Families of objects and morphisms parametrized
by $U=\Spec\K$ are just objects and morphisms in $\A$,
so we take $\fF_\A(\Spec\K)=\A$ in~(i).
\item Part (ii) is necessary for the 2-functor of {\it groupoids
of exact sequences} in $\fF_\A(U)$ to be a $\K$-{\it stack}.
\item Part (iii) requires algebraic families of elements
of $\A$ to have locally constant classes in $K(\A)$. Roughly,
this means the kernel of $K_0(\A)\ra K(\A)$ includes all
continuous variations, and so cannot be `too small'.

The condition that $[X]=0$ in $K(\A)$ implies $X\cong 0$
and Assumption \ref{aa8ass} both mean the kernel of
$K_0(\A)\ra K(\A)$ cannot be `too large'.
\item Part (iv) says if we pull $X,Y\in\A$ back to
{\it constant families\/} $\pi^*(X),\pi^*(Y)$ over the base
scheme $\Hom(X,Y)$, then there is a {\it tautological morphism\/}
$\th_{X,Y}:\pi^*(X)\ra\pi^*(Y)$ taking the value $f$ over each
$f\in\Hom(X,Y)$. It will be needed in \cite[\S 6]{Joyc2} to
ensure families of configurations with constant objects $\si(J)$
but varying morphisms $\io,\pi(J,K)$ behave as expected.
\item {\it Examples} of data $\A,K(\A),\fF_\A$
satisfying Assumption \ref{aa7ass} will be given in
\S\ref{aa9} and \S\ref{aa10}. Readers may wish to refer
to them at this point.
\end{itemize}

\subsection{Moduli stacks of configurations}
\label{aa72}

We can now define {\it moduli stacks} $\fM(I,\pr)_\A$ of
$(I,\pr)$-configurations, and two other stacks $\fObj_\A,
\fExact_\A$. We will show they are $\K$-stacks in
Theorem~\ref{aa7thm1}.

\begin{dfn} We work in the situation of Assumption \ref{aa7ass}.
Define {\it contravariant\/ $2$-functors} $\fObj_\A,\fExact_\A:
\Sch_\K\ra\gro$ as follows. For $U\in\Sch_\K$, let $\fObj_\A(U)$
be the groupoid with {\it objects} $X\in\Obj(\fF_\A(U))$, and
{\it morphisms} $\phi:X\ra Y$ isomorphisms in~$\Mor(\fF_\A(U))$.

Let $\fExact_\A(U)$ be the groupoid with {\it objects}
$(X,Y,Z,\phi,\psi)$ for short exact sequences $0\ra X\,
\smash{{\buildrel\phi\over\longra}\,Y\,{\buildrel\psi\over
\longra}}\,Z\ra 0$ in $\fF_\A(U)$. Let the {\it morphisms}
in $\fExact_\A(U)$ be $(\al,\be,\ga):(X,Y,Z,\phi,\psi)\!\ra\!
(X',Y',Z',\phi',\psi')$ for $\al:X\!\ra\!X'$, $\be:Y\!\ra\!Y'$,
$\ga:Z\!\ra\!Z'$ isomorphisms in $\fF_\A(U)$ with $\phi'\ci\al\!
=\!\be\ci\phi$, $\psi'\ci\be\!=\!\ga\ci\psi$.

If $f:U\ra V$ is a 1-morphism in $\Sch_\K$ then $\fF_\A(f):
\fF_\A(V)\ra\fF_\A(U)$ induces functors $\fObj_\A(f):
\fObj_\A(V)\ra\fObj_\A(U)$ and $\fExact_\A(f):\fExact_\A(V)
\ra\fExact_\A(U)$ in the obvious way, since $\fF_\A(f)$ is
an exact functor. If $f:U\ra V$ and $g:V\ra W$ are scheme
morphisms, $\ep_{g,f}:\fF_\A(f)\ci\fF_\A(g)\ra\fF_\A(g\ci f)$
induces isomorphisms of functors $\ep_{g,f}:\fObj_\A(f)\ci
\fObj_\A(g)\ra\fObj_\A(g\ci f)$ and~$\ep_{g,f}:\fExact_\A(f)
\ci\fExact_\A(g)\ra\fExact_\A(g\ci f)$.

As in \S\ref{aa23}, this data defines the 2-functors
$\fObj_\A,\fExact_\A$. It is easy to verify $\fObj_\A,\fExact_\A$
are {\it contravariant\/ $2$-functors}. We call $\fObj_\A$ the
{\it moduli stack of objects in} $\A$, and $\fExact_\A$ the
{\it moduli stack of short exact sequences in}~$\A$.

Let $(I,\pr)$ be a finite poset. Define a {\it contravariant\/
$2$-functor} $\fM(I,\pr)_\A:\Sch_\K\ra\gro$ as follows. For
$U\in\Sch_\K$, let $\fM(I,\pr)_\A(U)$ be the groupoid with
{\it objects} $(I,\pr)$-configurations $(\si,\io,\pi)$ in
$\fF_\A(U)$, and {\it morphisms} isomorphisms of configurations
$\al:(\si,\io,\pi)\ra(\si',\io',\pi')$ in~$\fF_\A(U)$.

If $f:U\ra V$ is a scheme morphism $\fF_\A(f):\fF_\A(V)\!\ra
\!\fF_\A(U)$ is an exact functor, and so takes
$(I,\pr)$-configurations to $(I,\pr)$-configurations, and
isomorphisms of them to isomorphisms. Thus $\fF_\A(f)$
induces a functor~$\fM(I,\pr)_\A(f):\fM(I,\pr)_\A(V)\ra
\fM(I,\pr)_\A(U)$.

If $f:U\ra V$ and $g:V\ra W$ are scheme morphisms, $\ep_{g,f}:
\fF_\A(f)\ci\fF_\A(g)\ra\fF_\A(g\ci f)$ induces $\ep_{g,f}:
\fM(I,\pr)_\A(f)\ci\fM(I,\pr)_\A(g)\ra\fM(I,\pr)_\A(g\ci f)$
in the obvious way. It is easy to verify $\fM(I,\pr)_\A$ is a
{\it contravariant\/ $2$-functor}. We call $\fM(I,\pr)_\A$ the
{\it moduli stack of\/ $(I,\pr)$-configurations in}~$\A$.
\label{aa7def1}
\end{dfn}

It is usual in algebraic geometry to study moduli spaces not
of {\it all\/} coherent sheaves on a variety, but of sheaves
with a fixed Chern character or Hilbert polynomial. The analogue
for configurations $(\si,\io,\pi)$ is to fix the classes
$[\si(J)]$ in $K(\A)$ for $J\in\F_\sIp$. To do this we
introduce~$(I,\pr,\ka)$-{\it configurations}.

\begin{dfn} We work in the situation of Assumption \ref{aa7ass}.
Define
\begin{equation*}
\bar C(\A)=\bigl\{[X]\in K(\A):X\in\A\bigr\}\subset K(\A).
\end{equation*}
That is, $\bar C(\A)$ is the collection of classes in $K(\A)$ of
objects $X\in\A$. Note that $\bar C(\A)$ is {\it closed under
addition}, as $[X\op Y]=[X]+[Y]$. In \cite{Joyc2,Joyc3} we
shall make much use of $C(\A)=\bar C(\A)\sm\{0\}$. We think of
$C(\A)$ as the `positive cone' and $\bar C(\A)$ as the `closed
positive cone' in $K(\A)$, which explains the notation. For
$(I,\pr)$ a finite poset and $\ka:I\ra\bar C(\A)$, define an
$(I,\pr,\ka)$-{\it configuration} to be an $(I,\pr)$-configuration
$(\si,\io,\pi)$ with $[\si(\{i\})]=\ka(i)$ in $K(\A)$ for all~$i\in I$.

We will also use the following shorthand: we {\it extend\/ $\ka$ to
the set of subsets of\/} $I$ by defining $\ka(J)=\sum_{j\in J}\ka(j)$.
Then $\ka(J)\in\bar C(\A)$ for all $J\subseteq I$, as $\bar C(\A)$ is
closed under addition. If $(\si,\io,\pi)$ is an
$(I,\pr,\ka)$-configuration then $[\si(J)]=\ka(J)$ for all
$J\in\F_\sIp$, by Proposition~\ref{aa4prop}.
\label{aa7def2}
\end{dfn}

Here is the generalization of Definition \ref{aa7def1} to
$(I,\pr,\ka)$-configurations.

\begin{dfn} For $\al\in\bar C(\A)$, define $\fObj^\al_\A:\Sch_\K\ra
\gro$ as follows. For $U\in\Sch_\K$, let $\fObj_\A^\al(U)$ be the
full subcategory of $\fObj_\A(U)$ with objects $X\in\fObj_\A(U)$
such that $[\fObj_\A(p)X]\!=\!\al\!\in\!K(\A)$ for all
$p:\Spec\K\!\ra\!U$, so that~$\fObj_\A(p)X\!\in\!\Obj(
\fF_\A(\Spec\K))\!=\!\Obj(\A)$.

If $f:U\ra V$ and $g:V\ra W$ are scheme morphisms, $\fObj_\A(f)$
restricts to a functor $\fObj_\A^\al(f):\fObj_\A^\al(V)\ra
\fObj_\A^\al(U)$ and $\ep_{g,f}:\fObj_\A(f)\ci\fObj_\A(g)\ra
\fObj_\A(g\ci f)$ restricts to $\ep_{g,f}:\fObj_\A^\al(f)\ci
\fObj_\A^\al(g)\ra\fObj_\A^\al(g\ci f)$. We call $\fObj^\al_\A$
the {\it moduli stack of objects in $\A$ with class}~$\al$.

For $\al,\be,\ga\!\in\!\bar C(\A)$ with $\be\!=\!\al\!+\!\ga$,
define $\fExact^{\al,\be,\ga}_\A:\Sch_\K\!\ra\!\gro$ as follows.
For $U\in\Sch_\K$, let $\fExact^{\al,\be,\ga}_\A(U)$ be the full
subgroupoid of $\fExact_\A(U)$ with objects $(X,Y,Z,
\phi,\psi)$ for $X\!\in\!\Obj(\fObj^\al_\A(U))$, $Y\!\in\!\Obj(
\fObj^\be_\A(U))$ and $Z\!\in\!\Obj(\fObj^\ga_\A(U))$. Define
$\fExact^{\al,\be,\ga}_\A(f)$ and $\ep_{g,f}$ by restriction
from $\fExact_\A$. We call $\fExact^{\al,\be,\ga}_\A$ the
{\it moduli stack of short exact sequences in $\A$ with
classes}~$\al,\be,\ga$.

Now let $(I,\pr)$ be a finite poset and $\ka:I\ra\bar C(\A)$.
Define $\fM(I,\pr,\ka)_\A:\Sch_\K\ra\gro$ as follows. For
$U\in\Sch_\K$, let $\fM(I,\pr,\ka)_\A(U)$ be the full subgroupoid
of $\fM(I,\pr)_\A(U)$ with {\it objects} $(\si,\io,\pi)$ with
$\si(J)\!\in\!\Obj(\fObj_\A^{\ka(J)}(U))$ for all $J\!\in\!\F_\sIp$.
Define $\fM(I,\pr,\ka)_\A$ on morphisms and $\ep_{g,f}$ by
restricting $\fM(I,\pr)_\A$. Then $\fObj_\A^\al,
\fExact^{\al,\be,\ga}_\A,\fM(I,\pr,\ka)_\A$ are
{\it contravariant\/ $2$-functors}. We call $\fM(I,\pr,\ka)_\A$
the {\it moduli stack of\/ $(I,\pr,\ka)$-configurations in}~$\A$.
\label{aa7def3}
\end{dfn}

The basic idea here is that $\fObj^\al_\A$ contains information
on {\it families of objects} $X_u$ in $\A$ with $[X_u]=\al$ in
$K(\A)$ for $u$ in $U$, and isomorphisms between such families.
We prove the 2-functors of Definitions \ref{aa7def1} and
\ref{aa7def3} are $\K$-{\it stacks}.

\begin{thm} $\fObj_\A,\fExact_\A,\fM(I,\pr)_\A$ above are
$\K$-stacks, as in Definition \ref{aa2def7}, and\/
$\fObj^\al_\A,\fExact^{\al,\be,\ga}_\A,\fM(I,\pr,\ka)_\A$
are open and closed\/ $\K$-substacks of them, so that we
have the disjoint unions
\e
\begin{gathered}
\fObj_\A=\!\!\coprod_{\al\in\bar C(\A)}\!\!\fObj_\A^\al,\qquad
\fExact_\A=\!\!\coprod_{\substack{\al,\be,\ga\in\bar C(\A):\\
\be=\al+\ga}}\!\!\!\!\fExact_\A^{\al,\be,\ga},\\[-4pt]
\text{and}\qquad
\fM(I,\pr)_\A=\!\!\coprod_{\ka:I\ra\bar C(\A)}\!\!\fM(I,\pr,\ka)_\A.
\end{gathered}
\label{aa7eq2}
\e
\label{aa7thm1}
\end{thm}

\begin{proof} For the first part, we already know
$\fObj_\A,\fExact_\A,\fM(I,\pr)_\A$ are contravariant
2-functors $\Sch_\K\ra\gro$, and we must show they are
{\it stacks in groupoids}, that is, that Definition
\ref{aa2def7}(i)--(iii) hold. For $\fObj_\A$ this follows
immediately from $\fF_\A$ being a stack in exact categories,
as $\fObj_\A$ comes from $\fF_\A$ by omitting morphisms
which are not isomorphisms.

The proofs for $\fExact_\A,\fM(I,\pr)_\A$ are similar, and we
give only that for $\fM(I,\pr)_\A$. Let $\{f_a:U_a\ra V\}_{a\in A}$
be an open cover of $V$ in the site $\Sch_\K$. For Definition
\ref{aa2def7}(i), let $(\si,\io,\pi),(\si',\io',\pi')\in\Obj
(\fM(I,\pr)_\A(V))$ and morphisms
\begin{equation*}
\phi_a:\fM(I,\pr)_\A(f_a)(\si,\io,\pi)\ra\fM(I,\pr)_\A
(f_a)(\si',\io',\pi')\;\>\text{for $a\in A$}
\end{equation*}
satisfy \eq{aa2eq4}. Then for $J\in\F_\sIp$, applying
Definition \ref{aa2def7}(i) for $\fF_\A$ to the family of
morphisms $\phi_a(J):\fF_\A(f_a)(\si(J))\ra\fF_\A(f_a)(\si'(J))$
for $a\in A$ gives $\eta(J):\si(J)\ra\si'(J)$ in $\Mor(\fF_\A(V))$
with $\fF_\A(f_a)\eta(J)=\phi_a(J)$ for all~$a\in A$.

Moreover, Definition \ref{aa2def7}(ii) implies $\eta(J)$ is
unique. Since the $\phi_a(J)$ are isomorphisms, gluing the
$\phi_a(J)^{-1}$ in the same way yields an inverse for
$\eta(J)$, so $\eta(J)$ is an isomorphism in $\fF_\A(V)$.
By \eq{aa4eq3} and functoriality of $\fF_\A(f_a)$ we have
\begin{align*}
\fF_\A(f_a)\bigl(\eta(K)\ci\io(J,K)\bigr)&=
\phi_a(K)\!\ci\!\fF_\A(f_a)\bigl(\io(J,K)\bigr)\!= & \text{for all $a\in A$}&\\
\fF_\A(f_a)\bigl(\io'(J,K)\bigr)\ci\phi_a(J)&=
\fF_\A(f_a)\bigl(\io'(J,K)\ci\eta(J)\bigr), & \text{and $(J,K)\!\in\!\G_\sIp$.}&
\end{align*}
Therefore $\eta(K)\ci\io(J,K)=\io'(J,K)\ci\eta(J)$ for all
$(J,K)\in\G_\sIp$ by Definition \ref{aa2def7}(ii)
for $\fF_\A$. Similarly $\eta(K)\ci\pi(J,K)=\pi'
(J,K)\ci\eta(J)$ for all $(J,K)\in\H_\sIp$, so
$\eta:(\si,\io,\pi)\ra(\si',\io',\pi')$ is an
isomorphism of configurations by~\eq{aa4eq3}.

That is, $\eta\!\in\!\Mor\bigl(\fM(I,\pr)_\A(V)\bigr)$.
Clearly $\fM(I,\pr)_\A(f_a)\eta\!=\!\phi_a$ for all $a\!\in\!A$.
This proves Definition \ref{aa2def7}(i) for $\fM(I,\pr)_\A$.
For (ii), let $\eta,\ze:(\si,\io,\pi)\!\ra\!(\si',\io',\pi')$
lie in $\Mor\bigl(\fM(I,\pr)_\A(V)\bigr)$ with $\fM(I,\pr)_\A(f_a)
\eta\!=\!\fM(I,\pr)_\A(f_a)\ze$ for all $a\!\in\!A$. Then for
$J\!\in\!F_\sIp$ we have $\fF_\A(f_a)\eta(J)\!=\!\fF_\A(f_a)\ze(J)$
for $a\!\in\!A$, so Definition \ref{aa2def7}(ii) for $\fF_\A$ gives
$\eta(J)\!=\!\ze(J)$, and thus~$\eta\!=\!\ze$.

A more complicated proof using Definition \ref{aa2def7}(i)--(iii)
for $\fF_\A$ shows (iii) holds for $\fM(I,\pr)_\A$, so $\fM(I,\pr)_\A$
is a $\K$-{\it stack}. We find $\fObj^\al_\A,\fExact^{\al,\be,\ga}_\A,
\fM(I,\pr,\ka)_\A$ are $\K$-{\it substacks} of $\fObj_\A,\fExact_\A,
\fM(I,\pr)_\A$ by the same methods.

If $U$ is a connected, nonempty $\K$-scheme and $X\in\fObj_\A(U)$
then by Assumption \ref{aa7ass}(iii) the map $\Hom(\Spec\K,U)\ra
K(\A)$ given by $u\mapsto[\fF_\A(u)X]$ is locally constant on
$\Hom(\Spec\K,U)$, so $[\fF_\A(u)X]\equiv\al$ for some $\al\in
\bar C(\A)$, and therefore $X\in\fObj_\A^\al(U)$. Hence
$\fObj_\A(U)=\coprod_{\al\in\bar C(\A)}\fObj_\A^\al(U)$. It
follows easily that $\fObj_\A^\al$ is an {\it open and closed\/}
$\K$-substack of $\fObj_\A$, and the first equation of \eq{aa7eq2}
holds. The proofs for $\fExact^{\smash{\al,\be,\ga}}_\A$ and
$\fM(I,\pr,\ka)_\A$ are similar.
\end{proof}

Write $\M(I,\pr)_\A\!=\!\fM(I,\pr)_\A(\K)$, $\M(I,\pr,\ka)_\A
\!=\!\fM(I,\pr,\ka)_\A(\K)$ for the sets of geometric points
of $\fM(I,\pr)_\A$ and $\fM(I,\pr,\ka)_\A$, as in Definition
\ref{aa2def8}. We show that these are the sets of {\it
isomorphism classes of\/ $(I,\pr)$- and\/
$(I,\pr,\ka)$-configurations} in $\A$. This justifies
calling $\fM(I,\pr)_\A$, $\fM(I,\pr,\ka)_\A$ moduli
stacks of $(I,\pr)$- and $(I,\pr,\ka)$-configurations.

\begin{prop} In the situation above, $\M(I,\pr)_\A$ and\/
$\M(I,\pr,\ka)_\A$ are the sets of isomorphism classes of\/
$(I,\pr)$- and\/ $(I,\pr,\ka)$-configurations in~$\A$.
\label{aa7prop1}
\end{prop}

\begin{proof} By definition, $\M(I,\pr)_\A$ is the isomorphism
classes in the groupoid $\fM(I,\pr)_\A(\Spec\K)$. By Assumption
\ref{aa7ass}(i), objects of $\fM(I,\pr)_\A(\Spec\K)$ are
$(I,\pr)$-configurations in $\A$, and morphisms are isomorphisms
of configurations. Thus $\fM(I,\pr)_\A(\K)=\M(I,\pr)_\A$ is the
set of isomorphism classes of $(I,\pr)$-configurations. Similarly,
$\fM(I,\pr,\ka)_\A(\K)$ has objects $(I,\pr,\ka)$-configurations
and morphisms their isomorphisms, and the result follows.
\end{proof}

\subsection{Morphisms of moduli stacks}
\label{aa73}

We shall now define families of natural 1-morphisms between
the $\K$-stacks of \S\ref{aa72}. As in \S\ref{aa23}, a
1-{\it morphism of\/ $\K$-stacks} $\phi:\fF\ra\fG$ is a
{\it natural transformation} between the 2-functors $\fF,\fG$.
For each $U\in\Sch_\K$ we must provide a functor $\fF(U)\ra\fG(U)$,
satisfying some obvious natural conditions. In all our examples
these conditions hold trivially, as each $\phi(U)$ is a `forgetful
functor' omitting part of the structure, so we shall not bother to
verify them.

\begin{dfn} Define 1-morphisms $\fb,\fm,\fe:\fExact_\A\ra\fObj_\A$
as follows. For $U\in\Sch_\K$, let $\fb(U),\fm(U),\fe(U):\fExact_\A(U)
\ra\fObj_\A(U)$ act on $(X,Y,Z,\phi,\psi)$ in
$\Obj(\fExact_\A(U))$ and $(\al,\be,\ga)$ in $\Mor(\fExact_\A(U))$ by
\begin{gather*}
\fb(U):(X,Y,Z,\phi,\psi)\mapsto X,\;\>
\fb(U):(\al,\be,\ga)\mapsto\al,\;\>
\fm(U):(X,Y,Z,\phi,\psi)\mapsto Y,\\
\fm(U):(\al,\be,\ga)\mapsto\be,\quad
\fe(U):(X,Y,Z,\phi,\psi)\mapsto Z\;\>\text{and}\;\>
\fe(U):(\al,\be,\ga)\mapsto\ga.
\end{gather*}
Then $\fb,\fm,\fe$ project to the {\it beginning}, {\it
middle} and {\it end\/} objects $X,Y,Z$ respectively in
$0\ra X\ra Y\ra Z\ra 0$. For $\al,\be,\ga\in\bar C(\A)$
with $\be=\al+\ga$, these restrict to
\begin{equation*}
\fb:\fExact_\A^{\al,\be,\ga}\ra\fObj_\A^\al,\;\>
\fm:\fExact_\A^{\al,\be,\ga}\ra\fObj_\A^\be,\;\>
\fe:\fExact_\A^{\al,\be,\ga}\ra\fObj_\A^\ga.
\end{equation*}

For $(I,\pr)$ a finite poset, $\ka:I\ra\bar C(\A)$, and
$J\in\F_\sIp$, define a 1-morphism $\bs\si(J):\fM(I,\pr)_\A
\ra\fObj_\A$, where $\bs\si(J)(U):\fM(I,\pr)_\A(U)\ra\fObj_\A(U)$
acts as $\bs\si(J)(U):(\si,\io,\pi)\mapsto\si(J)$ on objects and
$\bs\si(J)(U):\al\mapsto\al(J)$ on morphisms, for $U\in\Sch_\K$.
It restricts to~$\bs\si(J):\fM(I,\pr,\ka)_\A\ra\fObj_\A^{\ka(J)}$.

As in Definition \ref{aa5def1} $\F_\sJp\!\subseteq\!\F_\sIp$,
$\G_\sJp\!\subseteq\!\G_\sIp$ and $\H_\sJp\!\subseteq\!\H_\sIp$.
Define the $(J,\pr)$-{\it subconfiguration $1$-morphism}
$S(I,\pr,J):\fM(I,\pr)_\A\ra\fM(J,\pr)_\A$ by
\begin{equation*}
S(I,\pr,J)(U):(\si,\io,\pi)\mapsto(\si\vert_{\F_\sJp},
\io\vert_{\G_\sJp},\pi\vert_{\H_\sJp}),\quad
\al\mapsto\al\vert_{\F_\sJp}
\end{equation*}
on $(\si,\io,\pi)\!\in\!\Obj(\fM(I,\pr)_\A(U))$ and $\al\!\in\!
\Mor(\fM(I,\pr)_\A(U))$, for $U\!\in\!\Sch_\K$. It restricts
to~$S(I,\pr,J):\fM(I,\pr,\ka)_\A\ra\fM(J,\pr,\ka\vert_J)_\A$.

Now let $(I,\pr)$ and $(K,\tl)$ be finite posets, and
$\phi:I\ra K$ a surjective map with $i\pr j$ implies
$\phi(i)\tl\phi(j)$ for $i,j\in I$. As in Definition
\ref{aa5def2}, pullback of subsets of $K$ to $I$ gives
injective maps $\phi^*:\F_\sKt\!\ra\!\F_\sIp$, $\phi^*:
\G_\sKt\!\ra\!\G_\sIp$ and $\phi^*:\H_\sKt\!\ra\!\H_\sIp$.
Define the {\it quotient\/ $(K,\tl)$-configuration $1$-morphism}
$Q(I,\pr,K,\tl,\phi):\fM(I,\pr)_\A\ra\fM(K,\tl)_\A$ by
\begin{equation*}
Q(I,\pr,K,\tl,\phi)(U):(\si,\io,\pi)\mapsto(\si\ci\phi^*,
\io\ci\phi^*,\pi\ci\phi^*),\quad \al\mapsto\al\ci\phi^*
\end{equation*}
on $(\si,\io,\pi)\!\in\!\Obj(\fM(I,\pr)_\A(U))$ and $\al\!\in\!
\Mor(\fM(I,\pr)_\A(U))$, for~$U\!\in\!\Sch_\K$.

In the special case when $I=K$ and $\phi:I\ra I$ is the
identity map $\id_I$, write $Q(I,\pr,\tl)=Q(I,\pr,I,\tl,\id_I)$.
Given $\ka:I\ra\bar C(\A)$, define $\mu:K\ra\bar C(\A)$ by
$\mu(k)=\ka(\phi^{-1}(k))$. Then $Q(I,\pr,K,\tl,\phi)$
restricts to $Q(I,\pr,K,\tl,\phi):\fM(I,\pr,\ka)_\A\ra
\fM(K,\tl,\mu)_\A$. When $I=K$ and $\phi=\id_I$ we have $\mu=\ka$,
so~$Q(I,\pr,\tl):\fM(I,\pr,\ka)_\A\ra\fM(I,\tl,\ka)_\A$.
\label{aa7def4}
\end{dfn}

Each of these 1-morphisms $\psi$ induces a map $\psi_*$ on
the {\it sets of geometric points} of the $\K$-stacks, as in
Definition \ref{aa2def8}. Following Proposition \ref{aa7prop1},
it is easy to show these do the obvious things.

\begin{prop} In the situation above, the induced maps act as
\begin{gather*}
\bs\si(J)_*:\M(I,\pr)_\A\!\ra\!\fObj_\A(\K),\;\>
S(I,\pr,J)_*:\M(I,\pr)_\A\ra\M(J,\pr)_\A,\\
Q(I,\pr,K,\tl,\phi)_*:\M(I,\pr)_\A\!\ra\!\M(K,\tl)_\A,\;\>
\bs\si(J)_*:[(\si,\io,\pi)]\!\mapsto\![\si(J)],\\
S(I,\pr,J)_*:[(\si,\io,\pi)]\!\mapsto\![(\si',\io',\pi')],\;\>
Q(I,\pr,K,\tl,\phi)_*:[(\si,\io,\pi)]\!\mapsto\!
[(\ti\si,\ti\io,\ti\pi)]
\end{gather*}
on geometric points, where $(\si,\io,\pi)$ is an
$(I,\pr)$-configuration in $\A$, $(\si',\io',\pi')$ its
$(J,\pr)$-subconfiguration, and\/ $(\ti\si,\ti\io,\ti\pi)$
its quotient\/ $(K,\tl)$-configuration from~$\phi$.
\label{aa7prop2}
\end{prop}

\subsection{1-isomorphisms of moduli stacks}
\label{aa74}

We conclude this section by proving that a number of
1-morphisms above are 1-{\it isomorphisms}. To show
$\phi:\fF\ra\fG$ is a 1-isomorphism, we must show that
the functor $\phi(U):\fF(U)\ra\fG(U)$ is an {\it
equivalence of categories} for each $U\in\Sch_\K$.
That is, we must prove two things:
\begin{itemize}
\setlength{\itemsep}{0pt}
\setlength{\parsep}{0pt}
\item[(a)] $\phi(U):\Obj(\fF(U))\ra\Obj(\fG(U))$ induces a
{\it bijection} between {\it isomorphism classes of objects}
in $\fF(U)$ and $\fG(U)$; and
\item[(b)] $\phi(U):\Mor(\fF(U))\ra\Mor(\fG(U))$ induces
for all $X,Y\in\Obj(\fF(U))$ a {\it bijection}~$\Hom(X,Y)
\ra\Hom\bigl(\phi(U)X,\phi(U)Y\bigr)$.
\end{itemize}

\begin{prop} {\rm(i)} Let\/ $(I,\bu)$ be a finite poset with\/
$i\bu j$ if and only if\/ $i=j$, and\/ $\ka:I\ra\bar C(\A)$. Then
the following are $1$-isomorphisms:
\begin{equation*}
\prod_{i\in I}\bs\si(\{i\}):\fM(I,\bu)_\A\ra\prod_{i\in I}\fObj_\A,\;\>
\prod_{i\in I}\bs\si(\{i\}):\fM(I,\bu,\ka)_\A\ra\prod_{i\in I}
\fObj_\A^{\ka(i)}.
\end{equation*}
{\rm(ii)} Let\/ $(\{i,j\},\pr)$ be a poset with\/ $i\pr j$ and\/
$\ka\!:\!\{i,j\}\!\ra\!\bar C(\A)$. Define $1$-morphisms
$\Pi:\fM(\{i,j\},\pr)_\A\!\ra\!\fExact_\A$ and\/ $\Pi:\fM(\{i,j\},
\pr,\ka)_\A\!\ra\!\fExact_\A^{\ka(i),\ka(\{i,j\}),\ka(j)}$~by
\begin{gather*}
\Pi(U):(\si,\io,\pi)\longmapsto\bigl(\si(\{i\}),\si(\{i,j\}),
\si(\{j\}),\io(\{i\},\{i,j\}),\pi(\{i,j\},\{j\})\bigr)\\
\text{and}\quad
\Pi(U):\al\longmapsto\bigl(\al(\{i\}),\al(\{i,j\}),\al(\{j\})\bigr)
\end{gather*}
for all\/ $(\si,\io,\pi)\in\Obj(\fM(\{i,j\},\pr)_\A(U))$,
$\al\in\Mor(\fM(\{i,j\},\pr)_\A(U))$ and\/ $U\in\Sch_\K$.
Both of these $\Pi$ are $1$-isomorphisms.
\label{aa7prop3}
\end{prop}

\begin{proof} For (i) the proof for both 1-morphisms is the
same, so we consider only the first. Let $U\in\Sch_\K$, and
$(\si,\io,\pi)\in\Obj(\fM(I,\bu)_\A(U))$. Suppose $J\subseteq
K\subseteq I$, and set $L=K\sm J$. Then $J,K,L\in\F_\sIb$,
and $(J,K),(L,K)\in\G_\sIb$. Considering the diagram
\e
\begin{gathered}
\xymatrix@R=20pt{
0 \ar@/^/[r] & \si(J) \ar@/^/[l]
\ar@/^/[rr]^{\io_{\si(J)}} \ar[d]_{\id_{\si(J)}}
&& \si(J)\op\si(L) \ar@/^/[ll]^{\pi_{\si(J)}}
\ar@/^/[rr]^{\pi_{\si(L)}} \ar@{-->}[d]
&& \si(L) \ar@/^/[ll]^{\io_{\si(L)}}
\ar[d]^{\id_{\si(L)}} \ar@/^/[r] & 0 \ar@/^/[l] \\
0 \ar@/^/[r] & \si(J) \ar@/^/[l] \ar@/^/[rr]^{\io(J,K)}
&& \si(K) \ar@/^/[ll]^{\pi(K,J)} \ar@/^/[rr]^{\pi(K,L)}
&& \si(L) \ar@/^/[ll]^{\io(L,K)} \ar@/^/[r] & 0, \ar@/^/[l]
}
\end{gathered}
\label{aa7eq5}
\e
and using Definitions \ref{aa2def1}(iii) and \ref{aa4def} gives
a canonical isomorphism $\si(J)\op\si(L)\ra\si(K)$ making
\eq{aa7eq5} commute. By induction on $\md{J}$ we construct
canonical isomorphisms $\bigop_{i\in J}\si(\{i\})\ra\si(J)$
for all $J\subseteq I$, with $\io(J,K),\pi(J,K)$ corresponding
to projections from or to subfactors in the direct sums.

Let $(X_i)_{i\in I}$ lie in $\Obj(\prod_{i\in I}\fObj_\A(U))$.
Setting $\si(J)=\bigop_{i\in J}X_i$, so that $\si(\{i\})=X_i$,
and taking $\io(J,K),\pi(J,K)$ to be the natural projections
gives $(\si,\io,\pi)$ in $\Obj(\fM(I,\bu)_\A(U))$ with
$\prod_{i\in I}\bs\si(\{i\})(U)(\si,\io,\pi)=(X_i)_{i\in I}$. Hence
\e
\ts\prod_{i\in I}\bs\si(\{i\})(U):\Obj\bigl(\fM(I,\bu)_\A(U)
\bigr)\ra\Obj\bigl(\ts\prod_{i\in I}\fObj_\A(U)\bigr)
\label{aa7eq6}
\e
is surjective. Suppose $(\si,\io,\pi),(\si',\io',\pi')$
lie in $\Obj(\fM(I,\bu)_\A(U))$ with images $(X_i)_{i\in I}$,
$(X_i')_{i\in I}$ under $\prod_{i\in I}\bs\si(\{i\})(U)$. A
morphism $(f_i)_{i\in I}:(X_i)_{i\in I}\ra(X_i')_{i\in I}$ is
a collection of isomorphisms $f_i:X_i\ra X_i'$ in~$\fF_\A(U)$.

Any such $(f_i)_{i\in I}$ extends uniquely to a morphism
$\al:(\si,\io,\pi)\ra(\si',\io',\pi')$, where $\al(\{i\})=f_i$
for $i\in I$, and more generally $\al(J)$ corresponds to
$\bigop_{i\in J}f_i:\bigop_{i\in J}X_i\ra\bigop_{i\in J}X_i'$ under
the canonical isomorphisms $\bigop_{i\in J}X_i\ra\si(J)$ and
$\bigop_{i\in J}X_i'\ra\si'(J)$. Therefore the following
map is a bijection:
\begin{equation*}
\ts\prod_{i\in I}\bs\si(\{i\})(U):
\Hom\bigl((\si,\io,\pi),(\si',\io',\pi')\bigr)\ra
\Hom\bigl((X_i)_{i\in I},(X_i')_{i\in I}\bigr).
\end{equation*}
Together with surjectivity of \eq{aa7eq6}, this shows
$\prod_{i\in I}\bs\si(\{i\})$ is a 1-isomorphism.

For (ii), both functors $\Pi(U)$ are actually
{\it isomorphisms of categories}, not just equivalences.
This is because the only data `forgotten' by $\Pi(U)$ is
$\si(\emptyset)$, $\io(\emptyset,J)$ and $\pi(J,\emptyset)$
for $J\subseteq\{i,j\}$ on objects $(\si,\io,\pi)$, and
$\al(\emptyset)$ on morphisms $\al$. But by definition
$\si(\emptyset)=0$, so $\io(\emptyset,J)=0$, $\pi(J,\emptyset)=0$
and $\al(\emptyset)=0$, and there are unique choices for the
forgotten data. Thus $\Pi$ is a 1-isomorphism.

(Note: this assumes there is a {\it prescribed zero object\/}
0 in $\fF_\A(U)$, and that $\si(\emptyset)=0$ is part of the
definition of configuration in $\fF_\A(U)$. If instead
0 in $\fF_\A(U)$ is defined {\it only up to isomorphism},
then $\Pi(U)$ `forgets' a choice of 0 in $\si(\emptyset)$,
but is still an equivalence of categories.)
\end{proof}

Our final result extends Theorem \ref{aa5thm1} to moduli
stacks of configurations.

\begin{thm} Let\/ $(J,\ls)$ and\/ $(K,\tl)$ be finite posets
and\/ $L\in\F_\sKt$, with\/ $J\cap(K\sm L)=\emptyset$.
Suppose $\psi:J\ra L$ is surjective with\/ $i\ls j$ implies
$\psi(i)\tl\psi(j)$. Set\/ $I=J\cup(K\sm L)$, and define a
partial order $\pr$ on $I$ by
\begin{equation*}
i\pr j\quad\text{for $i,j\in I$ if}\quad
\begin{cases}
i\ls j, & i,j\in J, \\
i\tl j, & i,j\in K\sm L, \\
\psi(i)\tl j, & i\in J,\quad j\in K\sm L, \\
i\tl\psi(j), & i\in K\sm L,\quad j\in J.
\end{cases}
\end{equation*}
Then $J\in\F_\sIp$ with\/ $\pr\vert_J=\ls$. Define $\phi:I\ra K$ by
$\phi(i)=\psi(i)$ if\/ $i\in J$ and\/ $\phi(i)=i$ if\/ $i\in K\sm L$.
Then $\phi$ is surjective, with\/ $i\pr j$ implies~$\phi(i)\tl\phi(j)$.

Let\/ $\ka:I\ra\bar C(\A)$, and define $\mu:K\ra K(\A)$ by
$\mu(k)=\ka(\phi^{-1}(k))$. Then the following $1$-morphism
diagrams commute and are Cartesian squares:
\e
\begin{gathered}
\xymatrix@C=30pt@R=17pt{
\fM(I,\pr)_\A \ar[r]_{S(I,\pr,J)}
\ar@<-4ex>[d]^{Q(I,\pr,K,\tl,\phi)}
& \fM(J,\ls)_\A \ar@<4ex>[d]_{Q(J,\ls,L,\tl,\psi)} \\
\fM(K,\tl)_\A \ar[r]^{S(K,\tl,L)\,\,\,} & \fM(L,\tl)_\A,
}
\quad
\xymatrix@C=30pt@R=17pt{
\fM(I,\pr,\ka)_\A \ar[r]_{S(I,\pr,J)}
\ar@<-4.5ex>[d]^{Q(I,\pr,K,\tl,\phi)}
& \fM(J,\ls,\ka)_\A \ar@<4.5ex>[d]_{Q(J,\ls,L,\tl,\psi)} \\
\fM(K,\tl,\mu)_\A \ar[r]^{S(K,\tl,L)\,\,\,} & \fM(L,\tl,\mu)_\A.
}
\end{gathered}
\label{aa7eq7}
\e
\label{aa7thm2}
\end{thm}

\begin{proof} We give the proof for the first square of
\eq{aa7eq7} only, as the second is the same. Let $U\in\Sch_\K$
and $(\si,\io,\pi)\in\Obj(\fM(I,\pr)_\A(U))$. Write
\e
\begin{gathered}
\!\!\!\!\!(\ti\si,\ti\io,\ti\pi)\!=\!Q(I,\pr,K,\tl,\phi)(U)(\si,\io,\pi),\;
(\check\si,\check\io,\check\pi)\!=\!S(K,\tl,L)(U)(\ti\si,\ti\io,\ti\pi),\\
\!\!\!\!\!(\si',\io',\pi')\!=\!S(I,\pr,J)(U)(\si,\io,\pi),\;
(\hat\si,\hat\io,\hat\pi)\!=\!Q(J,\ls,L,\tl,\psi)(U)(\si',\io',\pi').
\end{gathered}
\label{aa7eq8}
\e
Then $(\si,\io,\pi)$ is an $(I,\pr)$-configuration in $\fF_\A(U)$,
$(\ti\si,\ti\io,\ti\pi)$ its quotient $(K,\tl)$-configuration
from $\phi$, and $(\si',\io',\pi')$ its $(J,\ls)$-subconfiguration.
Lemma \ref{aa5lem1} gives $(\check\si,\check\io,\check\pi)\!=\!
(\hat\si,\hat\io,\hat\pi)$, and the analogue for morphisms also
holds. Thus
\begin{equation*}
S(K,\tl,L)(U)\ci Q(I,\pr,K,\tl,\phi)(U)=
Q(J,\ls,L,\tl,\psi)(U)\ci S(I,\pr,J)(U)
\end{equation*}
as functors $\fM(I,\pr)_\A(U)\ra\fM(L,\tl)_\A(U)$. Since this
holds for all $U\in\Sch_\K$, the first square in \eq{aa7eq7}
commutes. In fact it {\it strictly commutes}, that is, the
1-morphisms $S(K,\tl,L)\ci Q(I,\pr,K,\tl,\phi)$ and
$Q(J,\ls,L,\tl,\psi)\ci S(I,\pr,J)$ are not just
2-{\it isomorphic}, but {\it equal}.

As in Definition \ref{aa2def9} we now get a 1-morphism
\begin{equation*}
\rho:\fM(I,\pr)_\A\ra\fM(K,\tl)_\A\t_{\fM(L,\tl)_\A}\fM(J,\ls)_\A,
\end{equation*}
unique up to 2-isomorphism, making \eq{aa2eq9} commute, and
the first square of \eq{aa7eq7} is Cartesian if $\rho$ is a
1-isomorphism. From Definition \ref{aa2def9}, we deduce that
{\it objects} of $\bigl(\fM(K,\tl)_\A\t_{\fM(L,\tl)_\A}
\fM(J,\ls)_\A\bigr)(U)$ are triples
\begin{equation*}
\bigl((\ti\si,\ti\io,\ti\pi),(\si',\io',\pi'),
\al:(\check\si,\check\io,\check\pi)\,{\smash{\buildrel\cong\over
\longra}}\,(\hat\si,\hat\io,\hat\pi)\bigr),
\end{equation*}
where $(\ti\si,\ti\io,\ti\pi)\in\Obj(\fM(K,\tl)_\A(U))$,
$(\si',\io',\pi')\in\Obj(\fM(J,\ls)_\A(U))$, and
$(\check\si,\check\io,\check\pi)$, $(\hat\si,\hat\io,\hat\pi)$
are as in \eq{aa7eq8}, and~$\al\in\Mor(\fM(L,\tl)_\A(U))$.

We can define $\rho(U)$ explicitly on objects $(\si,\io,\pi)$
and morphisms $\al$ by
\begin{align*}
&\rho(U):(\si,\io,\pi)\longmapsto
\bigl((\ti\si,\ti\io,\ti\pi),(\si',\io',\pi'),
\id_{(\check\si,\check\io,\check\pi)}\bigr),\\
&\rho(U):\al\longmapsto(\al\ci\phi^*,\al\vert_{\F_\sJl}),
\end{align*}
where $(\ti\si,\ti\io,\ti\pi)$, $(\si',\io',\pi')$ and
$(\check\si,\check\io,\check\pi)$ are as in \eq{aa7eq8},
and $(\check\si,\check\io,\check\pi)=(\hat\si,\hat\io,\hat\pi)$.
Now Theorem \ref{aa5thm1} shows $\rho(U)$ is an {\it equivalence
of categories}. As this holds for all $U\in\Sch_\K$, $\rho$ is a
1-isomorphism, so the first square of \eq{aa7eq7} is Cartesian.
\end{proof}

Proposition \ref{aa7prop3} and Theorem \ref{aa7thm2} will
be used in \S\ref{aa81} to prove $\fM(I,\pr)_\A$ and
$\fM(I,\pr,\ka)_\A$ are {\it algebraic} $\K$-stacks.

\section{Algebraic $\K$-stacks of configurations}
\label{aa8}

So far we have only shown that the moduli stacks $\fM(I,\pr)_\A$,
$\fM(I,\pr,\ka)_\A$ of \S\ref{aa72} are $\K$-{\it stacks}, which
is quite a weak, categorical concept. We now impose some additional
assumptions, which will enable us to prove that $\fM(I,\pr)_\A$,
$\fM(I,\pr,\ka)_\A$ are {\it algebraic} $\K$-stacks {\it locally
of finite type}, and that various morphisms between them are {\it
representable} or {\it of finite type}.

\begin{ass} Let Assumption \ref{aa7ass} hold for $\A,K(\A),
\fF_\A$. Suppose the $\K$-stacks $\fObj^\al_\A,
\fExact^{\al,\be,\ga}_\A$ of \S\ref{aa72} are {\it algebraic}
({\it Artin}) $\K$-{\it stacks, locally of finite type}.
Suppose the following 1-morphisms of \S\ref{aa73} are of
{\it finite type}:
\e
\begin{gathered}
\fm:\fExact^{\al,\be,\ga}_\A\ra\fObj_\A^\be,\quad
\fb\t\fe:\fExact^{\al,\be,\ga}_\A\ra\fObj_\A^\al\t\fObj_\A^\ga.
\end{gathered}
\label{aa8eq1}
\e
\label{aa8ass}
\end{ass}

This list of assumptions is motivated firstly because
they hold for the examples the author is interested in,
given in \S\ref{aa9} and \S\ref{aa10}, and secondly
as the results of \S\ref{aa81}--\S\ref{aa83} that we use
them to prove will be essential for the theory of
invariants `counting' (semi)stable configurations in
$\A$ to be developed in~\cite{Joyc2,Joyc3}.

\subsection{Moduli stacks of configurations are algebraic $\K$-stacks}
\label{aa81}

The moduli stacks $\fM(I,\pr)_\A$, $\fM(I,\pr,\ka)_\A$ of
\S\ref{aa72} are {\it algebraic} $\K$-stacks.

\begin{thm} Let Assumptions \ref{aa7ass} and \ref{aa8ass} hold for
$\A,K(\A),\fF_\A$. Then the $\fM(I,\pr)_\A,\fM(I,\pr,\ka)_\A$ of\/
\S\ref{aa72} are algebraic $\K$-stacks, locally of finite type.
\label{aa8thm1}
\end{thm}

\begin{proof} The $\fM(I,\pr)_\A$ case follows from the
$\fM(I,\pr,\ka)_\A$ case by \eq{aa7eq2}, so we prove the
$\fM(I,\pr,\ka)_\A$ case. When $\md{I}=0,1$ or 2 there
are four cases:
\begin{itemize}
\setlength{\itemsep}{0pt}
\setlength{\parsep}{0pt}
\item[(a)] $I=\emptyset$. Then $\fM(I,\pr,\ka)_\A$ is
1-isomorphic to~$\Spec\K$.
\item[(b)] $I=\{i\}$. By Proposition \ref{aa7prop3}(i),
$\fM(I,\pr,\ka)_\A$ is 1-isomorphic to~$\fObj_\A^{\ka(i)}$.
\item[(c)] $I=\{i,j\}$ with $a\pr b$ if and only if $a=b$.
Then by Proposition \ref{aa7prop3}(i) $\fM(I,\pr,\ka)_\A$ is
1-isomorphic to~$\fObj_\A^{\ka(i)}\t\fObj_\A^{\ka(j)}$.
\item[(d)] $I=\{i,j\}$ with $i\pr j$. By Proposition
\ref{aa7prop3}(ii) $\fM(I,\pr,\ka)_\A$ is 1-isomorphic
to~$\fExact_\A^{\ka(i),\ka(\{i,j\}),\ka(j)}$.
\end{itemize}
In each case $\fM(I,\pr,\ka)_\A$ is an {\it algebraic $\K$-stack
locally of finite type}, by Assumption \ref{aa8ass} in (b)--(d).
So the theorem holds when~$\md{I}\le 2$.

Next we prove the case that $\pr$ is a {\it total order}, that
is, $i\pr j$ or $j\pr i$ for all $i,j\in I$. Then $(I,\pr)$ is
canonically isomorphic to $(\{1,\ldots,n\},\le)$ for $n=\md{I}$.
Suppose by induction that $\fM(I,\pr,\ka)_\A$ is an algebraic
$\K$-stack locally of finite type for all total orders $(I,\pr)$
with $\md{I}\le n$. From above this holds for $n=2$, so take
$n\ge 1$. Let $(I,\pr)$ be a total order with $\md{I}=n+1$.

Let $i$ be $\pr$-minimal in $I$, and $j$ be $\pr$-minimal in
$I\sm\{i\}$, which defines $i,j$ uniquely as $\pr$ is a total
order. Let $J=\{i,j\}$ and $L=\{l\}$ be a one point set with
$l\notin I$, set $K=\{l\}\cup I\sm\{i,j\}$, and define $\tl$
on $K$ by $a\tl b$ if either $a=l$, or $a,b\in I\sm\{i,j\}$
with $a\pr b$. Then $\tl$ is a total order on $K$, with minimal
element $l$. Define $\mu:K\ra\bar C(\A)$ by $\mu(l)=\ka(i)+\ka(j)$,
and $\mu(a)=\ka(a)$ for $a\in I\sm\{i,j\}$. Define $\phi:I\ra K$ by
$\phi(a)=l$ for $a=i,j$ and $\phi(a)=a$ otherwise. Define
$\psi:\{i,j\}\ra\{l\}$ by~$\psi(a)=l$.

Theorem \ref{aa7thm2} now applies, and shows that the diagram
\e
\begin{gathered}
\xymatrix@R=11pt{
\fM(I,\pr,\ka)_\A \ar[rr]_{S(I,\pr,\{i,j\})\qquad}
\ar[d]_{Q(I,\pr,K,\tl,\phi)}
&& \fM(\{i,j\},\pr,\ka\vert_{\{i,j\}})_\A
\ar[d]^{Q(\{i,j\},\pr,\{l\},\tl,\psi)} \\
\fM(K,\tl,\mu)_\A \ar[rr]^{S(K,\tl,\{l\})\quad} &&
\fM(\{l\},\tl,\mu\vert_{\{l\}})_\A
}
\end{gathered}
\label{aa8eq2}
\e
is commutative, and a Cartesian square. The two right hand
corners are algebraic $\K$-stacks locally of finite type
by (b), (d) above, and the bottom left hand corner is by
induction as $\tl$ is a total order and $\md{K}=n$. Hence
$\fM(I,\pr,\ka)_\A$ is an {\it algebraic $\K$-stack
locally of finite type}, by properties of Cartesian squares
in \S\ref{aa23}. By induction, this holds whenever $\pr$ is
a {\it total order}.

Now let $(I,\pr)$ be a finite poset. Define $S_\pr=\bigl\{(i,j)
\in I\t I:i\npr j$ and $j\npr i\bigr\}$, and let $n_\pr=\md{S_\pr}$.
Let $\ls$ be a total order on $I$ which dominates $\pr$. Then
$\ls$ dominates $\pr$ by $n_\pr$ steps. If $n_\pr>0$, by
Proposition \ref{aa6prop1} there exists $\tl$ on $I$ such that
$\ls$ dominates $\tl$ by $n_\pr-1$ steps, so that $n_\tl=n_\pr-1$,
and $\tl$ dominates $\pr$ by one step. By Lemma \ref{aa6lem2},
there exist unique $i,j\in I$ with $i\tl j$ but~$i\npr j$.

Suppose by induction that $\fM(I,\pr,\ka)_\A$ is an algebraic
$\K$-stack locally of finite type whenever $n_\pr\le n$, for
$n\ge 0$. When $n=0$ this implies $\pr$ is a {\it total order},
so the first step $n=0$ holds from above. Let $(I,\pr),\ka$
have $n_\pr=n+1$. Then from above there is a partial order
$\tl$ on $I$ dominating $\pr$ by one step, so that $n_\tl=n$,
and unique $i,j\in I$ with $i\tl j$ but $i\npr j$. Define
$K=I$, $J=L=\{i,j\}$, $\phi:I\ra K$ and $\psi:J\ra L$ to
be the identity maps, and~$\mu=\ka$.

Theorem \ref{aa7thm2} now applies, and shows that the diagram
\e
\begin{gathered}
\xymatrix@R=11pt{
\fM(I,\pr,\ka)_\A \ar[rr]_{S(I,\pr,\{i,j\})\qquad}
\ar[d]_{Q(I,\pr,\tl)}
&& \fM(\{i,j\},\pr,\ka\vert_{\{i,j\}})_\A
\ar[d]^{Q(\{i,j\},\pr,\tl)} \\
\fM(I,\tl,\ka)_\A \ar[rr]^{S(I,\tl,\{i,j\})\qquad}
&& \fM(\{i,j\},\tl,\ka\vert_{\{i,j\}})_\A
}
\end{gathered}
\label{aa8eq3}
\e
is commutative, and a {\it Cartesian square}. The two right
hand corners are algebraic $\K$-stacks locally of finite type
by (c), (d) above, and the bottom left hand corner is by
induction, as $n_\tl=n$. Hence $\fM(I,\pr,\ka)_\A$ is an
{\it algebraic $\K$-stack locally of finite type}. By
induction, this completes the proof.
\end{proof}

The underlying idea in this proof is that
$\fM(I,\pr,\ka)_\A$ is 1-isomorphic to a complicated
{\it multiple fibre product}, constructed from many
copies of the $\fObj_\A^\al$ and $\fExact_\A^{\al,\be,\ga}$.
As the class of algebraic $\K$-stacks locally of finite
type is {\it closed under $1$-isomorphisms and fibre
products}, and $\fObj_\A^\al$ and $\fExact_\A^{\al,\be,\ga}$
lie in this class by Assumption \ref{aa8ass}, we see that
$\fM(I,\pr,\ka)_\A$ also lies in this class.

\subsection{Representable and finite type morphisms}
\label{aa82}

Next we show some 1-morphisms from \S\ref{aa73} are {\it
representable}, or of {\it finite type}. We begin with
1-morphisms involving {\it two point posets}.

\begin{prop} In the situation above, let\/ $\pr,\tl$ be
partial orders on $\{i,j\}$ with\/ $a\pr b$ only if\/ $a=b$
and\/ $i\tl j$, and let\/ $\ka:\{i,j\}\ra\bar C(\A)$. Then
\begin{itemize}
\setlength{\itemsep}{0pt}
\setlength{\parsep}{0pt}
\item[{\rm(a)}] $Q(\{i,j\},\pr,\tl):\fM(\{i,j\},\pr,\ka)_\A
\ra\fM(\{i,j\},\tl,\ka)_\A$ is representable and of finite type.
\item[{\rm(b)}] $\bs\si(\{i,j\})\!:\!\fM(\{i,j\},\tl,\ka)_\A\!\ra\!
\fObj_\A^{\ka(\{i,j\})}$ is representable and finite type.
\item[{\rm(c)}] $\bs\si(\{i\})\!\t\!\bs\si(\{j\}):\fM(\{i,j\},\tl,
\ka)_\A\!\ra\!\fObj_\A^{\ka(i)}\!\t\!\fObj_\A^{\ka(j)}$ is of
finite type.
\end{itemize}
\label{aa8prop1}
\end{prop}

\begin{proof} For (a), consider the equality of 1-morphisms
\e
\bs\si(\{i\})\t\bs\si(\{j\})=\bigl(\bs\si(\{i\})\t\bs\si(\{j\})\bigr)
\ci Q(\{i,j\},\pr,\tl)
\label{aa8eq4}
\e
acting $\fM(\{i,j\},\pr,\ka)_\A\!\!\ra\!\fObj_\A^{\ka(i)}\!\t\!
\fObj_\A^{\ka(j)}$ or $\fM(\{i,j\},\pr)_\A\!\!\ra\!\fObj_\A\!\t\!
\fObj_\A$. By Proposition \ref{aa7prop3}(i), the l.h.s.\ of
\eq{aa8eq4} is a 1-isomorphism, and so is representable and
finite type. But if $\fF\,\smash{{\buildrel\phi\over\longra}\,
\fG\,{\buildrel\psi\over\longra}}\,\fH$ are 1-morphisms of
algebraic $\K$-stacks and $\psi\ci\phi$ is representable and
finite type, then $\phi$ is too by \cite[Lem.~3.12(c)(ii) \&
Rem.~4.17(1)]{LaMo}. Hence $Q(\{i,j\},\pr,\tl)$ is representable
and of finite type.

In (b) it is easy to see that $\bs\si(\{i,j\})=\fm\ci\Pi$, where
$\Pi$ is the 1-isomorphism of Proposition \ref{aa7prop3}(ii) with
$\tl$ in place of $\pr$, and
\e
\fm:\fExact_\A^{\ka(i),\ka(i)+\ka(j),\ka(j)}\longra
\fObj_\A^{\ka(i)+\ka(j)}
\label{aa8eq5}
\e
is as in Definition \ref{aa7def4}. Since $\Pi$ is a 1-isomorphism,
as \eq{aa8eq5} is of finite type by Assumption \ref{aa8ass},
$\bs\si(\{i,j\})$ in (b) is of finite type. To show it is
representable, we must show \eq{aa8eq5} is representable.

From \cite[Cor.~8.1.1]{LaMo} we deduce the following necessary
and sufficient condition for a 1-morphism $\phi:\fF\ra\fG$
of {\it algebraic} $\K$-stacks to be {\it representable}:
for all $U\in\Sch_\K$ and all $X\in\Obj(\fF(U))$, the map
\e
\phi(U):\Hom(X,X)\ra\Hom(\phi(U)X,\phi(U)X)
\label{aa8eq6}
\e
induced by the functor $\phi(U):\fF(U)\ra\fG(U)$ should be
injective. Since \eq{aa8eq6} is a group homomorphism, it is
enough that $\phi(U)\al=\id_{\phi(U)X}$ implies~$\al=\id_X$.

Let $(\al,\be,\ga)\!:\!(X,Y,Z,\phi,\psi)\!\ra\!(X,Y,Z,\phi,\psi)$
in $\Mor(\fExact_\A^{\ka(i),\ka(\{i,j\}),\ka(j)}(U))$ or
$\Mor(\fExact_\A(U))$. Then $\fm(U)(X,Y,Z,\phi,\psi)=Y$ and
$\fm(U)(\al,\be,\ga)=\be$, so to show $\fm$ is representable
we must prove that $\be=\id_Y$ implies $\al=\id_X$ and $\ga=\id_Z$.
By definition of $\Mor(\fExact_\A(U))$ we have
\begin{equation*}
\phi\ci\al=\be\ci\phi=\id_Y\ci\phi=\phi\ci\id_X
\quad\text{and}\quad
\id_Z\ci\psi=\psi\ci\id_Y=\psi\ci\be=\ga\ci\psi.
\end{equation*}
As $\phi$ is injective and $\psi$ surjective these imply $\al=\id_X$
and $\ga=\id_Z$. Thus \eq{aa8eq5} is representable, proving (b).
Finally, in (c) we have $\bs\si(\{i\})\t\bs\si(\{j\})=(\fb\t\fe)\ci
\Pi$. But $\fb\t\fe$ is of finite type by Assumption \ref{aa8ass},
and $\Pi$ is a 1-isomorphism. So $\bs\si(\{i\})\t\bs\si(\{j\})$
is of finite type.
\end{proof}

Using this and inductive methods as in Theorem \ref{aa8thm1},
we show:

\begin{thm} If Assumptions \ref{aa7ass} and \ref{aa8ass} hold then
\begin{itemize}
\setlength{\itemsep}{0pt}
\setlength{\parsep}{0pt}
\item[{\rm(a)}] In Definition \ref{aa7def4}, the following are
representable and of finite type:
\begin{align}
&Q(I,\pr,\tl):\fM(I,\pr,\ka)_\A\ra\fM(I,\tl,\ka)_\A,
\label{aa8eq7}\\
&Q(I,\pr,\tl):\fM(I,\pr)_\A\ra\fM(I,\tl)_\A,
\label{aa8eq8}\\
&Q(I,\pr,K,\tl,\phi):\fM(I,\pr,\ka)_\A\ra\fM(K,\tl,\mu)_\A,
\label{aa8eq9}
\end{align}
and also $Q(I,\pr,K,\tl,\phi):\fM(I,\pr)_\A\ra\fM(K,\tl)_\A$
is representable.
\item[{\rm(b)}] $\bs\si(I):\fM(I,\pr,\ka)_\A\ra\fObj^{\ka(I)}_\A$
is representable and of finite type.
\item[{\rm(c)}] $\bs\si(I):\fM(I,\pr)_\A\ra\fObj_\A$ is representable.
\item[{\rm(d)}] $\prod_{i\in I}\bs\si(\{i\}):\fM(I,\pr,\ka)_\A\ra
\prod_{i\in I}\fObj_\A^{\ka(i)}$ and
\hfil\break
$\prod_{i\in I}\bs\si(\{i\}):\fM(I,\pr)_\A\ra\prod_{i\in I}\fObj_\A$
are of finite type.
\end{itemize}
\label{aa8thm2}
\end{thm}

\begin{proof} For \eq{aa8eq7}, first suppose $\tl$ dominates
$\pr$ by one step. Then as in the proof of Theorem \ref{aa8thm1},
$Q(I,\pr,\tl)$ fits into a Cartesian square \eq{aa8eq3}.
The right hand morphism $Q(\{i,j\},\pr,\tl)$ in \eq{aa8eq3} is
representable and finite type by Proposition \ref{aa8prop1}(a).
Hence the left hand morphism \eq{aa8eq7} is representable and
finite type.

When $\tl$ dominates $\pr$ by $s$ steps, by Proposition
\ref{aa6prop1} we may write $Q(I,\pr,\tl)$ as the composition
of $s$ 1-morphisms $Q(I,\ls_r,\ls_{r-1})$ with $\ls_{r-1}$
dominating $\ls_r$ by one step. By \cite[Lem.~3.12(b)]{LaMo}
compositions of representable or finite type 1-morphisms are
too, so \eq{aa8eq7} is representable and finite type by the
first part.

From \eq{aa7eq2} we have $\fM(I,\tl)_\A=\coprod_\ka
\fM(I,\tl,\ka)_\A$, and \eq{aa8eq8} coincides with
\eq{aa8eq7} over the open substack $\fM(I,\tl,\ka)_\A$ of
$\fM(I,\tl)_\A$. For a 1-morphism $\phi:\fF\ra\fG$ to be
representable or finite type is a local condition on $\fG$.
Thus, \eq{aa8eq8} is representable and finite type as
\eq{aa8eq7} is.

Next we prove (b). Suppose by induction that (b) holds
whenever $(I,\pr)$ is a {\it total order} with $\md{I}\le n$.
Follow the middle part of the proof of Theorem \ref{aa8thm1}.
Now $Q(\{i,j\},\pr,\{l\},\tl,\psi)$ in \eq{aa8eq2} is identified
with $\bs\si(\{i,j\})$ in Proposition \ref{aa8prop1}(b) by the
1-isomorphism $\fM(\{l\},\tl,\mu\vert_{\{l\}})_\A\ra
\fObj_\A^{\ka(\{i,j\})}$ of Proposition \ref{aa7prop3}(i). Thus
$Q(\{i,j\},\pr,\{l\},\tl,\psi)$ in \eq{aa8eq2} is representable
and finite type as $\bs\si(\{i,j\})$ is, so $Q(I,\pr,K,\tl,\phi)$
in \eq{aa8eq2} is representable and finite type.

But $\bs\si(I)=\bs\si(K)\ci Q(I,\pr,K,\tl,\phi):
\fM(I,\pr,\ka)_\A\ra\fObj^{\ka(I)}_\A$, where $\bs\si(K)$
is representable and finite type by induction, and
$Q(I,\pr,K,\tl,\phi)$ is representable and finite
type from above. Thus $\bs\si(I)$ is representable and
finite type by \cite[Lem.~3.12(b)]{LaMo}. By induction,
this proves (b) whenever $\pr$ is a total order.

For the general case, let $(I,\pr)$ be a finite poset, $\tl$
a total order on $I$ dominating $\pr$, and $\ka:I\ra\bar C(\A)$.
Then $\bs\si(I)=\bs\si(I)\ci Q(I,\pr,\tl):\fM(I,\pr,\ka)_\A
\ra\fObj^{\ka(I)}_\A$, where the second $\bs\si(I)$
acts on $\fM(I,\tl,\ka)_\A$, and is representable and
finite type as $\tl$ is a total order. But $Q(I,\pr,\tl)$
is representable and finite type by \eq{aa8eq7} in (a).
Therefore the composition $\bs\si(I)$ is, giving~(b).

For $\phi:\fF\ra\fG$ to be representable is a local
condition on both $\fF$ and $\fG$. By \eq{aa7eq2} we see
that $\bs\si(I):\fM(I,\pr)_\A\ra\fObj_\A$ locally coincides
with $\bs\si(I):\fM(I,\pr,\ka)_\A\ra\fObj^{\ka(I)}_\A$,
which is representable by (b). This proves~(c).

We can now complete (a). If $\fF\,\smash{{\buildrel\phi\over
\longra}\,\fG\,{\buildrel\psi\over\longra}}\,\fH$ are 1-morphisms
of algebraic $\K$-stacks and $\psi\ci\phi$ is representable
and finite type, then $\phi$ is too by \cite[Lem.~3.12(c)(ii)
\& Rem.~4.17(1)]{LaMo}. For \eq{aa8eq9} we have
\begin{equation*}
\bs\si(I)=\bs\si(K)\ci Q(I,\pr,K,\tl,\phi):\fM(I,\pr,\ka)_\A
\ra\fObj_\A^{\ka(I)}.
\end{equation*}
As $\bs\si(I)$ is representable and finite type by (b), we
see \eq{aa8eq9} is. Similarly, $Q(I,\pr,K,\tl,\phi):\fM(I,\pr)_\A
\!\ra\!\fM(K,\tl)_\A$ is representable using~(c).

Finally we prove (d). When $\md{I}=0,1$ the $\prod_{i\in I}\bs\si(\{i\})$
are 1-isomorphisms, so of finite type. Suppose by induction that the
first line of (d) is finite type whenever $(I,\pr)$ is a {\it total
order} with $\md{I}\le n$, for $n\ge 1$. Let $(I,\pr)$ be a total
order with $\md{I}=n+1$, and define $J,K,L,\tl,\mu$ as in the proof
of Theorem \ref{aa8thm1}. Since \eq{aa8eq2} is Cartesian and
$\bs\si(\{a\})=\bs\si(\{a\})\ci Q(I,\pr,K,\tl,\phi):
\fM(I,\pr,\ka)_\A\ra\fObj_\A^{\ka(a)}$ for all $a\in I\sm\{i,j\}$,
we have a Cartesian square
\e
\begin{gathered}
\xymatrix@R=20pt{
\fM(I,\pr,\ka)_\A \ar[rr]
\ar@{}@<-.8ex>[rr]_{S(I,\pr,\{i,j\})\t
\prod_{a\in I\sm\{i,j\}}\bs\si(\{a\})}
\ar@<-4ex>[d]^(0.47){Q(I,\pr,K,\tl,\phi)}
&& \fM(\{i,j\},\pr,\ka\vert_{\{i,j\}})_\A\t
\prod_{a\in I\sm\{i,j\}}\fObj_\A^{\ka(a)}
\ar@<14ex>[d]_{\substack{Q(\{i,j\},\pr,\{l\},\tl,\psi)\t\\
\prod_{a\in I\sm\{i,j\}}\id_{\fObj_\A^{\ka(a)}}}}
\\
\fM(K,\tl,\mu)_\A \ar[rr]
\ar@{}@<.8ex>[rr]^{S(K,\tl,\{l\})\t
\prod_{a\in I\sm\{i,j\}}\bs\si(\{a\})} &&
\fM(\{l\},\tl,\mu\vert_{\{l\}})_\A\t
\prod_{a\in I\sm\{i,j\}}\fObj_\A^{\ka(a)}.
}
\end{gathered}
\label{aa8eq10}
\e

Now $\bs\si(\{l\}):\fM(\{l\},\tl,\mu\vert_{\{l\}})_\A\ra
\fObj_\A^{\mu(l)}$ is a 1-isomorphism by Proposition
\ref{aa7prop3}(i). Composing with this
identifies the bottom morphism of \eq{aa8eq10} with
$\prod_{a\in K}\bs\si(\{a\})$, which is finite type by
induction. Hence the bottom morphism in \eq{aa8eq10} is
finite type, so the top morphism in \eq{aa8eq10} is too.
But $\bs\si(\{i\})\!\t\!\bs\si(\{j\})$ is finite type by
Proposition \ref{aa8prop1}(c). Composing with the top
morphism of \eq{aa8eq10} gives $\prod_{a\in I}\bs\si(\{a\})$,
which is therefore finite type. So by induction, the first
line of (d) is finite type for $(I,\pr)$ a total order.

Suppose $(I,\pr)$ is a finite poset, $\tl$ a total order on $I$
dominating $\pr$, and $\ka:I\ra\bar C(\A)$. Then $\prod_{i\in I}
\bs\si(\{i\})\!=\!\prod_{i\in I}\bs\si(\{i\})\ci Q(I,\pr,\tl):
\fM(I,\pr,\ka)_\A\!\ra\!\prod_{i\in I}\fObj^{\ka(i)}_\A$,
where the second $\prod_{i\in I}\bs\si(\{i\})$ acts on
$\fM(I,\tl,\ka)_\A$, and is finite type as $\tl$ is
a total order. As $Q(I,\pr,\tl)$ is finite type by (a),
composition gives the first line of (d), and the second
line follows as for~\eq{aa8eq8}.
\end{proof}

\subsection{The moduli spaces $\fM(X,I,\pr)_\A$ and
$\fM(X,I,\pr,\ka)_\A$}
\label{aa83}

We can now form two further classes of moduli spaces of
$(I,\pr)$-configurations $(\si,\io,\pi)$ with $\si(I)=X$,
for $X$ a fixed object in~$\A$.

\begin{dfn} In the situation above, let $X\in\Obj\A$. Assumption
\ref{aa7ass}(i) identifies $X$ with a 1-morphism $\psi:\Spec\K\ra
\fObj_\A^{[X]}$ or $\fObj_\A$. For $(I,\pr)$ a finite poset and
$\ka:I\ra\bar C(\A)$ with $\ka(I)=[X]$ in $\bar C(\A)$, define
\e
\begin{split}
\fM(X,I,\pr)_\A&=\fM(I,\pr)_\A\t_{\bs\si(I),\fObj_\A,\psi}\Spec\K\\
\text{and}\quad
\fM(X,I,\pr,\ka)_\A&=\fM(I,\pr,\ka)_\A\t_{\bs\si(I),\fObj^{\ka(I)}_\A,
\psi}\Spec\K.\\
\end{split}
\label{aa8eq11}
\e
As $\fM(I,\pr)_\A$, $\fM(I,\pr,\ka)_\A$ are algebraic $\K$-stacks
locally of finite type by Theorem \ref{aa8thm1}, Theorem
\ref{aa8thm2}(b) implies that $\fM(X,I,\pr,\ka)_\A$ is
{\it represented by an algebraic $\K$-space of finite type},
and Theorem \ref{aa8thm2}(c)  that $\fM(X,I,\pr)_\A$ is {\it
represented by an algebraic $\K$-space locally of finite type}.
Write $\Pi_X\!:\!\fM(X,I,\pr)_\A\!\ra\!\fM(I,\pr)_\A$ and
$\Pi_X\!:\!\fM(X,I,\pr,\ka)_\A\!\ra\!\fM(I,\pr,\ka)_\A$ for the
1-morphisms of stacks from the fibre products. Write $\M(X,I,\pr)_\A
\!=\!\fM(X,I,\pr)_\A(\K)$ and $\M(X,I,\pr,\ka)_\A\!=\!\fM(X,I,
\pr,\ka)_\A(\K)$ for their sets of geometric points.
\label{aa8def}
\end{dfn}

For the examples of \S\ref{aa9} and \S\ref{aa10} the
$\fM(X,I,\pr,\ka)_\A$ are actually {\it represented by
quasiprojective $\K$-schemes}. The reason for this is
that \eq{aa8eq1} are {\it quasiprojective} 1-morphisms.
Replacing finite type with quasiprojective 1-morphisms
in the proofs of \S\ref{aa82} shows that $\bs\si(I):
\fM(I,\pr,\ka)_\A\ra\fObj^{\ka(I)}_\A$ is {\it
representable and quasiprojective}, implying that
$\fM(X,I,\pr,\ka)_\A$ is represented by a
quasiprojective $\K$-scheme in Definition \ref{aa8def}.
Here is the analogue of Proposition~\ref{aa7prop1}.

\begin{prop} In Definition \ref{aa8def}, $\M(X,I,\pr)_\A$
and\/ $\M(X,I,\pr,\ka)_\A$ are naturally identified with
the sets of isomorphism classes of\/ $(I,\pr)$- and\/
$(I,\pr,\ka)$-configurations\/ $(\si,\io,\pi)$ in $\A$
with $\si(I)=X$, modulo isomorphisms\/ $\al:(\si,\io,\pi)
\ra(\si',\io',\pi')$ of\/ $(I,\pr)$-configurations
with\/~$\al(I)=\id_X$.
\label{aa8prop2}
\end{prop}

\begin{proof} By Definition \ref{aa2def9}, we find the
groupoid $\fM(X,I,\pr)_\A(\Spec\K)$ has objects
$\bigl((\ti\si,\ti\io,\ti\pi),\id_{\Spec\K},\phi\bigr)$,
for $(\ti\si,\ti\io,\ti\pi)$ an $(I,\pr)$-configuration in
$\A$, and $\phi:\ti\si(I)\ra X$ an isomorphism in $\A$.
Given such a $\bigl((\ti\si,\ti\io,\ti\pi),\id_{\Spec\K},
\phi\bigr)$, define an $(I,\pr)$-configuration $(\si,\io,\pi)$
by $\si(I)=X$, $\si(J)=\ti\si(J)$ for $J\ne I$, and
\begin{equation*}
\io(J,K)\!=\!\begin{cases} \ti\io(J,K), & J\!\ne\! I\!\ne\! K,\\
\phi\!\ci\!\ti\io(J,K),\!\!\! & J\!\ne\! I\!=\!K,\\
\id_X, & J\!=\!I\!=\!K,\end{cases}\;\>
\pi(J,K)\!=\!\begin{cases} \ti\pi(J,K), & J\!\ne\! I\!\ne\! K,\\
\ti\pi(J,K)\!\ci\!\phi^{-1},\!\!\! & J\!=\!I\!\ne\! K,\\
\id_X, & J\!=\!I\!=\!K.\end{cases}
\end{equation*}

Define $\be:(\ti\si,\ti\io,\ti\pi)\ra(\si,\io,\pi)$ by
$\be(J)\!=\!\id_{\ti\si(J)}$ if $I\ne J\in\F_\sIp$, and
$\be(I)\!=\!\phi$. Then $\be$ is an isomorphism of
$(I,\pr)$-configurations. Moreover,
\begin{equation*}
\bigl(\be,\id_{\id_{\Spec\K}}\bigr):
\bigl((\ti\si,\ti\io,\ti\pi),\id_{\Spec\K},\phi\bigr)\longra
\bigl((\si,\io,\pi),\id_{\Spec\K},\id_X\bigr)
\end{equation*}
is an isomorphism in $\fM(X,I,\pr)_\A(\Spec\K)$.
So each object of $\fM(X,I,\pr)_\A\ab
(\Spec\K)$ is isomorphic to some $((\si,\io,\pi),
\id_{\Spec\K},\id_X)$ for an $(I,\pr)$-configuration
$(\si,\io,\pi)$ with $\si(I)\!=\!X$. Isomorphisms
\begin{equation*}
\bigl(\al,\id_{\id_{\Spec\K}}\bigr):
\bigl((\si,\io,\pi),\id_{\Spec\K},\id_X\bigr)\longra
\bigl((\si',\io',\pi'),\id_{\Spec\K},\id_X\bigr)
\end{equation*}
between two elements of this form come from
isomorphisms $\al:(\si,\io,\pi)\ra(\si',\io',\pi')$ of
$(I,\pr)$-configurations with $\al(I)=\id_X$. Hence,
$\M(X,I,\pr)_\A$ is naturally identified with the set
of isomorphism classes of $(I,\pr)$-configurations
$(\si,\io,\pi)$ in $\A$ with $\si(I)=X$, modulo
isomorphisms $\al:(\si,\io,\pi)\ra(\si',\io',\pi')$
with $\al(I)=\id_X$. The proof for $\M(X,I,\pr,\ka)_\A$
is the same.
\end{proof}

\section{Coherent sheaves on a projective scheme}
\label{aa9}

Let $\K$ be an algebraically closed field, $P$ a projective
$\K$-scheme, and $\coh(P)$ the {\it abelian category of
coherent sheaves} on $P$. We shall apply the machinery of
\S\ref{aa7}--\S\ref{aa8} to $\A=\coh(P)$. Section \ref{aa91}
defines the data $\A,K(\A),\fF_\A$ required by Assumption
\ref{aa7ass}. Then \S\ref{aa92} proves that Assumption
\ref{aa7ass} holds, and \S\ref{aa93}--\S\ref{aa94} that
Assumption \ref{aa8ass} holds, for this data.

Thus by \S\ref{aa7}--\S\ref{aa8}, we have well-defined
{\it moduli stacks of configurations of coherent sheaves}
$\fM(I,\pr)_{\coh(P)}$, $\fM(I,\pr,\ka)_{\coh(P)}$, which
are {\it algebraic $\K$-stacks, locally of finite type},
and many 1-morphisms between them, some of which are
{\it representable} or {\it of finite type}. For background
on {\it coherent\/} and {\it quasicoherent sheaves} see
Hartshorne \cite[\S II.5]{Hart} or
Grothendieck~\cite[\S I.0.5]{Grot2}.

\subsection{Definition of the data $\A,K(\A),\fF_\A$}
\label{aa91}

Our first two examples define the data $\A,K(\A),\fF_\A$
of Assumption \ref{aa7ass} for the {\it abelian category of
coherent sheaves} on a smooth projective $\K$-scheme $P$.
The assumption that $P$ is {\it smooth\/} will be relaxed
in Example~\ref{aa9ex2}.

\begin{ex} Let $\K$ be an algebraically closed field and
$P$ a smooth projective $\K$-scheme, and take $\A$ to be
the {\it abelian category $\coh(P)$ of coherent sheaves}
on $P$. Then one may define the {\it Chern character}
${\mathop{\rm ch}}:K_0(\A)\ra H^{\rm even}(P,\Z)$, a
homomorphism of abelian groups. Let $K(\A)$ be the
quotient of $K_0(\A)$ by $\Ker({\mathop{\rm ch}})$. Then
$\mathop{\rm ch}$ identifies $K(\A)$ with a subgroup
of~$H^{\rm even}(P,\Z)$.

Motivated by \cite[\S 2.4.4]{LaMo}, for $U\in\Sch_\K$ define
$\fF_\A(U)=\fF_{\coh(P)}(U)$ to be the  category of {\it
finitely presentable quasicoherent sheaves\/} on $P\t U$,
as in \cite[\S I.0.5]{Grot2}, which are {\it flat over\/} $U$,
as in \cite[I.0.6.7]{Grot2}. This is a full additive subcategory
of the abelian category $\qcoh(P\t U)$ of quasicoherent sheaves on
$P\t U$, closed under extensions by \cite[Prop.~IV.2.1.8]{Grot2},
so it is an exact category.

If $f:U\!\ra\!V$ is a morphism in $\Sch_\K$ then so is
$\id_P\!\t\!f:P\!\t\!U\!\ra\!P\!\t\!V$. Define a
functor $\fF_{\coh(P)}(f):\fF_{\coh(P)}(V)\!\ra\!
\fF_{\coh(P)}(U)$ by pullback $(\id_P\t f)^*$ of
sheaves and their morphisms along $\id_P\t f$. That is,
if $X\in\fF_{\coh(P)}(V)$ then $\fF_{\coh(P)}(f)X$ is
the {\it inverse image sheaf\/} $(\id_P\t f)^*(X)$ on
$P\t U$, as in~\cite[p.~110]{Hart}.

Then $\fF_{\coh(P)}(f)X=(\id_P\t f)^*(X)$ is quasicoherent
by \cite[I.0.5.1.4]{Grot2} or \cite[Prop.~II.5.8(a)]{Hart},
finitely presentable by \cite[I.0.5.2.5]{Grot2}, and flat
over $U$ by \cite[Prop.~IV.2.1.4]{Grot2}. Thus
$\fF_{\coh(P)}(f)X\in\fF_{\coh(P)}(U)$, as we need. Also,
$(\id_P\t f)^*$ takes exact sequences of quasicoherent
sheaves on $P\t V$ flat over $V$ to exact sequences of
quasicoherent sheaves on $P\t U$ flat over $U$ by
\cite[Prop.~IV.2.1.8(i)]{Grot2}. Hence $\fF_{\coh(P)}(f)$ is
an exact functor.

Here is a slightly subtle point. Inverse images come from
a universal construction, and so are given not uniquely,
but only {\it up to canonical isomorphism}. So there could
be {\it many possibilities} for $(\id_P\t f)^*(X)$. To define
$\fF_{\coh(P)}(f)X$ we {\it choose} an inverse image
$(\id_P\t f)^*(X)$ in an arbitrary way for each $X\in
\fF_{\coh(P)}(U)$, using the axiom of choice. Let
$f:U\ra V$, $g:V\ra W$ be scheme morphisms, and
$X\in\fF_{\coh(P)}(W)$. Then $(\id_P\t f)^*\bigl((\id_P\t g)^*
(X)\bigr)$ and $\bigl(\id_P\t(g\ci f)\bigr)^*(X)$ are both
inverse images of $X$ on $P\t U$, which are canonically
isomorphic, but {\it may not\/} be the same.

That is, $\fF_{\coh(P)}(f)\ci\fF_{\coh(P)}(g)(X)$ and
$\fF_{\coh(P)}(g\ci f)(X)$ may be different, but there
is a canonical isomorphism $\ep_{g,f}(X):\fF_{\coh(P)}(f)
\ci\fF_{\coh(P)}(g)(X)\ra\fF_{\coh(P)}(g\ci f)(X)$ for all
$X\in\fF_{\coh(P)}(W)$. These make up an isomorphism
of functors $\ep_{g,f}:\fF_{\coh(P)}(f)\ci\fF_{\coh(P)}(g)
\ra\fF_{\coh(P)}(g\ci f)$, that is, a 2-morphism in
$\excat$. From \S\ref{aa23}, these 2-morphisms $\ep_{g,f}$
are the last piece of data we need to define the 2-functor
$\fF_{\coh(P)}$. It is straightforward to show the
definition of a {\it contravariant\/ $2$-functor}
\cite[App.~B]{Gome} is satisfied.
\label{aa9ex1}
\end{ex}

Here are some remarks on this definition:
\begin{itemize}
\setlength{\itemsep}{0pt}
\setlength{\parsep}{0pt}
\item The most obvious way to define $\fF_{\coh(P)}(U)$ would
be to use {\it coherent sheaves} on $P\t U$ flat over $U$.
However, this turns out to be a bad idea, as coherence is
not well-behaved over non-noetherian schemes.

In particular, {\it inverse images} of coherent sheaves may
not be coherent, so the functors $\fF_{\coh(P)}(f)$ would
not be well-defined. Instead, we use {\it finitely
presentable quasicoherent sheaves}, which are the same
as coherent sheaves on noetherian $\K$-schemes, and behave
well under inverse images, etc.
\item If $f:U\ra V$ is not flat then $(\id_P\t f)^*:
\qcoh(P\t V)\ra\qcoh(P\t U)$ is {\it not\/} an exact functor.
We only claim above that exact sequences on $P\t V$ {\it flat
over} $V$ lift to exact sequences on~$P\t U$.
\end{itemize}

We supposed $P$ {\it smooth\/} in Example \ref{aa9ex1} to
make the {\it Chern character} well-defined for coherent
sheaves on $P$. For $P$ not smooth there may be problems
with this, so we need a different choice for $K(\coh(P))$.
We cannot take $K(\coh(P))=K_0(\coh(P))$, as Assumption
\ref{aa7ass}(iii) would not hold. Instead, in the next
example we define $K(\coh(P))$ using {\it Hilbert polynomials}.

\begin{ex} Let $\K$ be an algebraically closed field,
$P$ a projective $\K$-scheme, not necessarily smooth,
and $\O_P(1)$ a {\it very ample invertible sheaf\/}
on $P$, so that $(P,\O_P(1))$ is a {\it polarized\/
$\K$-scheme}. Following \cite[\S 1.2]{HuLe} and \cite[Ex.s
III.5.1 \& III.5.2]{Hart}, define the {\it Hilbert
polynomial\/} $p_X$ for $X\in\coh(P)$ by
\begin{equation*}
p_X(n)=\ts\sum_{i=0}^{\dim P}(-1)^i\dim_\K H^k\bigl(P,X(n)\bigr)
\quad\text{for $n\in\Z$,}
\end{equation*}
where $X(n)=X\ot\O_P(1)^n$, and $H^*(P,\cdot)$ is sheaf
cohomology on $P$. Then
\e
p_X(n)=\ts\sum_{i=0}^{\dim P}a_in^i/i!
\quad\text{for $a_0,\ldots,a_{\dim P}\in\Z$,}
\label{aa9eq1}
\e
by \cite[p.~10]{HuLe}. So $p_X(t)$ is a polynomial with
rational coeffients, written $p_X(t)\in\Q[t]$, with
degree no more than~$\dim P$.

If $0\!\ra\!X\!\ra\!Y\!\ra\!Z\!\ra\!0$ is exact then the
long exact sequence in sheaf cohomology implies that
$p_Y=p_X+p_Z$. Therefore the map $X\mapsto p_X$ factors
through the Grothendieck group $K_0(\coh(P))$, and there
is a unique group homomorphism $p:K_0(\coh(P))\ra\Q[t]$
with $p([X])=p_X$ for all $X\in\coh(P)$. Set $\A=\coh(P)$,
and let $K(\A)$ be the quotient of $K_0(\coh(P))$ by the
kernel of $p$. Then $K(\A)$ is isomorphic to the image of
$p$ in $\Q[t]$, which by \eq{aa9eq1} lies in a sublattice
of $\Q[t]$ isomorphic to $\Z^{1+\dim P}$. Now define
$\fF_{\coh(P)}$ as in Example \ref{aa9ex1}. Since this
does not use $K(\A)$ or the fact that $P$ is smooth,
no changes are needed, and $\fF_\A:\Sch_\K\ra\excat$
is a {\it contravariant\/ $2$-functor}.
\label{aa9ex2}
\end{ex}

\subsection{Verifying Assumption \ref{aa7ass}}
\label{aa92}

We now show that the examples of \S\ref{aa91} satisfy
Assumption \ref{aa7ass}. The main point to verify is
that $\fF_{\coh(P)}$ is a {\it stack in exact categories}.
The proof uses results of Grothendieck \cite{Grot3}, as in
Laumon and Moret-Bailly~\cite[\S 3.4.4]{LaMo}.

\begin{thm} Examples \ref{aa9ex1} and \ref{aa9ex2} satisfy
Assumption~\ref{aa7ass}.
\label{aa9thm1}
\end{thm}

\begin{proof} Let $\K$ be an algebraically closed field,
and $P$ a projective $\K$-scheme, not necessarily smooth.
We first verify the condition that if $X\in\A$ and $[X]=0$
in $K(\A)$ then $X\cong 0$. Let $K(\coh(P))$ be as in
either Example \ref{aa9ex1} or Example \ref{aa9ex2},
and suppose $X\in\coh(P)$ with $[X]=0$ in $K(\A)$. In
both cases, this implies the Hilbert polynomial $p_X$
of $X$ is zero. Now Serre's vanishing theorem shows
that for $n\gg 0$, the tautological map $H^0(X(n))\ot\O_P
(-n)\ra X$ is surjective, and $\dim H^0(X(n))=p_X(n)=0$,
so~$X\cong 0$.

We have shown $\fF_{\coh(P)}$ is a {\it contravariant\/
$2$-functor} $\Sch_\K\ra\excat$. We shall prove it is a
{\it stack in exact categories}, that is, that Definition
\ref{aa2def7}(i)--(iii) hold. Grothendieck
\cite[Cor.~VIII.1.2]{Grot3} proves Definition
\ref{aa2def7}(i),(ii) hold, using only the assumption
that the sheaves involved are {\it quasicoherent}.

Let $\{f_i:U_i\ra V\}_{i\in I}$ be an open cover of $V$ in the site
$\Sch_\K$, and let $X_i\in\Obj(\fF_{\coh(P)}(U_i))$ and $\phi_{ij}:
\fF_{\coh(P)}(f_{ij,j})X_j\ra\fF_{\coh(P)}(f_{ij,i})X_i$ be as in
Definition \ref{aa2def7}(iii). Grothendieck \cite[Cor.~VIII.1.3]{Grot3}
constructs $X\in\qcoh(P\t V)$ and isomorphisms $\phi_i:(\id_P\t f_i)^*
(X)\ra X_i$ in $\Mor(\qcoh(P\t U_i))$ satisfying Definition
\ref{aa2def7}(iii). This $X$ is finitely presentable by
\cite[Prop.~VIII.1.10]{Grot3} and flat over $V$ by
\cite[Cor.~IV.2.2.11(iii)]{Grot2}, so $X\in\fF_{\coh(P)}(V)$,
Definition \ref{aa2def7}(iii) holds for $\fF_{\coh(P)}$, and
$\fF_{\coh(P)}$ is a stack in exact categories.

It remains to verify Assumption \ref{aa7ass}(i)--(iii).
As $P\cong P\t\Spec\K$ is {\it noetherian}, Hartshorne
\cite[Prop.~II.5.7]{Hart} implies that $X\in\qcoh(P\t\Spec\K)$
is coherent if and only if it is finitely presentable, and
flatness over $\Spec\K$ is trivial. Hence
$\fF_{\coh(P)}(\Spec\K)=\coh(P)$, identifying $P$ and
$P\t\Spec\K$, and Assumption \ref{aa7ass}(i) holds. Part
(ii) follows from~\cite[Prop.~IV.2.2.7]{Grot2}.

Let $U\in\Sch_\K$ and $X\in\Obj(\fF_{\coh(P)}(U))$, and write
$X_u=\fF_{\coh(P)}(u)X$ in $\coh(P)$ for $u\in\Hom(\Spec\K,U)$.
Then we can regard $X$ as a flat family of $X_u\in\coh(P)$,
depending on $u\in\Hom(\Spec\K,U)$. But Chern classes in Example \ref{aa9ex1},
and by \cite[Prop.~III.7.9.11]{Grot2} Hilbert polynomials in
Example \ref{aa9ex2}, are both locally constant in flat
families. So Assumption \ref{aa7ass}(iii) holds.

Finally, let $X,Y\in{\coh(P)}$. Choose a basis $e_1,\ldots,e_n$
for $\Hom(X,Y)$, and let $z_1,\ldots,z_n:\Hom(X,Y)\ra\K$
be the corresponding coordinates. Regarding $\Hom(X,Y)$
as an {\it affine $\K$-scheme}, $z_1,\ldots,z_n$ become
sections of $\O_{\Hom(X,Y)}$. Write $\pi:\Hom(X,Y)\!\ra\!\Spec\K$,
$\pi_1:P\!\t\!\Hom(X,Y)\!\ra\!P$ and $\pi_2:P\!\t\!\Hom(X,Y)\!\ra\!
\Hom(X,Y)$ for the natural projections. Define
\begin{equation*}
\th_{X,Y}=\ts\sum_{i=1}^n\pi_2^*(z_i)\cdot\pi_1^*(e_i):
\pi_1^*(X)\longra\pi_1^*(Y).
\end{equation*}

Here $\pi_1^*(X)=\fF_{\coh(P)}(\pi)X$ and $\pi_1^*(Y)=\fF_{\coh(P)}(\pi)Y$
are coherent sheaves on $P\t\Hom(X,Y)$, $\pi_1^*(e_i):
\pi_1^*(X)\ra\pi_1^*(Y)$ is a morphism of sheaves on
$P\t\Hom(X,Y)$, $\pi_2^*(z_i)$ is a section of
$\O_{P\t\Hom(X,Y)}$, and `$\,\cdot\,$' multiplies sections
of $\O_{P\t\Hom(X,Y)}$ and sheaf morphisms on $P\t\Hom(X,Y)$.
Thus $\th_{X,Y}:\fF_{\coh(P)}(\pi)X\ra\fF_{\coh(P)}(\pi)Y$ is
a morphism in $\fF_{\coh(P)}(\Hom(X,Y))$, and clearly \eq{aa7eq1}
holds for all $f\in\Hom(X,Y)$. This verifies
Assumption~\ref{aa7ass}(iv).
\end{proof}

\subsection{Showing stacks are algebraic and locally of finite type}
\label{aa93}

Here and in \S\ref{aa94} we show that Examples \ref{aa9ex1},
\ref{aa9ex2} satisfy Assumption~\ref{aa8ass}.

\begin{thm} In Examples \ref{aa9ex1} and \ref{aa9ex2},
$\fObj_{\coh(P)}^\al,\fExact_{\coh(P)}^{\smash{\al,\be,\ga}}$
are algebraic $\K$-stacks, locally of finite type, for
all\/~$\al,\be,\ga$.
\label{aa9thm2}
\end{thm}

\begin{proof} Laumon and Moret-Bailly \cite[Th.~4.6.2.1]{LaMo}
prove $\fObj_{\coh(P)}$ is an algebraic $\K$-stack, locally of
finite type. The corresponding result for $\fObj_{\coh(P)}^\al$
then follows from Theorem \ref{aa7thm1}. An important part of the
proof is to construct an {\it atlas} for $\fObj_{\coh(P)}$, and
we sketch how this is done using {\it Quot-schemes}.

For $N\ge 0$, write $\Quot_P(\op^N\O_P)$ for {\it
Grothendieck's Quot-scheme} \cite[\S 3.2]{Grot1},
\cite[\S 2.2]{HuLe}. This is the {\it moduli scheme of
quotient sheaves} $(Z,\psi)$ of the coherent sheaf $\op^N\O_P$
on $P$, where $\psi:\op^N\O_P\ra Z$ is surjective. By Assumption
\ref{aa7ass}(iii) it may be written as a disjoint union
\begin{equation*}
\ts\Quot_P(\op^N\O_P)=\coprod_{\al\in \bar C(\coh(P))}
\ts\Quot_P(\op^N\O_P,\al),
\end{equation*}
where $\Quot_P(\op^N\O_P,\al)$ is the subscheme of
$(Z,\psi)$ with $[Z]=\al$ in~$K(\coh(P))$.

Considered as a $\K$-stack, for $U\in\Sch_\K$,
$\Quot_P(\op^N\O_P,\al)(U)$ has {\it objects\/} pairs
$(Z,\psi)$ for $Z\in\fObj_{\coh(P)}^\al(U)$ and
$\psi:\op^N\O_{P\t U}\ra Z$ a surjective morphism in
$\fF_{\coh(P)}(U)$, that is, the r.h.s.\ of some short
exact sequence. {\it Morphisms} $\be:(Z,\psi)\ra(Z',\psi')$
are isomorphisms $\be:Z\ra Z'$ in $\fF_{\coh(P)}(U)$
with~$\psi'=\be\ci\psi$.

Grothendieck \cite[\S 3.2]{Grot1} shows $\Quot_P(\op^N\O_P,\al)$
is represented by a projective $\K$-scheme. Let $\Quot_P^\ci
(\op^N\O_P,\al)$ be the open $\K$-substack of $\Quot_P(\op^N
\O_P,\al)$ with objects $(Z,\psi)$ over $U\in\Sch_\K$ such that
$(\pi_U)_*(\psi):\op^N\O_U\ra(\pi_U)_*(Z)$ is an isomorphism and
$R^p(\pi_U)_*(Z)=0$ for all $p>0$, in the notation of
\cite[\S III.8]{Hart}, where $\pi_U:P\t U\ra U$ is the projection.
Then $\Quot_P^\ci(\op^N\O_P,\al)$ is represented by a quasiprojective
$\K$-scheme, and so is of finite type.

Write $\al\mapsto\al(n)$ for the unique automorphism
of $K(\coh(P))$ such that if $X\in\coh(P)$ with $[X]=\al$
in $K(\coh(P))$, then $[X(n)]=\al(n)$. For all $N,n\ge 0$,
define a 1-morphism $\Pi_{N,n}:\Quot_P^\ci(\op^N\O_P,
\al(n))\ra\fObj_{\coh(P)}^\al$ by $\Pi_{N,n}:(Z,\psi)
\mapsto Z(-n)$ on objects and $\Pi_{N,n}:\be\mapsto\be(-n)$
on morphisms, where $Z(-n)$ and $\be(-n):Z(-n)\ra Z'(-n)$
are the twists of $Z$ and $\be:Z\ra Z'$ by the lift of the
invertible sheaf $\O_P(-n)$ to $P\t U$. Define a 1-morphism
\begin{equation*}
\Pi^\al=\ts\coprod_{N,n\ge 0}\Pi_{N,n}:\coprod_{N,n\ge 0}\Quot_P^\ci
(\op^N\O_P,\al(n))\longra\fObj_{\coh(P)}^\al.
\end{equation*}
As in \cite[p.~30]{LaMo}, $\Pi^\al$ is smooth and surjective,
so it is an atlas for $\fObj_{\coh(P)}^\al$. Since
$\coprod_{N,n\ge 0}\Quot_P^\ci(\op^N\O_P,\al(n))$ is
locally of finite type, so is~$\fObj_{\coh(P)}^\al$.

The proof for $\fExact_{\coh(P)}^{\smash{\al,\be,\ga}}$ is mostly
a straightforward generalization of that for $\fObj_{\coh(P)}$. The
difficult part is to construct an atlas for
$\fExact_{\coh(P)}^{\smash{\al,\be,\ga}}$, which we do by a method
explained to me by Bernd Siebert. The important point is that
Quot-schemes work not just for $\K$-schemes, but for
$S$-{\it schemes} over a general locally noetherian base scheme~$S$.

As above $\Quot_P^\ci(\op^N\O_P,\be)$ is represented by a
quasiprojective $\K$-scheme $S_N^\be$. Thus the groupoid
$\Quot_P(\op^N\O_P,\be)(S_N^\be)$ is equivalent to
$\Hom(S_N^\be,S_N^\be)$ in $\Sch_\K$. Let $(Y_N^\be,\xi_N^\be)$
in $\Quot_P(\op^N\O_P,\be)(S_N^\be)$ be identified with
$\id_{S_N^\be}$. Then $(Y_N^{\smash{\be}},\xi_N^{\smash{\be}})$ is
the {\it universal quotient sheaf\/} of $\op^N\O_P$ with class~$\be$.

Now $Y_N^\be$ is a quasicoherent sheaf on $P\t S_N^\be$,
flat over $S_N^\be$. Regard $P\t S_N^\be$ with the
projection $P\t S_N^\be\ra S_N^\be$ as a {\it projective
$S_N^\be$-scheme}. Then $Y_N^\be$ becomes a {\it coherent
sheaf\/} on the $S_N^\be$-scheme $P\t S_N^\be$. By
\cite[\S 3.2]{Grot1} we can therefore form the
Quot-scheme $\Quot_{P\t S_N^\be/S_N^\be}(Y_N^\be,\ga)$.
It is the {\it moduli stack of quotient sheaves of\/
$Y_N^{\smash{\be}}$ on the $S_N^{\smash{\be}}$-scheme
$P\t S_N^{\smash{\be}}$ with class} $\ga\in K(\coh(P))$, and
is {\it represented by a projective $S_N^{\smash{\be}}$-scheme}.

Interpreting $\Quot_{P\t S_N^\be/S_N^\be}(Y_N^\be,\ga)$ as a
$\K$-stack rather than an $S_N^\be$-stack, $\Quot_{P\t S_N^\be
/S_N^\be}(Y_N^{\smash{\be}},\ga)(U)$ for $U\in\Sch_\K$ has
{\it objects} $(g,Z,\psi)$, where $g:U\ra S_N^\be$ is a
morphism in $\Sch_\K$, $Z\in\fObj_{\coh(P)}^\ga(U)$ and
$\psi:\fF_{\coh(P)}(g)Y_N^\be\ra Z$ is surjective in
$\Mor(\fF_{\coh(P)}(U))$, and {\it morphisms\/} $\tau:
(g,Z,\psi)\ra(g,Z',\psi')$, where $\tau:Z\ra Z'$ is an
isomorphism in $\Mor(\fF_{\coh(P)}(U))$ with~$\psi'=\tau\ci\psi$.

For $N,n\ge 0$, define $\Pi_{N,n}:\Quot_{P\t S_N^{\be(n)}/
S_N^{\be(n)}}(Y_N^{\be(n)},\ga(n))\!\ra\!\fExact_{\coh(P)}^{
\smash{\al,\be,\ga}}$ on {\it objects} by $\Pi_{N,n}(U):
(g,Z,\psi)\!\mapsto\!(X,(\fF_{\coh(P)}(g)Y_N^{\smash{\be(n)}})
(-n),Z(-n),\phi,\psi(-n))$, where $\phi:X\!\ra\!(\fF_{\coh(P)}(g)
Y_N^{\be(n)})(-n)$ is a kernel for $\psi(-n)$. This involves a
{\it choice} of $X,\phi$, but one which is {\it unique up to
canonical isomorphism}. If $\tau:(g,Z,\psi)\ra(g,Z',\psi')$ in
$\Quot_{P\t S_N^\be/S_N^\be}(Y_N^\be,\ga)$ and $X,\phi$, $X',\phi'$
are choices made for $\Pi_{N,n}(U)(g,Z,\psi),\Pi_{N,n}(U)(g,Z',\psi')$,
we define $\Pi_{N,n}$ on {\it morphisms} by $\Pi_{N,n}(U):\tau\mapsto
(\ze,\id_{(\fF_{\coh(P)}(g)Y_N^{\be(n)})(-n)},\tau(-n))$,
where $\ze:X\ra X'$ is the unique isomorphism with
$\phi=\phi'\ci\ze$. Then $\Pi_{N,n}$ is a 1-morphism. Define
\begin{equation*}
\Pi^{\al,\be,\ga}\!=\!\ts\coprod_{N,n\ge 0}\Pi_{N,n}:\coprod_{N,n\ge 0}
\Quot_{P\t S_N^{\be(n)}/S_N^{\be(n)}}\bigl(Y_N^{\be(n)},\ga(n)
\bigr)\!\ra\!\fExact_{\coh(P)}^{\al,\be,\ga}.
\end{equation*}
As in \cite[p.~30]{LaMo}, $\Pi^{\al,\be,\ga}$ is smooth and
surjective, so it is an atlas.
\end{proof}

\subsection{Showing $1$-morphisms are of finite type}
\label{aa94}

Next we prove that the 1-morphisms of \eq{aa8eq1} are of
finite type. For the 1-morphisms $\fm$ in the next
proposition, this is because the fibres of $\fm$ are
essentially Quot-schemes of quotient sheaves with
fixed Hilbert polynomial. Thus by Grothendieck's
construction they are {\it projective $\K$-schemes}.

\begin{prop} In Examples \ref{aa9ex1} and \ref{aa9ex2},
$\fm:\fExact^{\smash{\al,\be,\ga}}_{\coh(P)}\ra\fObj_{\coh(P)}^{
\smash{\be}}$ is a finite type $1$-morphism.
\label{aa9prop2}
\end{prop}

\begin{proof} Consider the commutative diagram
\e
\begin{gathered}
\xymatrix@R=15pt{
\coprod_{N,n\ge 0}\Quot_{P\t S_N^{\be(n)}/S_N^{\be(n)}}
\bigl(X_N^{\be(n)},\ga(n)\bigr)
\ar[rr]_{\qquad\qquad\Pi^{\al,\be,\ga}}
\ar[d]_{\coprod_{N,n\ge 0}\pi_{N,n}^{\be,\ga}}
&& \fExact_{\coh(P)}^{\al,\be,\ga} \ar[d]^\fm \\
\coprod_{N,n\ge 0}S_N^{\be(n)}\!=\!\coprod_{N,n\ge 0}
\Quot_P^\ci\bigl(\op^N\O_P,\be(n)\bigr)
\ar[rr]^{\qquad\qquad\qquad\Pi^\be}&& \fObj_{\coh(P)}^\be,
}
\end{gathered}
\label{aa9eq2}
\e
where $\pi_{N,n}^{\be,\ga}$ is the morphism making
$\Quot_{\smash{P\t S_N^{\be(n)}/S_N^{\be(n)}}}
\bigl(X_N^{\be(n)},\ga(n)\bigr)$ into an $S_N^{\be(n)}$-scheme.
It is {\it projective}, as $\Quot_{\smash{P\t S_N^{\be(n)}/
S_N^{\be(n)}}}\bigl(X_N^{\be(n)},\ga(n)\bigr)$ is projective
over $S_N^{\be(n)}$. Thus $\pi_{N,n}^{\be,\ga}$ is of finite
type, and so is $\coprod_{N,n\ge 0}\pi_{N,n}^{\be,\ga}$ in
\eq{aa9eq2}. As $\Pi^{\al,\be,\ga},\Pi^\be$ are atlases,
\eq{aa9eq2} implies $\fm$ is of finite type.
\end{proof}

For the 1-morphism $\fb\t\fe$ of \eq{aa8eq1}, the fibre
over $(X,Z)$ in $\coh(P)\t\coh(P)$ is the stack of
{\it isomorphism classes of exact sequences}
$0\ra X\,\smash{{\buildrel\phi\over\longra}\,Y\,
{\buildrel\psi\over\longra}}\,Z\ra 0$ in $\coh(P)$.
Such sequences are classified by $\Ext^1(Z,X)$, so the
fibre of $\fb\t\fe$ over $(X,Z)$ should be the {\it
quotient stack}~$[\Ext^1(Z,X)/\Aut(X)\t\Aut(Z)]$.

As $\Ext^1(Z,X)$ is a finite-dimensional $\K$-vector
space, this fibre is finite type, so $\fb\t\fe$ should
be of finite type. The proof below does not use this
argument, but depends on facts about Quot-schemes
which encode the same ideas.

\begin{thm} In Examples \ref{aa9ex1} and \ref{aa9ex2}, the
$1$-morphism $\fb\t\fe:\fExact^{\smash{\al,\be,\ga}}_{\coh(P)}\ab
\ra\ab
\fObj_{\coh(P)}^\al\t\fObj_{\coh(P)}^\ga$ is of finite type.
\label{aa9thm3}
\end{thm}

\begin{proof} Use the notation of Theorem \ref{aa9thm2}. Then
$\fObj_{\coh(P)}^\al\t\fObj_{\coh(P)}^\ga$ is covered by open
substacks of the form
\e
V_{M,N,n}^{\al,\ga}\!=\!(\Pi_{M,n}^\al\!\t\!\Pi_{N,n}^\ga)
\bigl(\Quot_P^\ci(\op^M\O_P,\al(n))\!\t\!
\Quot_P^\ci(\op^N\O_P,\ga(n))\bigr),
\label{aa9eq3}
\e
for $M,N,n\ge 0$. Let $W_{M,N,n}^{\smash{\al,\be,\ga}}=
(\fb\t\fe)^{-1}(V_{M,N,n}^{\smash{\al,\ga}})$ be the inverse
image of $V_{M,N,n}^{\smash{\al,\ga}}$ in $\fExact^{\smash{
\al,\be,\ga}}_{\coh(P)}$ under $\fb\t\fe$. As $\Quot_P^\ci
(\op^M\O_P,\al(n))\!\t\!\Quot_P^\ci(\op^N\O_P,\ga(n))$ is an
atlas for $V_{M,N,n}^{\smash{\al,\ga}}$ represented by a
quasiprojective scheme, $V_{M,N,n}^{\smash{\al,\ga}}$ is
of finite type. Also $S_{M+N}^{\smash{\be(n)}}$ is a
quasiprojective scheme representing $\Quot_P^\ci(\op^{M+N}
\O_P,\be(n))$, and $(Y_{M+N}^{\smash{\be(n)}},\xi_{M+N}^{\smash{
\be(n)}})$ is the universal quotient sheaf of $\op^N\O_P$
on~$S_{M+N}^{\smash{\be(n)}}$.

Form the Quot-scheme $Q_{M+N}^{\smash{\be(n),\ga(n)}}=
\Quot_{P\t S_{M+N}^{\smash{\be(n)}}/S_{M+N}^{\smash{\be(n)}}}
(Y_{M+N}^{\smash{\be(n)}},\ga(n))$, as in Theorem \ref{aa9thm2}.
It is projective over $S_{M+N}^{\smash{\be(n)}}$, and
$S_{M+N}^{\smash{\be(n)}}$ is of finite type, so
$Q_{M+N}^{\smash{\be(n),\ga(n)}}$ is of finite type.
We have a projection
\e
\Pi_{M+N,n}^{\al,\be,\ga}:Q_{M+N}^{\be(n),\ga(n)}
\longra\fExact_{\coh(P)}^{\al,\be,\ga},
\label{aa9eq4}
\e
forming part of an atlas for $\fExact_{\coh(P)}^{\smash{\al,\be,\ga}}$.
We shall show \eq{aa9eq4} covers~$W_{M,N,n}^{\smash{\al,\be,\ga}}$.

Let $U\in\Sch_\K$, and let $(X,Y,Z,\phi,\psi)\in W_{M,N,n}^{
\smash{\al,\be,\ga}}(U)$. By definition this means that
$(X,Y,Z,\phi,\psi)\in\fExact_{\coh(P)}^{\smash{\al,\be,\ga}}
(U)$ and $(X,Z)\in V_{M,N,n}^{\smash{\al,\ga}}(U)$. As
$V_{M,N,n}^{\smash{\al,\ga}}$ is defined as an {\it image}
in \eq{aa9eq3}, roughly speaking this means that $X$ lifts
to $\Quot_P^\ci(\op^M\O_P,\al(n))(U)$ and $Z$ lifts
to~$\Quot_P^\ci(\op^N\O_P,\ga(n))(U)$.

However, these lifts need exist only {\it locally} in the
\'etale topology on $U$. That is, there exists an open cover
$\{f_i:U_i\ra U\}_{i\in I}$ of $U$ in the site $\Sch_\K$ and
objects $(X_i,\la_i)\in\Quot_P^\ci(\op^M\O_P,\al(n))(U_i)$
and $(Z_i,\nu_i)\in\Quot_P^\ci(\op^N\O_P,\ga(n))(U_i)$
with $X_i(-n)=\fF_{\coh(P)}(f_i)(X)$, $Z_i(-n)=\fF_{\coh(P)}
(f_i)(Z)$ for all $i\in I$. Set $Y_i\!=\!(\fF_{\coh(P)}(f_i)Y)(n)$,
$\phi_i\!=\!(\fF_{\coh(P)}(f_i)\phi)(n)$, $\psi_i\!=\!(\fF_{\coh(P)}
(f_i)\psi)(n)$.

Refining the cover $\{f_i:U_i\ra U\}_{i\in I}$ if necessary,
we can construct $\mu_i$ to make a commutative diagram
\begin{equation*}
\xymatrix@R=15pt{
0 \ar[r]
& \op^M\O_{P\t U_i} \ar[r] \ar[d]^{\la_i}
& \op^{M+N}\O_{P\t U_i} \ar[r] \ar[d]^{\mu_i}
& \op^N\O_{P\t U_i} \ar[r] \ar[d]^{\nu_i} & 0 \\
0 \ar[r] & X_i \ar[r]^{\phi_i} & Y_i \ar[r]^{\psi_i}
& Z_i \ar[r] & 0
}
\end{equation*}
in $\qcoh(P\t U_i)$, with exact rows and surjective columns,
such that $(Y_i,\mu_i)\in\Quot_P^\ci(\op^{M+N}
\O_P,\be(n))(U_i)$, for all~$i\in I$.

As the $\K$-scheme $S_{M+N}^{\smash{\be(n)}}$ represents
$\Quot_P^\ci(\op^{M+N}\O_P,\be(n))$ with universal quotient
sheaf $(Y_{M+N}^{\smash{\be(n)}},\xi_{M+N}^{\smash{\be(n)}})$,
there is a unique morphism $g_i:U_i\ra S_{M+N}^{\smash{\be(n)}}$
and an isomorphism $\eta_i:Y_i\ra\fF_{\coh(P)}(g_i)
Y_{M+N}^{\smash{\be(n)}}$ with $\eta_i\ci\mu_i=\fF_{\coh(P)}
(g_i)\xi_{M+N}^{\smash{\be(n)}}$. One can then show that
the image of $(g_i,Z_i,\psi_i\ci\eta_i^{-1})$ under
$\Pi_{N,n}^{\al,\be,\ga}(U_i)$ is isomorphic to
$(X_i(-n),Y_i(-n),Z_i(-n),\phi_i(-n),\psi_i(-n))=
\fExact_{\coh(P)}^{\smash{\al,\be,\ga}}(f_i)(X,Y,Z,\phi,\psi)$.

We have shown that for any $U\in\Sch_\K$ and $(X,Y,Z,\phi,\psi)\in
W_{M,N,n}^{\smash{\al,\be,\ga}}(U)$, there exists an open cover
$\{f_i:U_i\ra U\}_{i\in I}$ of $U$ in the site $\Sch_\K$ and
objects $(g_i,Z_i,\psi_i\ci\eta_i^{-1})\in Q_{M+N}^{\smash{
\be(n),\ga(n)}}(U_i)$ such that
\begin{equation*}
\Pi_{N,n}^{\be,\ga}(U_i)(g_i,Z_i,\psi_i\ci\eta_i^{-1})\cong
\fExact_{\coh(P)}^{\al,\be,\ga}(f_i)(X,Y,Z,\phi,\psi)
\end{equation*}
in $\fExact_{\coh(P)}^{\al,\be,\ga}(U_i)$, for all $i\in I$. This
proves that \eq{aa9eq4} covers~$W_{M,N,n}^{\smash{\al,\be,\ga}}$.

But from above $Q_{M+N}^{\smash{\be(n),\ga(n)}}$ is
of finite type, so $W_{M,N,n}^{\smash{\al,\be,\ga}}$
is of finite type, and therefore the 1-morphism
$\fb\t\fe:W_{M,N,n}^{\smash{\al,\be,\ga}}\ra
V_{M,N,n}^{\smash{\al,\ga}}$ is of finite type. As the
$V_{M,N,n}^{\smash{\al,\ga}}$ cover $\fObj_{\coh(P)}^\al
\t\fObj_{\coh(P)}^\ga$ and $W_{M,N,n}^{\smash{\al,\be,\ga}}$
is the inverse image of $V_{M,N,n}^{\smash{\al,\ga}}$ under
$\fb\t\fe$, this shows $\fb\t\fe$ is of finite type,
and the proof is complete.
\end{proof}

The last three results now prove:

\begin{thm} Examples \ref{aa9ex1} and \ref{aa9ex2} satisfy
Assumption~\ref{aa8ass}.
\label{aa9thm4}
\end{thm}

Theorems \ref{aa9thm1} and \ref{aa9thm4} show that we may
apply the results of \S\ref{aa7} and \S\ref{aa8} to Examples
\ref{aa9ex1} and \ref{aa9ex2}. This yields large classes of
{\it moduli stacks of\/ $(I,\pr)$-configurations of coherent
sheaves} $\fM(I,\pr)_{\coh(P)},\fM(I,\pr,\ka)_{\coh(P)}$ on
a projective $\K$-scheme $P$, which are {\it algebraic
$\K$-stacks, locally of finite type}. It also gives many
1-{\it morphisms} $S(I,\pr,J)$, $Q(I,\pr,K,\tl,\phi),\ldots$
between these moduli stacks, various of which are
{\it representable} or {\it of finite type}.

\section{Representations of quivers and algebras}
\label{aa10}

Finally we consider configurations in some more large classes
of examples of abelian categories, {\it representations of
quivers} $Q$ and of {\it finite-dimensional algebras} $A$.
After introducing quivers and their representations in
\S\ref{aa101}, section \ref{aa102} defines the data $\A,K(\A),
\fF_\A$ of Assumption \ref{aa7ass} for five related families of
examples. Sections \ref{aa103} and \ref{aa104} then prove that
Assumptions \ref{aa7ass} and \ref{aa8ass} hold for each case,
so that the results of \S\ref{aa7} and \S\ref{aa8} apply.

\subsection{Introduction to quivers}
\label{aa101}

Here are the basic definitions in quiver theory, taken from
Benson \cite[\S 4.1]{Bens}. We fix an algebraically closed
field $\K$ throughout.

\begin{dfn} A {\it quiver\/} $Q$ is a finite directed graph.
That is, $Q$ is a quadruple $(Q_0,Q_1,b,e)$, where $Q_0$ is
a finite set of {\it vertices}, $Q_1$ is a finite set of
{\it arrows}, and $b,e:Q_1\ra Q_0$ are maps giving the
{\it beginning} and {\it end\/} of each arrow.

The {\it path algebra} $\K Q$ is an associative algebra
over $\K$ with basis all {\it paths of length\/} $k\ge 0$,
that is, sequences of the form
\e
v_0\,{\buildrel a_1\over\longra}\, v_1\ra\cdots\ra
v_{k-1}\,{\buildrel a_k\over\longra}\,v_k,
\label{aa10eq1}
\e
where $v_0,\ldots,v_k\in Q_0$, $a_1,\ldots,a_k\in Q_1$,
$b(a_i)=v_{i-1}$ and $e(a_i)=v_i$. Multiplication is given
by composition of paths {\it in reverse order}.

Each $v\in Q_0$ determines a basis element \eq{aa10eq1} with
$k=0$, $v_0=v$, and the identity in $\K Q$ is $1=\sum_{v\in Q_0}v$.
Each $a\in Q_1$ determines a basis element $b(a)\,{\smash{\buildrel a
\over\longra}}\,e(a)$ in $\K Q$ with $k=1$. For brevity we refer to
this element as $a$. Note that $\K Q$ is finite-dimensional if
and only if $Q$ has {\it no oriented cycles}.

For $n\ge 0$, write $\K Q_{(n)}$ for the vector subspace
of $\K Q$ with basis all paths of length $k\ge n$. Then
$\K Q_{(n)}$ is a {\it two-sided ideal\/} in $\K Q$. A
{\it quiver with relations} $(Q,I)$ is defined to
be a quiver $Q$ together with a two-sided ideal
$I$ in $\K Q$ such that $I\subseteq\K Q_{(2)}$. Then
$\K Q/I$ is an associative $\K$-algebra.
\label{aa10def1}
\end{dfn}

\begin{dfn} Let $Q=(Q_0,Q_1,b,e)$ be a quiver. A {\it
representation} of $Q$ consists of finite-dimensional
$\K$-vector spaces $X_v$ for each $v\in Q_0$, and linear
maps $\rho_a:X_{b(a)}\ra X_{e(a)}$ for each $a\in Q_1$.
Representations of $Q$ are in 1-1 correspondence with
{\it finite-dimensional left\/ $\K Q$-modules} $(X,\rho)$,
as follows.

Given $X_v,\rho_a$, define $X=\bigop_{v\in Q_0}X_v$,
and a linear $\rho:\K Q\ra\End(X)$ taking \eq{aa10eq1}
to the linear map $X\ra X$ acting as $\rho_{a_k}\ci
\rho_{a_{k-1}}\ci\cdots\ci\rho_{a_1}$ on $X_{v_0}$,
and 0 on $X_v$ for $v\ne v_0$. Then $(X,\rho)$ is
a left $\K Q$-module. Conversely, any such $(X,\rho)$
comes from a unique representation of $Q$, taking $X_v$
for $v\in Q_0$ to be the 1-eigenspace of $\rho(v)$ in $X$,
and $\rho_a$ for $a\in Q_1$ to be the restriction
of $\rho(a):X\ra X$ to~$X_{b(a)}$.

We generally write representations of $Q$ as left
$\K Q$-modules $(X,\rho)$. A {\it morphism of representations}
$\phi:(X,\rho)\ra(Y,\si)$ is a linear map $\phi:X\ra Y$ with
$\phi\ci\rho(\ga)=\si(\ga)\ci\phi$ for all $\ga\in\K Q$.
Equivalently, $\phi$ defines linear maps $\phi_v:X_v\ra Y_v$
for all $v\in Q_0$ with $\phi_{e(a)}\ci\rho_a=\si_a\ci\phi_{b(a)}$
for all~$a\in Q_1$.

A representation $(X,\rho)$ of $Q$ is called {\it nilpotent\/}
if $\rho(\K Q_{(n)})=\{0\}$ in $\End(X)$ for some $n\ge 0$. Let
$(Q,I)$ be a {\it quiver with relations}. A {\it representation} of
$(Q,I)$ is a representation $(X,\rho)$ of $Q$ with $\rho(I)=\{0\}$.
Then $X$ is a representation of the quotient algebra~$\K Q/I$.

Write $\modKQ$ for the {\it category of representations} of
$Q$, and $\nilKQ$ for the full subcategory of {\it nilpotent\/}
representations of $Q$. If $(Q,I)$ is a {\it quiver with
relations}, write $\modKQI$ for the category of representations
of $(Q,I)$, and $\nilKQI$ for the full subcategory of {\it
nilpotent\/} representations of $(Q,I)$. It is easy to show
all of these are abelian categories, of finite length. If $Q$
has {\it no oriented cycles} then $\modKQ=\nilKQ$, since
$\K Q_{(n)}=0$ for $n>\md{Q_1}$. If $\K Q_{(n)}\subseteq I$
for some $n\ge 2$ then~$\modKQI=\nilKQI$.
\label{aa10def2}
\end{dfn}

We consider the {\it Grothendieck groups} of~$\modKQ,\ldots,\nilKQI$.

\begin{dfn} Let $Q=(Q_0,Q_1,b,e)$ be a quiver and $(X,\rho)$
a representation of $Q$. Write $\N^{Q_0}$ and $\Z^{Q_0}$ for
the sets of maps $Q_0\!\ra\!\N$ and $Q_0\!\ra\!\Z$. Define the
{\it dimension vector} $\bdim(X,\rho)\!\in\!\N^{Q_0}\subset\Z^{Q_0}$
of $(X,\rho)$ by $\bdim(X,\rho):v\mapsto\dim_\K X_v$. This
induces a surjective group homomorphism $\bdim:K_0(\modKQ)
\!\ra\!\Z^{Q_0}$. The same applies to $\nilKQ,\ldots,\nilKQI$.
As $\modKQ,\ldots,\nilKQI$ have finite length $K_0(\modKQ),
\ldots,K_0(\nilKQI)$ are the free abelian groups with bases
isomorphism classes of simple objects in~$\modKQ,\ldots,\modKQI$.

A nilpotent representation $(X,\rho)$ is simple if $X_v\cong\K$
for some $v\in Q_0$, and $X_w=0$ for $w\ne v$, and $\rho_a=0$
for all $a\in Q_1$. So simple objects in $\nilKQ,\nilKQI$ up to
isomorphism are in 1-1 correspondence with $Q_0$, and $\bdim:
K_0(\nilKQ),K_0(\nilKQI)\ra\Z^{Q_0}$ is an isomorphism. When
$Q$ has {\it oriented cycles}, there are usually {\it many}
simple objects in $\modKQ,\modKQI$, and $K_0(\modKQ),K_0(\modKQI)$
are much larger than~$\Z^{Q_0}$.
\label{aa10def3}
\end{dfn}

Quivers are used to study the representations of
{\it finite-dimensional algebras}.

\begin{dfn} Let $A$ be a {\it finite-dimensional\/ $\K$-algebra},
and $\modA$ the category of {\it finite-dimensional left\/
$A$-modules} $(X,\rho)$, where $X$ is a finite-dimensional
$\K$-vector space and $\rho:A\ra\End(X)$ an algebra morphism.
Then $\modA$ is an abelian category of finite length.
Following Benson \cite[\S 2.2, Def.~4.1.6, Prop.~4.1.7]{Bens}
one defines a quiver with relations $(Q,I)$ called the
{\it Ext-quiver} of $A$, whose vertices $Q_0$ correspond to
isomorphism classes of simple objects in $\modA$, with a
natural {\it equivalence of categories} between $\modA$ and
$\modKQI$. So, the representations of $A$ can be understood
in terms of those of~$(Q,I)$.
\label{aa10def4}
\end{dfn}

\subsection{Definition of the data $\A,K(\A),\fF_\A$}
\label{aa102}

In five examples we define the data of Assumption \ref{aa7ass}
for the abelian categories $\modKQ,\nilKQ,\modKQI,\nilKQI,\modA$
of \S\ref{aa101}, respectively. The main ideas are all in Example
\ref{aa10ex1}, with minor variations in Examples
\ref{aa10ex2}--\ref{aa10ex5}. We fix an {\it algebraically closed
field\/} $\K$ throughout.

\begin{ex} Let $Q=(Q_0,Q_1,b,e)$ be a quiver. Take $\A=\modKQ$,
the abelian category of representations of $Q$. Define $K(\modKQ)$
to be the quotient of $K_0(\modKQ)$ by the kernel of $\bdim:
K_0(\modKQ)\ra\Z^{Q_0}$. Then $\bdim$ induces an isomorphism
$K(\modKQ)\cong\Z^{Q_0}$. We shall {\it identify} $K(\modKQ)$
and $\Z^{Q_0}$, so that for $X\in\modKQ$ the class $[X]$ in
$K(\modKQ)$ is~$\bdim X$.

Motivated by King \cite[Def.~5.1]{King}, for $U\in\Sch_\K$
define $\fF_\modKQ(U)$ to be the category with {\it objects}
$(X,\rho)$ for $X$ a locally free sheaf of finite rank on $U$
and $\rho:\K Q\ra\Hom(X,X)$ a $\K$-algebra homomorphism, and
{\it morphisms} $\phi:(X,\rho)\ra(Y,\si)$ to be morphisms
of sheaves $\phi:X\ra Y$ with $\phi\ci\rho(\ga)=\si(\ga)\ci\phi$
in $\Hom(X,Y)$ for all~$\ga\in\K Q$.

Now define $\A_U$ to be the category with objects
$(X,\rho)$ for $X$ a quasicoherent sheaf on $U$ and
$\rho:\K Q\ra\Hom(X,X)$ a $\K$-algebra homomorphism,
and morphisms $\phi$ as above. It is easy to show
$\A_U$ is an abelian category, and $\fF_\modKQ(U)$ an
exact subcategory of $\A_U$. Thus $\fF_\modKQ(U)$ is an
exact category.

For $f:U\ra V$ in $\Mor(\Sch_\K)$, define a functor
$\fF_\modKQ(f):\fF_\modKQ(V)\ra\fF_\modKQ(U)$ by
$\fF_\modKQ(f):(X,\rho)\mapsto(f^*(X),f^*(\rho))$
on objects $(X,\rho)$ and $\fF_\modKQ(f):\phi\mapsto f^*(\phi)$
on morphisms $\phi:(X,\rho)\ra(Y,\si)$, where $f^*(X)$ is
the {\it inverse image sheaf\/} and $f^*(\rho)(\ga)=
f^*(\rho(\ga)):f^*(X)\ra f^*(X)$ for $\ga\in\K Q$ and
$f^*(\phi):f^*(X)\ra f^*(Y)$ are pullbacks of morphisms
between inverse images.

Since $f^*(\O_V)\cong\O_U$, inverse images of locally free
sheaves of finite rank are also locally free of finite rank,
so $\fF_\modKQ(f)$ is a functor $\fF_\modKQ(V)\ra\fF_\modKQ(U)$.
As locally free sheaves on $V$ are flat over $V$, Grothendieck
\cite[Prop.~IV.2.1.8(i)]{Grot2} implies $\fF_\modKQ(f)$ is an
exact functor.

As in Example \ref{aa9ex1}, defining $\fF_\modKQ(f)(X,\rho)$
involves making a {\it choice} for $f^*(X)$, arbitrary up
to canonical isomorphism. When $f:U\ra V$, $g:V\ra W$ are
morphisms in $\Sch_\K$, the canonical isomorphisms yield
an isomorphism of functors $\ep_{g,f}:\fF_\modKQ(f)\ci
\fF_\modKQ(g)\ra\fF_\modKQ(g\ci f)$, that is, a 2-morphism
in $\excat$. The $\ep_{g,f}$ complete the definition of
$\fF_\modKQ$, and one can readily show that the definition
of a 2-functor \cite[App.~B]{Gome} holds.
\label{aa10ex1}
\end{ex}

Here is how to extend this to {\it nilpotent\/} representations.

\begin{ex} Take $\A=\nilKQ$, the abelian category of nilpotent
representations of a quiver $Q$. Then $K_0(\nilKQ)\cong\Z^{Q_0}$.
Set~$K(\nilKQ)=K_0(\nilKQ)$.

For $U\in\Sch_\K$ let $\fF_\nilKQ(U)$ be the full exact
subcategory of $\fF_\modKQ(U)$ with objects $(X,\rho)$
such that there exists an open cover $\{f_n:U_n\ra U\}_{n\in\N}$
of $U$ in the site $\Sch_\K$ for which $(X_n,\rho_n)=\fF_\modKQ
(f_n)(X,\rho)$ satisfies $\rho_n(\K Q_{(n)})=\{0\}$ in
$\Hom(X_n,X_n)$, for all $n\in\N$. Define $\fF_\nilKQ(f)$
and $\ep_{g,f}$ as in Example \ref{aa10ex1}, but restricting
to $\fF_\nilKQ(U)$, $\fF_\nilKQ(V)$, $\fF_\nilKQ(W)$.
Then $\fF_\nilKQ:\Sch_\K\ra\excat$ is a {\it
contravariant\/ $2$-functor}.
\label{aa10ex2}
\end{ex}

The point here is that $(X,\rho)\in\modKQ$ is nilpotent if
$\rho(\K Q_{(n)})=0$ for some $n\in\N$. But in a {\it family}
of nilpotent representations $(X_u,\rho_u)$ parametrized by $u$
in a base scheme $U$, this number $n$ could vary with $u$, and
might be {\it unbounded\/} on $U$. Thus it is not enough to
define $\fF_\nilKQ(U)$ as the subcategory of $(X,\rho)$ in
$\fF_\modKQ(U)$ with $\rho(\K Q_{(n)})=\{0\}$ for some $n\in\N$.
Instead, we cover $U$ by open sets $U_n$ with $\rho(\K Q_{(n)})=
\{0\}$ over $U_n$ for~$n=0,1,\ldots$.

The extension of Example \ref{aa10ex1} to quivers with
relations $(Q,I)$ is trivial.

\begin{ex} Let $(Q,I)$ be a {\it quiver with relations}. Take
$\A=\modKQI$, and define $K(\modKQI)\cong\Z^{Q_0}$ as in
Example \ref{aa10ex1}. For $U\in\Sch_\K$ define $\fF_\modKQI(U)$
to be the full exact subcategory of $\fF_\modKQ(U)$ with objects
$(X,\rho)$ such that $\rho(I)=\{0\}$ in $\Hom(X,X)$. Define
$\fF_\modKQI(f)$ and $\ep_{g,f}$ as in Example \ref{aa10ex1},
but restricting to $\fF_\modKQI(U),\ldots,\fF_\modKQI(W)$. Then
$\fF_\modKQI:\Sch_\K\ra\excat$ is a {\it contravariant\/ $2$-functor}.
\label{aa10ex3}
\end{ex}

We can also combine Examples \ref{aa10ex2} and~\ref{aa10ex3}.

\begin{ex} Let $(Q,I)$ be a {\it quiver with relations}.
Take $\A=\nilKQI$. Then $K_0(\nilKQI)\cong\Z^{Q_0}$. Set
$K(\nilKQI)=K_0(\nilKQI)$. For $U\in\Sch_\K$ define
$\fF_\nilKQI(U)$ to be the intersection of $\fF_\nilKQ(U)$ and
$\fF_\modKQI(U)$ in $\fF_\modKQ(U)$, which is a full exact
subcategory of $\fF_\modKQ(U)$. Let $\fF_\nilKQI(f)$ and
$\ep_{g,f}$ be as in Example \ref{aa10ex1}, but restricting to
$\fF_\nilKQI(U),\ldots,\fF_\nilKQI(W)$. Then $\fF_\nilKQI:
\Sch_\K\ra\excat$ is a {\it contravariant\/ $2$-functor}.
\label{aa10ex4}
\end{ex}

Here is the generalization to finite-dimensional algebras~$A$.

\begin{ex} Let $A$ be a {\it finite-dimensional\/ $\K$-algebra},
with {\it Ext-quiver} $Q$. Take $\A=\modA$, so that $K_0(\modA)
\cong\Z^{Q_0}$. Set $K(\modA)=K_0(\modA)$. Define a {\it
contravariant\/ $2$-functor} $\fF_\modA:\Sch_\K\ra\excat$ as
in Example \ref{aa10ex1}, replacing $\K Q$ by $A$ throughout.
\label{aa10ex5}
\end{ex}

\subsection{Verifying Assumption \ref{aa7ass}}
\label{aa103}

We show the examples of \S\ref{aa102} satisfy Assumption \ref{aa7ass},
following~\S\ref{aa92}.

\begin{thm} Examples \ref{aa10ex1}--\ref{aa10ex5} satisfy
Assumption~\ref{aa7ass}.
\label{aa10thm1}
\end{thm}

\begin{proof} In each example, if $(X,\rho)\in\A$ then
$[(X,\rho)]\in K(\A)$ corresponds to $\bdim(X,\rho)\in\Z^{Q_0}$,
so $[(X,\rho)]=0$ implies $\bdim(X,\rho)=0$ and hence $\dim X=0$,
so that $X=\{0\}$ and $(X,\rho)\cong 0$ in $\A$. This proves the
condition on $K(\A)$ in Assumption \ref{aa7ass}. For the rest of
the proof, we do Example \ref{aa10ex1} first.

We must prove Definition \ref{aa2def7}(i)--(iii) for $\fF_\modKQ$.
Let $\{f_i:U_i\!\ra\!V\}_{i\in I}$ be an open cover of $V$ in the
site $\Sch_\K$. For (i), let $(X,\rho),(Y,\si)\!\in\!\Obj(\fF_\modKQ(V))$
and $\phi_i:\fF_\modKQ(f_i)(X,\rho)\!\ra\!\fF_\modKQ(f_i)(Y,\si)$ for
$i\!\in\!I$ satisfy \eq{aa2eq4}. Applying \cite[Cor.~VIII.1.2]{Grot3}
to the family of sheaf morphisms $\phi_i:f_i^*(X)\ra f_i^*(Y)$ gives
a unique morphism $\eta:X\ra Y$ with $f_i^*(\eta)\!=\!\phi_i$. Let
$\ga\!\in\!\K Q$. Then
\e
\begin{gathered}
f_i^*\bigl(\eta\ci\rho(\ga)\bigr)=
f_i^*(\eta)\ci f_i^*\bigl(\rho(\ga)\bigr)=
\phi_i\ci f_i^*(\rho)(\ga)=\\
f_i^*(\si)(\ga)\ci\phi_i=
f_i^*\bigl(\si(\ga)\bigr)\ci f_i^*(\eta)=
f_i^*\bigl(\si(\ga)\ci\eta\bigr),
\end{gathered}
\label{aa10eq2}
\e
since $\phi_i:\bigl(f_i^*(X),f_i^*(\rho)\bigr)\ra
\bigl(f_i^*(Y),f_i^*(\si)\bigr)$ is a morphism in
$\fF_\modKQ(U_i)$. Using \eq{aa10eq2}, uniqueness
in \cite[Cor.~VIII.1.2]{Grot3} implies that
$\eta\ci\rho(\ga)=\si(\ga)\ci\eta$. As this holds for
all $\ga\in\K Q$, $\eta:(X,\rho)\ra(Y,\si)$ lies in
$\Mor(\fF_\modKQ(V))$. This proves Definition \ref{aa2def7}(i).
Part (ii) follows from~\cite[Cor.~VIII.1.2]{Grot3}.

For (iii), let $(X_i,\rho_i)\in\Obj(\fF_\modKQ(U_i))$ and
$\phi_{ij}:\fF_\modKQ(f_{ij,j})(X_j,\rho_j)\ra\fF_\modKQ
(f_{ij,i})(X_i,\rho_i)$ for $i,j\in I$ satisfy \eq{aa2eq5}.
Then \cite[Cor.~VIII.1.3]{Grot3} constructs $X\in\qcoh(V)$
and isomorphisms $\phi_i:f_i^*(X)\ra X_i$ satisfying
\eq{aa2eq6}, and \cite[Prop.~VIII.1.10]{Grot3} implies $X$
is locally free of finite rank. Using \cite[Cor.~VIII.1.2]{Grot3}
we construct $\rho$ from the $\rho_i$ such that $(X,\rho)\in\Obj
\bigl(\fF_\modKQ(V)\bigr)$ satisfies Definition \ref{aa2def7}(iii).
Thus $\fF_\modKQ$ is a {\it stack in exact categories}.

As a locally free sheaf of finite rank on the point $\Spec\K$
is just a finite-dimensional $\K$-vector space, $\fF_\modKQ
(\Spec\K)=\modKQ$, and Assumption \ref{aa7ass}(i) holds. Part
(ii) follows from \cite[Prop.~IV.2.2.7]{Grot2} as in
 Theorem~\ref{aa9thm1}.

An object $(X,\rho)\in\fF_\modKQ(U)$ is equivalent to
vector bundles $X_v$ on $U$ for $v\in Q_0$ and morphisms
$\phi_a:X_{b(a)}\ra X_{e(a)}$ for $a\in Q_1$. For each
$u\in\Hom(\Spec\K,U)$ the class $\bigl[\fF_\modKQ(u)
(X,\rho)\bigr]$ in $K(\modKQ)=\Z^{Q_0}$ is the map
$Q_0\ra\Z$ taking $v\in Q_0$ to the rank of $X_v$ at $u$.
Clearly, this is a locally constant function of $u$, so
Assumption \ref{aa7ass}(iii) holds. Part (iv) can be
verified as in Theorem \ref{aa9thm1}. Thus, the data of
Example \ref{aa10ex1} satisfies Assumption \ref{aa7ass}.
The modifications for Examples \ref{aa10ex2}--\ref{aa10ex5}
are all more-or-less trivial.
\end{proof}

\subsection{Verifying Assumption \ref{aa8ass}}
\label{aa104}

Here is the analogue of Theorem \ref{aa9thm2} for
Examples \ref{aa10ex1}--\ref{aa10ex5}. Note however that
we prove $\fObj_\A^\al,\fExact_\A^{\al,\be,\ga}$ are of
{\it finite type}, not just {\it locally} so.

\begin{thm} In Examples \ref{aa10ex1}--\ref{aa10ex5}, $\fObj_\A^\al$
and\/ $\fExact_\A^{\smash{\al,\be,\ga}}$ are algebraic $\K$-stacks
of finite type for all\/~$\al,\be,\ga$.
\label{aa10thm2}
\end{thm}

\begin{proof} We begin with Example \ref{aa10ex1}, so that
$\A=\modKQ$. Fix $\al\in\N^{Q_0}$. For each $v\in Q_0$ choose
a $\K$-vector space $A_v$ with $\dim A_v=\al(v)$. Define
$L=\prod_{a\in Q_1}A_{b(a)}^*\ot A_{e(a)}$. Then $L$ is
a finite-dimensional $\K$-vector space, and thus an {\it
affine $\K$-scheme}. Write $\GL(A_v)$ for the group of
automorphisms of $A_v$. Then $G=\prod_{v\in Q_0}\GL(A_v)$
acts naturally on $L=\prod_{a\in Q_1}A_{b(a)}^*\ot A_{e(a)}$,
so we may form the {\it quotient stack\/}~$[L/G]$.

Let $U\in\Sch_\K$ and $(X,\rho)\in\fObj_\modKQ^\al(U)$. Then $X$
decomposes naturally as $\bigop_{v\in Q_0}X_v$, for $X_v$ a locally
free sheaf over $U$ of rank $\al(v)$. As locally free sheaves are
locally trivializable we may choose an open cover $\{f_i:U_i
\ra U\}_{i\in I}$ of $U$ in the site $\Sch_\K$ such that $f_i^*(X_v)$
is trivial of rank $\al(v)$ on $U_i$ for all $i\in I$, $v\in Q_0$.
So we may choose isomorphisms $f_i^*(X_v)\cong A_v\t U_i$, unique
up to the action of $\Hom\bigl(U,\GL(A_v)\bigr)$. Using this and
the definitions of $\fObj^\al_\modKQ$ and $[L/G]$, one can construct
a 1-isomorphism $\fObj^\al_\modKQ\cong[L/G]$. Thus $\fObj_\modKQ^\al$
is an {\it algebraic $\K$-stack of finite type}, as $[L/G]$ is.

Now fix $\al,\ga$ and $\be=\al+\ga$ in $\N^{Q_0}$. For all $v\in Q_0$,
choose $\K$-vector spaces $A_v,B_v,C_v$ with $\dim A_v=\al(v)$,
$\dim B_v=\be(v)$, and $\dim C_v=\ga(v)$. Define
\e
\begin{gathered}
M=\Bigl\{\bigl(\Pi_ax_a,\Pi_ay_a,\Pi_az_a,\Pi_vp_v,\Pi_vq_v\bigr)\in
\prod_{a\in Q_1}(A_{b(a)}^*\!\ot\!A_{e(a)})\t\\
\prod_{a\in Q_1}\!\!(B_{b(a)}^*\!\ot\!B_{e(a)})\!\t\!\!\!
\prod_{a\in Q_1}\!\!(C_{b(a)}^*\!\ot\!C_{e(a)})\!\t\!\!\!
\prod_{v\in Q_0}\!\!(A_v^*\!\ot\!B_v)\!\t\!\!\!
\prod_{v\in Q_0}\!\!(B_v^*\!\ot\!C_v):\\
y_a\ci p_{b(a)}=p_{e(a)}\ci x_a
\;\>\text{and}\;\>
z_a\ci q_{b(a)}=q_{e(a)}\ci y_a
\;\>\text{for all $a\in Q_1$,}\\
\text{and}\;\> 0\ra A_v\,{\buildrel p_v\over\longra}\,
B_v\,{\buildrel q_v\over\longra}\,C_v\ra 0
\;\>\text{is exact for all $v\in Q_0$}\Bigr\}.
\end{gathered}
\label{aa10eq3}
\e

This defines $M$ as a subset of a finite-dimensional $\K$-vector
space, $N$ say. The third line of \eq{aa10eq3} is finitely many
quadratic equations in $N$. Exactness in the fourth line is
equivalent to $q_v\ci p_v=0$, more quadratic equations, together
with injectivity of $p_v$ and surjectivity of $q_v$, which are
open conditions. Thus, $M$ is a Zariski open subset of the
zeroes of finitely many polynomials in $N$, and is a {\it
quasiaffine $\K$-scheme}. Define $H=\prod_{v\in Q_0}(\GL(A_v)
\t\GL(B_v)\t\GL(C_v))$, as an {\it algebraic $\K$-group}. Then
$H$ has an obvious action on $M$, so we can form the {\it
quotient stack\/} $\bigl[M/H]$. A similar proof to the
$\fObj_\modKQ^\al$ case gives a 1-isomorphism
$\fExact_\modKQ^{\smash{\al,\be,\ga}}\cong[M/H]$. Thus
$\fExact_\modKQ^{\smash{\al,\be,\ga}}$ is an {\it algebraic
$\K$-stack of finite type}, as $[M/H]$ is. This proves the
theorem for Example~\ref{aa10ex1}.

Next we do Example \ref{aa10ex3}, so let $(Q,I)$ be a {\it quiver
with relations}. Let $A_v,L$ and $G$ be as above, and set $A=
\bigop_{v\in Q_0}A_v$. For each $\Pi_ax_a$ in $L$, define
$\rho_{\Pi_ax_a}:\K Q\ra\End(A)$ to be the unique algebra
homomorphism such that $\rho(v)=1$ on $A_v$ and 0 on $A_w$ for
$v\ne w\in Q_0$, and $\rho(a)=x_a:A_{b(a)}\ra A_{e(a)}$ on
$A_{b(a)}$, $\rho(a)=0$ on $A_v$ for $b(a)\ne v\in Q_0$, and
all $a\in Q_1$. Define
\begin{equation*}
L_I=\bigl\{\Pi_ax_a\in\ts\prod_{a\in Q_1}(A_{b(a)}^*\!\ot\!A_{e(a)}):
\rho_{\Pi_ax_a}(i)=0
\;\>\text{for all $i\in I$}\bigr\}.
\end{equation*}

Each $i\in I$ is a finite linear combination of basis
elements \eq{aa10eq1} of $\K Q$, so $\rho_{\Pi_ax_a}$
is a finite linear combination of the corresponding
$x_{a_k}\!\ci\!\cdots\!\ci\!x_{a_1}$ in $\End(A)$. Thus for
fixed $i\!\in\!I$, the map $\Pi_ax_a\!\mapsto\!\rho_{\Pi_ax_a}(i)$
is a polynomial on $L$ with values in $\End(A)$. Hence
$L_I$ is the zeroes of a collection of polynomials on $L$,
and so is a $G$-invariant affine $\K$-scheme. Modifying
the proof above shows $\fObj_\modKQI^\al$ is 1-isomorphic to
$[L_I/G]$, and so is {\it algebraic} and of {\it finite type}.

Similarly, for $\fExact_\modKQI^{\smash{\al,\be,\ga}}$ we
replace $M$ in \eq{aa10eq3} by $M_I$, where we add extra
conditions $\rho_{\Pi_ax_a}(i)=0$ in $\End(\bigop_{v\in Q_0}A_v)$,
$\rho_{\Pi_ay_a}(i)=0$ in $\End(\bigop_{v\in Q_0}B_v)$ and
$\rho_{\Pi_az_a}(i)=0$ in $\End(\bigop_{v\in Q_0}C_v)$ on
$\bigl(\Pi_ax_a,\ldots,\Pi_vq_v\bigr)$ for all $i\in I$.
Then $M_I$ is a quasiaffine $\K$-scheme invariant under $H$,
and $\fExact_\modKQI^{\smash{\al,\be,\ga}}$ is 1-isomorphic
to $[M_I/H]$, so it is {\it algebraic} and of {\it finite type}.

For Examples \ref{aa10ex2} and \ref{aa10ex4}, suppose $(X,\rho)\in
\Obj(\nilKQ)$. Then the vector subspaces $X_k=\rho(\K Q_{(k)})X$
of $X$ must decrease strictly in dimension until they become zero.
Hence $\rho(\K Q_{(m)})X=\{0\}$ for $m=\dim X$. Now let $U\in\Sch_\K$
and $(X,\rho)\in\fObj^\al_\modKQ(U)$. The same proof shows
$(X,\rho)\in\fObj^\al_\nilKQ(U)$ if and only if $\rho(\K Q_{(m)})=0$
in $\End(X)$ for $m=\sum_{v\in Q_0}\al(v)$. By a similar argument for
$\fExact_\modKQ^{\smash{\al,\be,\ga}}$, we see that
\begin{align*}
\fObj^\al_\nilKQ&=\fObj^\al_{\modKQ/\K Q_{(m)}},&\;
\fObj^\al_\nilKQI&=\fObj^\al_{\modKQ/(I+\K Q_{(m)})},\\
\fExact_\nilKQ^{\al,\be,\ga}&=\fExact_{\modKQ/\K Q_{(n)}
}^{\al,\be,\ga},&\;
\fExact_\nilKQI^{\al,\be,\ga}&=\fExact_{\modKQ/(I+\K Q_{(n)})
}^{\al,\be,\ga},
\end{align*}
for $m=\sum_{v\in Q_0}\al(v)$ and $n=\sum_{v\in Q_0}\be(v)$. Thus
the theorem for Examples \ref{aa10ex2} and \ref{aa10ex4} follows
from the Example \ref{aa10ex3} case, with $\K Q_{(n)},\K Q_{(m)},
I+\K Q_{(m)}$ or $I+\K Q_{(n)}$ in place of the ideal~$I$.

Finally, in the situation of Definition \ref{aa10def4} and
Example \ref{aa10ex5}, the equivalence of categories between
$\modA$ and $\modKQI$ extends easily to give an equivalence
of 2-functors $\fF_{\smash{\modA}}\ra\fF_{\smash{\modKQI}}$, and
hence 1-isomorphisms $\fObj_\modA^\al\ra\fObj_\modKQI^\al$ and
$\fExact_\modA^{\al,\be,\ga}\ra\fExact_\modKQI^{\al,\be,\ga}$.
Thus the theorem for Example \ref{aa10ex5} again follows from the
Example \ref{aa10ex3} case.
\end{proof}

In each of Examples \ref{aa10ex1}--\ref{aa10ex5}, Theorem
\ref{aa10thm3} shows $\fObj^\al_\A,\fExact^{\al,\be,\ga}_\A$ are
{\it algebraic} $\K$-stacks of {\it finite type}, and hence
{\it locally} of finite type. Since $\fExact^{\al,\be,\ga}_\A$
is of finite type, it follows immediately that the 1-morphisms
\eq{aa8eq1} are of {\it finite type}. Thus we prove:

\begin{thm} Examples \ref{aa10ex1}--\ref{aa10ex5} satisfy
Assumption~\ref{aa8ass}.
\label{aa10thm3}
\end{thm}

We can also strengthen Theorem \ref{aa8thm1} for the examples
of~\S\ref{aa102}.

\begin{thm} Let\/ $\A,K(\A),\fF_\A$ be as in any of
Examples \ref{aa10ex1}--\ref{aa10ex5}. Then $\fM(I,\pr,\ka)_\A$
is an algebraic $\K$-stack of finite type for all\/~$(I,\pr),\ka$.
\label{aa10thm4}
\end{thm}

\begin{proof} This follows as $\bs\si(I):\fM(I,\pr,\ka)_\A\ra
\fObj^{\ka(I)}_\A$ is of finite type by Theorem \ref{aa8thm2}, and
$\fObj^{\ka(I)}_\A$ is of finite type by Theorem~\ref{aa10thm2}.
\end{proof}

\medskip

\noindent{\small\sc The Mathematical Institute, 24-29 St. Giles,
Oxford, OX1 3LB, U.K.}

\noindent{\small\sc E-mail: \tt joyce@maths.ox.ac.uk}

\end{document}